\documentclass{article}
\usepackage[utf8]{inputenc}
\usepackage{amsmath}
\usepackage{amssymb}
\usepackage{graphicx}
\usepackage[english]{babel}
\usepackage[top=2cm, bottom=3cm, left=2cm, right=2cm]{geometry}
\usepackage{listings}
\usepackage{algorithm}
\usepackage{algorithmic}
\usepackage{multicol}
\usepackage{amsthm}
\usepackage{tikz}
\usepackage{tikz-network}
\usepackage[toc,page]{appendix}
\usepackage{hyperref}
\usepackage{subfig}
\usepackage{xspace}
\usepackage{todonotes}
\usetikzlibrary{shapes.geometric}
\usepackage{authblk}

\usepackage{thmtools}
\usepackage{thm-restate}

\usepackage{hyperref}

\usepackage{cleveref}

\newtheorem{theorem}{Theorem}[section]

\newtheorem{corollary}[theorem]{Corollary}

\newtheorem{definition}[theorem]{Definition}
\newtheorem{observation}[theorem]{Observation}

\newtheorem{claim}[theorem]{Claim}

\newtheorem{lemma}[theorem]{Lemma}

\newtheorem{remark}[theorem]{Remark}

\newcommand{\D}{\mathcal{D}}

\newcommand{\sk}{S}

\newcommand{\Awin}{${\mathcal A}$-win\xspace}
\newcommand{\Draw}{${\mathcal D}$raw\xspace}
\newcommand{\oA}{\mathcal A}
\newcommand{\oD}{\mathcal D}

\tikzstyle{v}=[circle,inner sep=0, minimum size =6 pt, line width = 1pt, draw=black, fill=black, text= white]
\tikzstyle{R}=[circle,inner sep=0, minimum size =7 pt, line width = 1pt, draw=red, fill = red]
\tikzstyle{B}=[circle,inner sep=0, minimum size =7 pt, line width = 1pt, draw=blue, fill = blue]
\tikzstyle{tri}=[triangle,inner sep=0, minimum size =7 pt, line width = 1pt, draw=black, fill = black]

\usetikzlibrary{shapes,arrows}
\tikzstyle{decision} = [diamond, draw, fill=yellow!20, 
    text width=10em, text centered, node distance=5cm, inner sep=0pt]
\tikzstyle{block} = [rectangle, draw, fill=yellow!20, 
    text width=10em, text centered, rounded corners, minimum height=4em]
\tikzstyle{line} = [draw, -latex']
    
\tikzstyle{sortieD} = [draw, ellipse,fill=blue!20, node distance=3cm,
    minimum height=2em]
 \tikzstyle{sortieA} = [draw, ellipse,fill=red!20, node distance=3cm,
    minimum height=2em]

\newcommand{\qedclaim}{\hfill $\diamond$ \medskip}

\newenvironment{proofclaim}{\noindent{\em Proof of the claim.}}{\qedclaim}

\title{The Maker-Maker domination game in forests\thanks{This work has been supported by the National Research Agency (ANR) projects P-GASE (ANR-21-CE48-0001-01).}}

\author[1]{Eric Duchêne}

\author[2]{Arthur Dumas}

\author[1]{Nacim Oijid}

\author[1]{Aline Parreau}

\author[3]{Eric Rémila}

\affil[1]{Univ. Lyon, Universit\'e Lyon 1, LIRIS UMR CNRS 5205, F-69621, Lyon, France.}

\affil[2]{Univ. Rennes, ENS Rennes, Bruz, France.}

\affil[3]{Univ. Lyon, UJM Saint‐Etienne, GATE UMR CNRS 5824, F-42023, Saint-Etienne, France.}

\date{\today}

\begin{document}

\maketitle

\begin{abstract}
    We study the Maker-Maker version of the domination game introduced in 2018 by Duchêne {\em et al}. 
    Given a graph, two players alternately claim vertices. The first player to claim a dominating set of the graph wins. As the Maker-Breaker version, this game is {\sf PSPACE}-complete on split and bipartite graphs. Our main result is a linear time algorithm to solve this game in forests. We also give a characterization of the cycles where the first player has a winning strategy.

\end{abstract}

%\todo{Eric (R), tu peux ajouter ton affiliation ?}

\section{Introduction}
Positional games have been introduced successively by Hales and Jewett in~\cite{Hales1963} and by Erd\H{o}s and Selfridge in~\cite{erdos}, and then widely studied in the literature (see the two books \cite{beck,positionalgames} for an overview). They are played on an hypergraph of vertex set $X$, with a finite set $\mathcal{F}\subseteq 2^X$ of hyperedges. The set $X$ is often called the {\it board} of the game, and an element of $\mathcal{F}$ a {\it winning set}. The game involves two players that alternately claim a previously unclaimed vertex of $X$. The winner of the game is defined according to a convention. The original one is called {\it Maker-Maker} (or strong convention), where both players have the same objective, i.e. filling a whole winning set with their own claimed vertices. Such games may end in a draw if each winning set contains one vertex claimed by each player. When considering positional games, the main issue consists in determining which player has a winning strategy. In particular, it is well-known that in the Maker-Maker convention, the second player has no winning strategy. Thus, resolving a Maker-Maker game consists in determining whether the first player has a winning strategy or whether it ends in a draw. As mentioned in~\cite{positionalgames}, despite this result, this convention has not been widely studied in the literature. The main reason is due to the hardness to tackle it, as the first player tries at the same time to fill a winning set while considering all the threats of his opponent. As a consequence, Maker-Maker instances often satisfy the so-called {\em extra set paradox}, which claims that adding new winning sets in the hypergraph is not necessarily better for the first player.

For all these reasons, the convention that has been mainly studied in the literature is the {\it Maker-Breaker convention}. In this convention, both players, called Maker and Breaker, have opposite objectives: Maker aims at filling a winning set while Breaker prevents her to do so. Constructive results are generally more affordable in this convention. In addition, the extra set paradox does not exist when playing in the Maker-Breaker convention: adding new winning sets is always better for Maker.\\

In the literature, there are many graph optimization problems that have been turned into positional games. One can cite for example the clique game, the connectivity game, the $H$-game~\cite{beck}, or the domination game~\cite{makerbreaker}. If a large part of the studies is devoted to the case where they are played on complete or random graphs, there is a more recent approach that consists in playing such games on a general graph. It generally yields to algorithmic results, both in terms of hardness proofs or the construction of polynomial time algorithms to compute the winner of the game~\cite{bonnetComplexity,papierAvecMilos}. One can also find results about the parameterized complexity of such games~\cite{bonnet1}. In addition, general algorithmic results about Maker-Breaker games played on $k$-uniforms hypergraphs (i.e., all the winning sets are of size $k$) have been given by Rahman and Watson~\cite{Rahman20206UniformMG} and Galliot~\cite{Florian}. Indeed, they respectively proved that determining the winner of a Maker-Breaker positional game is PSPACE-complete when $k=6$ and polynomial when $k=3$.\\

When switching to the Maker-Maker convention, by putting together the above result of Rahman and Watson with an argument of Byskov~\cite{byskov2004}, it has been derived that Maker-Maker games are PSPACE-complete on $7$-uniforms hypergraphs. Beyond this result, connections between the two types of conventions are not well established. Indeed, there are very few algorithmic results, even for particular positional games derived from optimization problems. The objective of this paper is to investigate the Maker-Maker domination game, by highlighting the similarities, the implications, and the differences with the results known in the Maker-Breaker convention. It is known that the Maker-Breaker domination game is PSPACE-complete even for bipartite graphs and split graphs, and polynomial for cographs and forests\cite{makerbreaker}. More precisely, it is shown that Maker has a winning strategy playing second on a forest if and only if it admits a perfect matching. When switching to the Maker-Maker convention, we will see that some of these complexity results still hold. Yet, in the case of forests, the polynomial characterization of the winning positions is far more complex than simply finding a perfect matching. The major result of the current paper consists in determining this characterization. To the best of our knowledge, this is the first time in the literature that a non-trivial algorithm is given to determine the winner of a game played according to the Maker-Maker convention. \\

The paper is organized as follows: in Section 2, we present the main definitions that will be useful in the Maker-Maker domination game. In Section 3, we give results about the Maker-Maker convention that are derived from the Maker-Breaker one, including the PSPACE-hardness result. Then, we fully solve the case of paths and cycles in Section 4. Sections 5 and 6 are devoted to the resolution of the Maker-Maker domination game on forests in linear time. As the proof is rather complex, we decided to split it into two parts, where Section 5 corresponds to the overview of the proof with the presentation of the algorithm, and Section 6 to the proof of the most technical elements.

\section{Preliminaries}

\subsection{Standard definitions of graph theory}

In this paper, we will only consider finite, undirected and simple graphs. A graph $G$ is defined by a couple $(V,E)$ where $V$ denotes the set of vertices and $E$ the edges of the graph. The  {\it (closed) neighborhood} of a vertex $x\in V$, denoted by $N[x]$, is the subset of vertices containing $x$ and all the vertices that are adjacent to $x$. A vertex $x$ is {\em universal} if $N[x]=V$. The {\it degree} of a vertex is the number of vertices adjacent to it. A {\it leaf} is a vertex of degree 1. If $x$ is a vertex, $G\setminus \{x\}$ denotes the graph on the vertex set $V\setminus \{x\}$ with all the edges of $E$ that are not incident to $x$. 

Let $X$ be a subset of vertices. $X$ is a {\em independent set} if there are no adjacent vertices in $X$. $X$ is a {\em cutset} if $G\setminus X$ is disconnected. The subgraph of $G$ {\em induced by $X$}, denoted by $G[X]$, is the graph with vertex set $X$ and edge set all the edges of $E$ whose both extremities are in $X$.
A {\em matching} $M$ is a subset of edges that are two by two not incident. If $|M|=|V|/2$ (i.e. if all the vertices appear in some edge of $M$) then $M$ is {\em perfect}.

Let $G_1=(V_1,E_1)$ and $G_2=(V_2,E_2)$ be two graphs on disjoint vertex sets. The {\em union} of $G_1$ and $G_2$ is the graph $G_1\cup G_2=(V_1\cup V_2,E_1\cup E_2)$. The {\em join} of $G_1$ and $G_2$ is the graph on vertex set $V_1\cup V_2$ with edge set $E_1\cup E_2 \cup E_{\times}$ where $E_{\times}=V_1\times V_2$.

A {\em cograph} is a graph that is either a single vertex or the union of two cographs or the join of two cographs. A {\em bipartite} graph is a graph whose vertex set can be partitioned into two independent sets. A {\em split} graph is a graph whose vertex set can be partitioned into two sets $K$ and $I$ where $K$ induces a complete graph and $I$ an independent set. A {\em path} is a graph whose vertex set is $\{v_1,....,v_n\}$ and $v_iv_j$ is an edge if $|i-j|=1$. A {\em cycle} is a path with the additional edge $v_1v_n$. A {\em forest} is a graph without any cycle. A {\em tree} is a connected forest. 

A vertex $x$ {\em dominates} a vertex $y$ if $y\in N[x]$. A subset of vertices $S$ dominates a vertex $y$ if there exists $x\in S$ that dominates $y$. A {\it dominating set} $S$ of $G$ is a subset of vertices that dominates all the vertices of the graph. The smallest size of a dominating set of $G$ is denoted by~$\gamma(G)$.

\subsection{The Maker-Maker domination game}

The Maker-Maker domination game is played on a graph $G=(V,E)$. Two players, Alice and Bob, alternately claim an unclaimed vertex of the graph, with Alice playing first. The game ends when the vertices claimed by one of players form a dominating set (in which case the corresponding player wins) or when all the vertices have been claimed and none of the players managed to claim a dominating set (in which case the game is a draw).

As a $2$-player finite perfect information game, if both players play optimally, one of the players has a winning strategy or the game ends a draw. Furthermore, since this game is a Maker-Maker positional game, it is well known that the second player, Bob, does not have a winning strategy (see for example \cite{positionalgames}). As a consequence, there are only two possible {\it outcomes} for the Maker-Maker domination game played on $G$: either Alice has a winning strategy, which will be denoted by $o(G)=\oA$, or both players can ensure a draw, which will be denoted by $o(G)=\oD$. The problem of deciding, given a graph $G$, if $o(G)=\oA$ or $o(G)=\oD$ is named the {\sc Maker-Maker Domination Game} problem. 

\paragraph{Positions.} A {\it position} $P$ of the game is a triple $(G,V_A, V_B)$ where $G = (V,E)$ is a graph and $V_A,V_B \subseteq V$ are two subsets of vertices such that $V_A\cap V_B=\emptyset$. The vertices in $V_A$ (respectively $V_B$) correspond to vertices claimed by Alice (resp. Bob). A vertex not in $V_A$ nor $V_B$ is {\it unclaimed}. 
If $X$ is a set of vertices, the {\em subposition induced by $X$} is the position $(G[X],V_A\cap X,V_B\cap X)$. 
If \{x, y\} is a pair of unclaimed vertices of $P$. The position $P_{x, y }$ is defined by $P_{x, y } = (G, V_A \cup \{x \}, V_B   \cup \{y\})$. 
A set of unclaimed vertices $S$ is a {\em winning set for Player t} if $S\cup V_t$ is a dominating set of $G$.

A {\it pointed position} is a position where the next player is specified, it will be denoted by a couple $(P,t)$ where $P$ is a position and $t \in \{A,B\}$ is the next player ($A$ stands for Alice and $B$ for Bob). 

The {\em outcome $o(P,t)$} of a pointed position $(P,t)$ is defined for positions $P$ satisfying the property that both $V_A$ and $V_B$ do not dominate $G$ simultaneously. We have $o(P,t)=\oA$,  if, the next player being  $t$, considering that Alice  has already claimed all the vertices $V_A$ and  Bob has claimed vertices $V_B$, there is a winning strategy for Alice, i.e. a strategy  which ensures Alice to dominate the graph $G$ before Bob. We will say in this case that the position is \Awin. 
Otherwise, we say that the position is \Draw and we denote it by $o(P,t)=\oD$. In this latter case, note that it also covers the case where Bob first dominates the graph, as the pointed position $P$ may not be balanced depending on the sets $V_A$ and $V_B$. Such cases will be considered in the upcoming proofs when considering the different possible sequences of moves (but their outcome will still be considered as a draw, as Alice will avoid them in her optimal sequence of moves).

Note that the starting position on the graph $G$ is the pointed position $((G,\emptyset, \emptyset),A)$ and we have $o(G)=o((G,\emptyset,\emptyset),A)$.

\paragraph{Ordering positions.} Two positions $P$ and $P'$ are said \emph{equivalent} if for any $t\in\{A,B\}$, $o(P,t)=o(P',t)$.
Two pointed positions $(P,t)$ and $(P',t)$ with the same first player are said \emph{equivalent} if $o(P,t)=o(P',t)$.

In addition, and by analogy with combinatorial game theory, it is standard to order outcomes, stating that $\oA > \oD$. With this convention, a pointed position $(P,t)$ is {\em better} for Alice (respectively Bob) than a pointed position $(P',t')$ if $o(P,t) \geq o(P',t')$ (respectively $o(P,t) \leq o(P',t'))$. 

As an illustration of these definitions, the two following observations ensure that:
\begin{itemize}
    \item it is always better for any player to start.
    \item if two positions have exactly the same winning sets, there are equivalent.
\end{itemize}

\begin{observation}\label{obs:betterstarts}
For any position $P$,  we have $o(P, A) \geq o(P, B)$.
\end{observation}
The proof of the above result derives from a standard stealing strategy.

\begin{observation}\label{obs:sameWS}
Let $P=(G,V_A,V_B)$ and $P'=(G',V'_A,V'_B)$ be two positions such that there exists a bijection $f:V(G)\setminus V_A\cup V_B \to V(G')\setminus V'_A\cup V'_B$ between the sets of unclaimed vertices that satisfies that $S$ is a winning set for Alice (respectively Bob) if and only if $f(S)$ is a winning set for Alice (respectively Bob). Then $P$ and $P'$ are equivalent.
\end{observation}

We  say that a vertex $x$  is \emph{forced} for Alice (resp.  Bob) in a position $(P,A)$ (resp. $(P,B)$, if whenever Alice (resp.  Bob) claims a vertex $y \neq x$ in $P$, the resulting position is \Draw  (resp.  \Awin).
Sometimes, a move is always better than another. This is for example the case when the neighborhood of one vertex contained another neighborhood. The following lemma gives a general formal framework which allows us to eliminate uninteresting moves. Note that this lemma or its derivatives are often used when studying positional games, sometimes not explicitly. For the sake of completeness, we give a proof here.

\begin{lemma}\label{lem:sommetdomine}
    Let $P = (G, V_A, V_B)$ be a position. Let $x,y$ be two unclaimed vertices. Assume that for both $t \in \{A,B\}$, $N[x]\setminus N[V_t]\subseteq N[y]\setminus N[V_t]$, where  $N[V_t] = \underset{v \in V_t}{\cup} N[v]$ .
Then there exists an optimal strategy in which $y$ is claimed before $x$.
\end{lemma}

Roughly speaking, the condition says that $y$ dominates a superset of the vertices dominated by $x$ for both players. The conclusion says that if there exists an optimal strategy (for any player) starting by claiming $x$, one can assume that there also exists an optimal strategy starting by claiming $y$.

\begin{proof}
    Let $P = (G, V_A, V_B)$ be a position. Suppose that Alice has a winning strategy $\mathcal{S}$ on $P$ in which $x$ is claimed before $y$. Up to consider more moves, consider a position $P'$ where the next vertex that will be claimed is $x$ or $y$. If it is $y$, there is nothing to do, so suppose it is $x$.
    
    If $x$ is claimed by Bob, by definition any strategy is loosing for Bob, thus they are all equivalent, and he can claim $y$ first instead without changing the outcome of the game. Therefore, we can suppose that Alice is going to claim $x$. Consider the following strategy $\mathcal{S}'$ for Alice.
    \begin{itemize}
        \item Instead of claiming $x$, claim $y$.
        \item While Bob claims a vertex $v \neq x$, she claims the vertex $w$ she would have claimed according to $\mathcal{S}$ if Bob has claimed $v$ when she had claimed $x$.
        \item If Bob claims $x$, she claims the vertex $w$ she would have claimed according to $\mathcal{S}$ if Bob had claimed $y$.
    \end{itemize}

As $\mathcal{S}$ was a winning strategy for Alice, at a certain moment of the game, she would obtain a set of vertices $S_A\supset V_A$ that is a dominating set of $G$, while $S_B\supset V_B$, the vertices claimed by Bob, is not.
If both $x$ and $y$ are in $S_A$, then $\mathcal{S'}$ will make both Alice and Bob claim exactly the same vertices, so Alice will win. Otherwise, by denoting by $S_A'$ and $S_B'$ the set of vertices claimed by Alice and Bob respectively after $\mathcal{S'}$, we have $S_A' = S_A \setminus \{x\} \cup \{y\}$ and $S_B' = S_B \setminus \{y\} \cup \{x\}$. 
We prove that $S_A'$ is still a dominating set
and $S_B'$ is not, which proves that  $\mathcal{S'}$ is a winning strategy for Alice. Indeed, for $S'_A$, if $u$ is dominated only by $x$ in $S_A$, then it means that $u\notin N[V_A]$ and thus $u\in N[x]\setminus N[V_A]$. Then by definition, $u\in N[y]\setminus N[V_A]$ and thus is dominated by $y$.
For $S_B'$, assume that it is a dominating set. Since $S_B$ is not dominating, it means that there exists $u\in N[x]$ not dominated in $S_B$. Thus, $u\notin N[y]$ but we also have $u\notin N[V_B]$ as $V_B \subset S_B$. Hence, $u$ is in $N[x]\setminus N[V_B]$ but not in $N[y]\setminus N[V_B]$, a contradiction.

The proof that if Bob has a drawing strategy is similar, by adding the fact that if Bob does not dominate the graph, then $S_A$ is not a dominating set, and thus neither is $S_A \setminus \{y\} \cup \{x\}$. 
\end{proof}

\paragraph{Union and decomposition of positions.} Let $P=(G,V_A,V_B)$ and $P'=(G',V'_A,V'_B)$ be two positions on disjoint sets of vertices \footnote{In this paper, the union of two positions will always be done on disjoint sets of vertices. It could happen, for simplicity reasons, that two vertices in different components have the same name. They must be considered as distinct vertices.}. The {\em union} of $P$ and $P'$, denoted by $P\cup P'$, is the position $(G\cup G',V_A\cup V'_A,V_B\cup V'_B)$.
%If $Q = P \cup P'$, then we say that $P' = Q \setminus P$. 
Note that one can remove a position where both players have a dominating set. It is a simple consequence of Observation~\ref{obs:sameWS} since the winning sets are the same.

\begin{observation}\label{obs:removedomcomp}
    Let $P=(G,V_A,V_B)$ and $P'=(G',V'_A,V'_B)$  be two positions on disjoint graphs. Assume that $V'_A$ and $V'_B$ are dominating sets of $G'$. Then $P\cup P'$ and $P$ are equivalent.
\end{observation}

One cannot in general determine the outcome of a position $P\cup P'$ knowing the outcome of $P$ and $P'$. However, when both positions are \Awin when Bob starts, the union is still \Awin:

\begin{observation} \label{obs:union}
Let $P$ and $P'$ be two positions such that $o(P, B) =  o(P', B) = \oA$. Then we have $o(P \cup P', B) = \oA$. 
\end{observation}

\begin{proof}
 Alice follows her strategies as second player in both $P$ and $P'$ until she dominates one of the component, say $P$. Then she just claims in $P'$ following her strategy. She might claim several times in a row in $P'$ (if Bob goes on claiming on $P$), but by Observation~\ref{obs:betterstarts}, it is always better for Alice to play first than second. Thus, if she has to claim twice in a row in $P'$, she will be able to dominate $P'$, before Bob does.\end{proof}

When considering a position, it could be useful to decompose it into several disjoint games. To do such a decomposition, the winning sets for both players should be the same. 
This is the case when a cut set is completely claimed and dominated by both players.

\begin{lemma}\label{lem:split}
Let $P=(G,V_A,V_B)$ be a position. Assume $V(G)$ can be partitioned into three sets $V_1,V_2,X$ such that:
\begin{itemize}
    \item There are no edges between $V_1$ and $V_2$;
    \item All the vertices of $X$ have been claimed: $X\subseteq V_A\cup V_B$.
    \item The vertices in $X$ are already dominated by $V_A\cap X$ and $V_B\cap X$.
\end{itemize}
Let $P_1$ and $P_2$ be the subpositions of $P$ induced by $V_1\cup X$ and $V_2\cup X$ respectively (vertices of $P_1$ and $P_2$ are disjoint).
Then the position $P$ and the position $P_1\cup P_2$ are equivalent.
\end{lemma}

\begin{proof} We consider the trivial bijective map $f$ between the unclaimed vertices of $P$ and $P_1\cup P_2$. Let $S$ be a set of unclaimed vertices of $P$ and $t\in \{A,B\}$. 
Assume first that $S$ is a winning set of $P$ for Player $t$
Let $S_1=f(S)\cap V_1$. We prove that $S_1$ is a winning set of $P_1$ for Player $t$. Let $u$ be a vertex of $P_1$ and $u'=f^{-1}(u)$. Either $u'\in X$ and thus is dominated by a vertex of $V_t\cap X$ in $P$ and thus $u$ is still dominated in $P_1$. Or $u'\notin X$, and then $u'$ is dominated by some vertex $s$ of $S\cup V_t$ in $P$. Since $u'\notin X$ and $X$ is a cutset, $s$ should be in $V_1\cup X$ and $f(s)$ is still in $S_1\cup V_t$.
Similarly, we can prove that $S_2=f(S)\cap V_2$ is a winning set of $P_2$ for Player t, and thus $f(S)$ is a winning set of $P_1\cup P_2$ for Player t. The reverse is easier: the union of two winning sets in $P_1$ and $P_2$ is clearly a winning set of $P$.
Therefore, Observation~\ref{obs:sameWS} applies and the two positions are equivalent.
\end{proof}

\subsection{Traps}
We will frequently use the notion of {\it trap} that is defined in this section. Roughly speaking, a $A$-trap (respectively $B$-trap) is a vertex of a game position such that, if it is not claimed by Alice (resp. Bob) by the end of the game, means that Alice (resp. Bob) will never build a dominating set. Formally, traps can be defined as follows:

\begin{definition}
Let $P = (G, V_A, V_B)$ be a position of the game. A {\it $A$-trap} (respectively {\it $B$-trap}) is an unclaimed vertex $v$ such that there exists a vertex $w$ with $N[w] \cap V\backslash V_B = \{v\}$ (resp. $N[w] \backslash V_A = \{v\}$)
\end{definition}

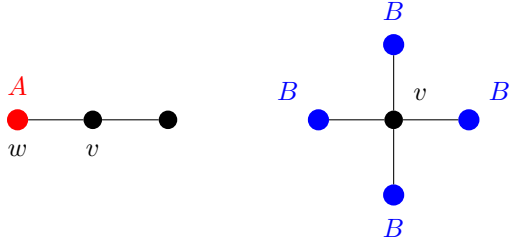
\begin{figure}
    \centering
    \begin{tikzpicture}

\draw (-1,0) node[R](v0) {} node[below=2mm] {$w$} node[above=2mm] {\color{red}$A$};
\draw (0,0) node[v](v1) {} node[below=2mm] {$v$};
\draw (1,0) node[v](v2) {} node[above=2mm] {};

\draw (4,0) node[v](s) {} node[above right=2mm] {$v$};
\draw (3,0) node[B](s1) {} node[above=2mm] {}node[above left=2mm] {\color{blue}$B$};
\draw (5,0) node[B](s2) {} node[above=2mm] {}node[above right=2mm] {\color{blue}$B$};
\draw (4,1) node[B](s3) {} node[above right=2mm] {} node[above=2mm] {\color{blue}$B$};
\draw (4,-1) node[B](s4) {} node[below right=2mm] {} node[below=2mm] {\color{blue}$B$};

\draw[] (v0) -- (v1) ;
\draw[] (v1) -- (v2) ;

\draw[] (s) -- (s1) ;
\draw[] (s) -- (s2) ;
\draw[] (s) -- (s3) ;
\draw[] (s) -- (s4) ;

    \end{tikzpicture}
    \caption{On the left, $v$ is a $B$-trap and $w$ may be isolated; on the right $v$ is a $A$-trap}
    \label{fig:trap}
\end{figure}

In this definition, $v$ corresponds to the vertex that must be claimed, and $w$ to the vertex that will not be dominated if $v$ is not claimed. Figure~\ref{fig:trap} illustrates the notion of traps.

The next lemma shows that if there is a $A$-trap in a position, one can consider that the next player will claim it immediately.

\begin{lemma}\label{lemma:trap}
Let $P$ be a position of the game and $v$ be a $A$-trap of $P$. Claiming $v$ is an optimal move for both players. Moreover, $o(P,B) = \oD$. 
\end{lemma}

\begin{proof}
Let $w$ such that $N[w] \cap V\backslash V_B = \{v\}$.
Suppose it is Bob's turn. By claiming $v$, Bob isolates the vertex $w$ as now $N[w]\subseteq V_B$. Therefore, Alice cannot dominate $w$ any longer and the outcome is $\oD$.

If it is Alice's turn, if she does not claim $v$, Bob claims it and once again isolates $w$. Therefore, $v$ is forced for Alice. Thus, $v$ is an optimal move. 
\end{proof}

\begin{corollary}\label{cor:2traps}
Let $P$ be a position of the game. If there exist two $A$-traps $v_1$ and $v_2$ of $P$ such that $v_1 \neq v_2$, then $o(P,A)=o(P,B)= \oD$.
\end{corollary}
\begin{proof}
Since there are two distinct $A$-traps, the two vertices $w_1$ and $w_2$ that might be isolated are necessarily distinct (otherwise $v_1$ and $v_2$ would both be in their neighborhood, contradicting the definition of a trap). Hence, even if Bob is not the first player, he will be able to claim in one of the two traps and hence isolate either $w_1$ or $w_2$.
\end{proof}

The next lemma proposes another example of a forced move for Alice when there exists an unclaimed $P_5$ in a position as a subgraph.

\begin{lemma}\label{lem:P5}
Let $(P,B)$ be a pointed position of the game with $P=(G,V_A,V_B)$. If there exists a path $G'=(v_1,v_2,v_3,v_4,v_5)$ such that $G'$ is a subgraph of $G$ with $V(G')\cap (V_A\cup V_B)=\emptyset$ and $v_2, v_3, v_4$ of degree $2$, then if Bob claims $v_3$, Alice is then forced to answer on $v_2$ or $v_4$.
\end{lemma}

\begin{proof}
If Bob claims $v_3$, then if Alice answers elsewhere than on $v_1$, $v_2$ or $v_4$, Bob claims his second move on $v_2$ and creates two $A$-traps in $v_1$ and $v_4$. Indeed, since the vertices $v_2$ and $v_3$ are of degree $2$, there remains only one way for Alice to dominate $v_2$ (i.e. by claiming $v_1$) and $v_3$ (i.e. by claiming $v_4$). By corollary~\ref{cor:2traps}, the resulting position is \Draw. If Alice claims $v_1$, then Bob claims $v_4$ and by symmetry creates two $A$-traps equivalent to the previous case. 
\end{proof}

\section{General results derived from the Maker-Breaker convention}

Whereas the Maker-Maker domination game has not been studied yet, the Maker-Breaker version has been defined in \cite{makerbreaker}. In this section, we recall some results in this convention that have consequences for the Maker-Maker convention. 

\subsection{Basics}
In Maker-Breaker games, only the winning condition differs from the Maker-Maker games. For more convenience and according to the terminology of \cite{makerbreaker}, we will call the players Dominator and Staller. Dominator wins if she claims a dominating set, otherwise Staller wins. In particular, there is no draw. In this convention, since the two players have different roles, one need to precise who the first player is. Note that both players have interest to start since one can prove that if a player has a winning strategy playing second, he also has a strategy playing first. 

When playing on the same graph, a winning strategy for Staller going second is a drawing strategy for Bob. Indeed, Bob will prevent Alice to claim a dominating set and thus Alice has no winning strategy.
The contrary is not true since Bob can threaten Alice to create a dominating set in the Maker-Maker convention whereas Staller cannot. As a counter-example, we will prove that Bob has a drawing strategy on the cycle on ten vertices (see Theorem~\ref{thm:cycle}) whereas Dominator has a winning strategy for this graph \cite{makerbreaker}.
 However, strategies for Dominator can sometimes be used by Alice (see for example Lemma~\ref{lemma pds=gamma}). Furthermore, if Alice manages to prevent Bob from claiming a dominating set, then she can play as Dominator in the Maker-Breaker convention in the rest of the game. Using this fact, one can prove that deciding the outcome in the Maker-Maker convention is as difficult as deciding the outcome in the Maker-Breaker convention.

\begin{theorem}
{\sc Maker-Maker Domination Game} is {\sf PSPACE}-complete even if $G$ is bipartite or split.
\end{theorem}

\begin{proof}
First note {\sc Maker-Maker Domination Game} is in {\sf PSPACE}, by application of Lemma~2.2 from Schaefer~\cite{schaefer}, since there are most $|V|$ turns and during each turn, at most $|V|$ moves are available.

We prove that it is {\sf PSPACE}-hard by a reduction from {\sc Maker-Breaker Domination Game}: given a graph $G$ and a first player (Dominator or Staller), who has a winning strategy ?
{\sc Maker-Breaker Domination Game} has been proven {\sf PSPACE}-complete in \cite{makerbreaker}, even if Staller starts, and even if the graph is split of bipartite.

We do the reduction as follows. Let $G = (V,E)$ be a graph. Consider $G' = (V \cup \{v_0\}, E)$ with $v_0$ a new isolated vertex. Note that the property of being split or bipartite is maintained by this operation. We prove that Dominator wins the {\sc Maker-Breaker Domination Game} on $G$ going second, if and only if Alice wins the {\sc Maker-Maker Domination Game} on $G'$ going first.

Suppose first that Dominator has a winning strategy $\mathcal{S}$ on $G$ going second. We define the following strategy for Alice on $G'$: first claim $v_0$, then apply $\mathcal{S}$. By hypothesis, this strategy is a winning strategy for Dominator, thus, the set of vertices claimed by Alice at the end of the game will dominate the graph. As Bob cannot dominate $v_0$, he cannot dominate before her, thus Alice wins.

Reciprocally, suppose that Staller has a winning strategy $\mathcal{S}$ on $G$ going first. We define the following strategy for Bob on $G'$. If Alice does not claim first $v_0$, claims it. Alice cannot dominate $v_0$ any longer, so the outcome is at least a \Draw. Otherwise, apply $\mathcal{S}$. By hypothesis, this strategy is a winning strategy for Staller, thus, the set of vertices claimed by Alice at the end of the game will not dominate $G$, and the outcome is \Draw.
\end{proof}

Deciding the outcome of the union of graphs when the outcome of each graph is known is trivial in Maker-Breaker convention (see \cite{makerbreaker}). Thus, to study a class of graphs, one just need to consider the connected case. In Maker-Maker, the situation is slightly different, one cannot in general say anything about the union.

\subsection{Pairing strategies}

A standard strategy in positional games is the {\it pairing strategy}. Let $\mathcal S$ be a set of disjoint pair of vertices. A pairing strategy using $\mathcal S$ consist in, whenever the opponent claims a vertex of a pair, answering by claiming the other one. This strategy ensures that the player who follows it takes at least one vertex in each pair, i.e. takes a {\em transversal} of $\mathcal S$. For the Maker-Breaker domination game, this strategy can be used by Dominator if any transversal of $\mathcal S$ is a dominating set. This corresponds to the definition of a {\em pairing dominating set.}

\begin{definition}[\cite{makerbreaker}]
Let $G = (V, E)$ be a graph. A subset of pair of vertices $\mathcal S=\{(u_1, v_1),..., (u_k, v_k)\}$ of $V$ is a {\em pairing dominating set} if all the vertices are distinct and if the intersections of the closed neighborhoods of each pair cover all the vertices of the graph:
$$V = \cup_{i=1}^k N[u_i]\cap N[v_i].$$
In other words, any transversal of $\mathcal S$ is a dominating set.
\end{definition}

\begin{lemma}[\cite{makerbreaker}]\label{lem:pds}
    If $G$ has a pairing dominating set, then in the Maker-Breaker convention, Dominator wins playing first or second. 
\end{lemma}

A particular case of pairing dominating sets are perfect matchings. If $G$ is a forest, both notions are equivalent. In general, the reverse of Lemma~\ref{lem:pds} is not true. However, for some classes of graphs like forests or cographs, the reverse holds: if Dominator has a winning strategy playing second (and thus first), then $G$ admits a pairing dominating set \cite{makerbreaker}. In the Maker-Maker convention, Alice can also use pairing strategies. A first possibility for Alice is when there exists a pairing dominating set of size $\gamma(G)$, the minimum size of a dominating set.

\begin{lemma}\label{lemma pds=gamma}
Let $G$ be a graph. If $G$ has a pairing dominating set of size $\gamma(G)$, then Alice has a winning strategy in $G$.
\end{lemma}

\begin{proof}
Let $G$ be a graph. Suppose $G$ has a pairing dominating set of size $\gamma(G)$. By claiming her $\gamma$ first moves in it, and claiming in the same pair as Bob if Bob claims a first vertex in a pair, Alice can dominate in $\gamma$ moves. Bob cannot dominate before by definition of $\gamma(G)$. Therefore, Alice has a winning strategy in $G$.
\end{proof}

This lemma can be applied in particular for connected cographs. Indeed, if $G$ is a connected cograph with a universal vertex, then Alice wins at her first move. Otherwise, $G$ has a pairing dominating set of size $2$ and $\gamma(G)=2$, thus by Lemma~\ref{lemma pds=gamma}, Alice has a winning strategy.
As mentioned previously, dealing with disconnected graphs is not easy in the 
Maker-Maker convention. Actually, we did not manage to determine the outcome of a general cograph and let it as an open problem. Note that determining the outcome for cographs in Maker-Breaker convention is polynomial (see \cite{makerbreaker}), but finding the minimum number of moves needed by Dominator to win is surprisingly open~\cite{fastwin}.

Another application of pairing strategies will be when there is some $B$-trap in a position. If it is Alice's turn, she can claim the trap and then play any strategy of Dominator, like a pairing strategy. For that, we will use a pairing dominating set that are taking into consideration the vertices already claimed. We call such a pairing dominating set a {$A$-pairing}.

\begin{definition}\label{pds-position}
Let $G = (V, E)$ be a graph and $P=(G,V_A,V_B)$ a position. 
A set of disjoint pairs of unclaimed vertices $\mathcal S=\{(u_1, v_1),..., (u_k, v_k)\}$ of $V$ is a {\em $A$-pairing} of position $P$ if for each transversal $T$ of $S$, the set  $V_A  \cup T$  dominates $V$. 
\end{definition}

If Alice can prevent Bob to dominate the graph, $A$-pairings are then enough for Alice. This can be used in particular where there are some $B$-traps in the game.

\begin{lemma}\label{lem:trap_pairing}
Let $P$ be a position with a $B$-trap and a $A$-pairing. Then, we have $o(P,A)=\oA$.
\end{lemma}

\begin{proof}
Let $P$ be a position containing a $B$-trap $v$ and a $A$-pairing $\mathcal S$.  Alice playing first can claim $v$. Now, Bob cannot dominate $G$ anymore. Therefore, by following a pairing strategy using $\mathcal S$, Alice will claim a transversal of $\mathcal S$. Thus, she will dominate the whole graph and win.
\end{proof}

\begin{lemma}  \label{lem:eyes_pairing}
Let $P $  be a position with two $B$-traps $v$ and $v'$ and a $A$-pairing that contains nor $v$ nor $v'$ in $P$.
Then, we have $o(P, B) = \oA$. 
\end{lemma}  

\begin{proof}
Alice follows a pairing strategy using $\mathcal S\cup \{v,v'\}$. Bob cannot dominate since she will claim either $v$ or $v'$ and thus, there will be a vertex not dominated by Bob. She will dominate the graph since she will claim a transversal of $\mathcal S$.
 \end{proof}

\subsection{Removing leaves}

The key ingredient to solve trees in the Maker-Breaker convention is to remove leaves using the following lemma:

\begin{lemma}[\cite{makerbreaker}]\label{lem:removeleavesMB}
    Let $G$ be a graph, $\ell$ be a vertex of degree 1 and $u$ be its unique neighbor. If $u$ has degree $2$, then $G$ and $G\setminus \{\ell,u\}$ have the same outcome in Maker-Breaker convention: Dominator has a winning strategy in $G$ if and only if she has a winning strategy in $G\setminus \{\ell,u\}$ (whoever is first).
\end{lemma}

Let $T$ be a tree. One can apply successively Lemma~\ref{lem:removeleavesMB} until
obtaining a {\em reduced} tree $T'$ where the unique neighbor of any leaf has degree at least 3. Then there are only three possibilities for $T'$. Either $T'$ is empty, which means that $T$ has a perfect matching. In this case, Dominator wins playing first or second. The second possibility is that $T'$ is a single vertex or a star. Then the first player has a winning strategy. The last possibility is that $T'$ contains two disjoint {\em cherries}. A cherry is a vertex $u$ connected to two leaves. In this last case, Staller wins playing first or second.
This operation of removing leaves of a tree until obtaining a reduced tree correspond to the algorithm given in \cite{makerbreaker} to solve trees, and thus forests using unions of trees.

The situation is a quite different in the Maker-Maker convention. Indeed, Lemma~\ref{lem:removeleavesMB} is not true anymore and one cannot reduce trees as in Maker-Breaker convention. Indeed, leaves are still playing an important role as shown by Figure~\ref{fig:Non-equivalence}. Actually, one can prove that it is always optimal for Bob to claim the unique neighbor of a leave:

\begin{figure}
    \centering
    \begin{tikzpicture}
    \draw (0,0) node[v](s1) {} node[above=2mm] {};
    \draw (1,0) node[v](s2) {} node[above=2mm] {};
    \draw (-1,0) node[v](A1) {} node[above=2mm] {};
    \draw (-2,0) node[v](A2) {} node[above=2mm] {};
    \draw (0,1) node[v](B1) {} node[above=2mm] {};
    \draw (0,2) node[v](B2) {} node[above=2mm] {};
    \draw (2,0) node[v](C1) {} node[above=2mm] {};
    \draw (3,0) node[v](C2) {} node[above=2mm] {};
    \draw (1,-1) node[v](D1) {} node[above=2mm] {};
    \draw (1,-2) node[v](D2) {} node[above=2mm] {};

\draw[] (s1) -- (s2) ;
\draw[] (s1) -- (A1) ;
\draw[] (A1) -- (A2) ;
\draw[] (s1) -- (B1) ;
\draw[] (B1) -- (B2) ;
\draw[] (s2) -- (C1) ;
\draw[] (C1) -- (C2) ;
\draw[] (s2) -- (D1) ;
\draw[] (D1) -- (D2) ;

    \end{tikzpicture}
    \caption{Example of a graph where Bob wins in Maker-Maker, but if we remove any leaf with its neighbor, Alice wins.}
    \label{fig:Non-equivalence}
\end{figure}
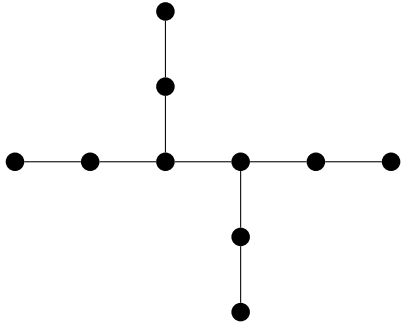

\begin{lemma}\label{lem:tail}

Let $P=(G, V_A, V_B)$ be a position with an unclaimed leaf $\ell$ for which its unique neighbor $u$ is also unclaimed.
 Then  $o(P , B ) =o(P_{\ell,u},B)$
  \footnote{In this equality,  it is assumed that $P_{\ell,u}$ is really a position, i.e. $V_B \cup \{u\}$ and $V_A \cup \{\ell\}$  do not both dominate the graph. In this latter  particular case, we obviously have $o(P , B ) = \oD$ .}.  

\end{lemma}

\begin{proof}
Suppose first that $(P_{\ell,u},B)$ is  \Draw. Bob can, in $(P, B)$, claims $u$ first. Then $\ell$ is a $A$-trap and by Lemma~\ref{lemma:trap}, Alice must claim it. Thus, the game is now $(P_{\ell,u},B)$ that is \Draw by hypothesis. Thus, $(P, B)$ is \Draw.

We prove the other implication by contraposition. Suppose $(P_{\ell,u},B)$ is \Awin. Let $\mathcal{S}$ be a winning strategy for Alice in $(P_{\ell,u},B)$. We consider a strategy $\mathcal{S'}$ for Alice in $(P,B)$ defined as follows:
\begin{itemize}
    \item If Bob claims a vertex in $\{\ell,u\}$, Alice claims the other one.
    \item Otherwise, Alice claims according to $\mathcal{S}$ (ignoring the moves on $\ell$ and $u$) until she dominates all the vertices of $V \backslash \{\ell,u\}$. If $u,\ell$ have been claimed at this moment, either Bob has claimed $u$ and Alice $\ell$, and she wins since the strategy played is the same as $\mathcal{S'}$, or Bob has claimed $\ell$ and Alice $u$, which dominates more vertices than $\ell$. Therefore, in both cases, her winning strategy in $\mathcal{S'}$ ensures her that she dominates first. Otherwise, when she dominates all the vertices of $V \backslash\{\ell,u\}$, Bob does not dominate in the game played on $(P_{\ell,u},B)$. Thus,  claiming $u$ does not make him dominate $G$. Therefore, by claiming at her next move a vertex in $\{\ell,u\}$, Alice wins. 
\end{itemize}

Thus, Alice has a winning strategy in $(P,B)$, finishing the proof.
\end{proof}

In other words, after the first move of Alice, one can always assume that Bob has claimed all the vertices that are adjacent to leaves of $G$, and that Alice has answered by claiming all the leaves.
This will be particularly important when dealing with trees.

\section{Paths and cycles}
In this section, we consider paths and cycles. The structure of the positions obtained after some moves will be basically union of paths where the extremities are claimed. Thus, we first deal with these paths and derive some general results on them that will help us to solve paths and cycles (and also forests in the next sections).

\subsection{Bounded Paths}

A {\em bounded path} is a path on at least four vertices where the four vertices at its extremities (the two leftmost and the two rightmost) are already claimed by Alice and Bob in such a way that the four vertices are already dominated by both players.

\begin{definition}
A bounded path of length $n$ is a position $(G,V_A,V_B)$  such that:
\begin{itemize}
    \item $G$  is a path $ ( v_{-1}, v_0, v_1, ...., v_n, v_{n+1}, v_{n+2})$;
    \item the unclaimed vertices are exactly vertices $v_1$ to $v_n$;
    \item exactly one vertex among $\{v_{-1},v_0\}$ (respectively $\{v_{n+1},v_{n+2}\}$) is in $V_A$, the other being in $V_B$.
\end{itemize}
\end{definition}

\begin{figure}[h]
    \centering
\scalebox{0.8}{\begin{tikzpicture}

\draw (-1,0) node[B](-1) {} node[above=2mm] {$B$};
\draw (0,0) node[R](0) {} node[above=2mm] {$A$};
\draw (6,0) node[R](6) {} node[above=2mm] {$A$};
\draw (7,0) node[B](7) {} node[above=2mm] {$B$};

\foreach \I in {1,...,5}{
            \node[v](\I) at (\I,0) {}; }

\draw[] (-1) -- (0) ;
\draw[] (0) -- (1) ;
\draw[] (1) -- (2) ;
\draw[] (2) -- (3) ;
\draw[] (3) -- (4) ;
\draw[] (4) -- (5) ;
\draw[] (5) -- (6) ;
\draw[] (6) -- (7) ;

\node at (3,-0.8) {\large $[Ao^5A]$};

\begin{scope}[shift={(10,0)}]

\draw (-1,0) node[R](-1) {} node[above=2mm] {$A$};
\draw (0,0) node[B](0) {} node[above=2mm] {$B$};
\draw (4,0) node[R](4) {} node[above=2mm] {$A$};
\draw (5,0) node[B](5) {} node[above=2mm] {$B$};

\foreach \I in {1,...,3}{
            \node[v](\I) at (\I,0) {}; }

\draw[] (-1) -- (0) ;
\draw[] (0) -- (1) ;
\draw[] (1) -- (2) ;
\draw[] (2) -- (3) ;
\draw[] (3) -- (4) ;
\draw[] (4) -- (5) ;

\node at (2,-0.8) {\large $[Bo^3A]$};

\end{scope}
\end{tikzpicture}}
    \caption{The bounded paths $[Ao^5A]$ and $[Bo^3A]$.}
    \label{fig: example bounded path}
\end{figure}
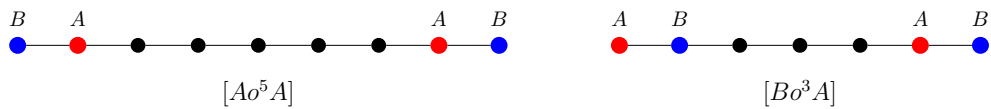

According to this definition, the knowledge of the label of $v_0$ and $v_{n+1}$ is sufficient to deduce the label of $v_{-1}$ and $v_{n+2}$. Therefore, for $t,t'$ in $\{A,B\}$, we will denote by $[t o^n t']$ the bounded path of size $n$ such that $v_0 \in V_{t}$ and $v_{n+1} \in V_{t'}$. See Figure~\ref{fig: example bounded path} for an illustration of $[Ao^5A]$ and $[Bo^3A]$.

In some situations, bounded paths can be considered as a neutral structure that preserves the outcome when adjoined to another position. The next lemma illustrates a first case in which this may occur. It will also lead to a natural resolution of paths.

\begin{lemma}  \label{lem:[A, B]}
For any position  $P$ and any integer $n$, $o(P, B)=o(P\cup  [A o^n B], B)$.
\end{lemma}

\begin{proof}
Let $P=(G,V_A,V_B)$ be a position. We do the proof by induction on the number of unclaimed vertices in $P\cup  [A o^n B]$. First note that if $n\le 1$, the result is true since $[A o^n B]$ is already dominated by both players. Hence, by Observation~\ref{obs:removedomcomp}, the result is true when the number of unclaimed vertices is at most $1$. 

Now assume there are at least two unclaimed vertices and $n \ge 2$. Suppose first that $(P, B) = \oD$. In this case, Bob starts by claiming $v_{n-1}$ in $P\cup  [A o^n B]$ creating a $A$-trap in $v_n$. By Lemma~\ref{lemma:trap}, Alice has to answer $v_n$. 
By induction hypothesis $o(P \cup  [A o^{n-2} B], B) = \oD$, which ensures that  $o(P \cup  [A o^{n} B], B) = \oD$. 
% Then the pointed position  (P\cup  [A o^n B], B) has the same outcome as $(P \cup  [A o^{n-2} B], B)$, i.e   .

Now assume that $(P, B)$ is \Awin and let $\mathcal{S}$ be a winning strategy for Alice in this pointed position. We give a strategy for Alice in $(P\cup [A o^n B],B)$. Let $v_B$ be the vertex claimed by Bob in $(P \cup  [A o^n B],B)$. Alice answers as follows:

\begin{itemize}
    \item If $v_B \in V(G)$ and if Alice does not dominate $P$, she claims the same vertex $v_A$ she would have answered following $\mathcal{S}$ if Bob had claimed $v_B$ in $(P, B)$. The resulting position has two vertices less, and is \Awin by induction hypothesis.
    \item If $v_B \in V(G)$ and Alice already dominates $P$, we necessarily have $n\ge 2$, otherwise, she would already have won. She claims $v_2$. The resulting position is better than $(P\cup [AoA] \cup [Ao^{n-2}B],B)$, which itself is better than $([AoA] \cup [Ao^{n-2}B],B)$, as Alice already dominates $P$. Thus, by induction hypothesis applied with $P' = [AoA]$, which is not dominated by Bob, the position is \Awin.
    \item If $v_B$ is a vertex of $[A o^n B]$, then $v_B$ corresponds to some vertex $v_k$, with $k\in\{1,\ldots ,n\}$. If $k=n$, Alice answers $v_{n-1}$. Otherwise, she claims $v_{k+1}$. Note that if $k=n$, the position obtained is better for Alice than the one obtained by $k = n-1$ (since she dominates a superset of vertices). Therefore, and without loss of generality, we will assume that $k < n$. By Lemma~\ref{lem:split}, the resulting position is then equivalent to $(P \cup [A o^{k-1} B]\cup [A o^{n-1-k} B], B)$ which has less unclaimed vertices than the original one. By induction hypothesis applied twice, it has the same outcome as $(P \cup [A o^{k-1} B] ,B)$, which also has the same outcome as $(P,B)$, which is \Awin by hypothesis.
\end{itemize}
This analysis ensures that $o(P\cup  [A o^n B], B) = \oA$, since  for each claim of Bob, there exists an answer of Alice  leading to an \Awin position. 
\end{proof}

The next lemma presents another situation where the adjunction of some bounded paths with particular constraints does not change a winning outcome for Alice. 

\begin{lemma}\label{lem:AABB}
Let  $P=(G, V_A, V_B) $ be a position, let $n$ be an integer such that $n \not \equiv 0 \mod 3$ and let $n'$ be an integer such that $n' \equiv 0 \mod 2$.
Then if  $o(P,B) = \oA$, then 
$o(P \cup  [Ao^n A] \cup [Bo^{n'} B] ,B) = \oA$.
\end{lemma}

\begin{proof}
%We assume that Alice has a winning strategy in $(P,B)$. We prove 
%if $n   \leq 2 $ and $n'  = 0$ then the result is trivial. Otherwise, we proceed by induction on the number of unclaimed vertices of $P \cup  [Ao^n A] \cup [Bo^{n'} B]$. 

% that Alice has a winning strategy in $((P \cup  [Ao^n A] \cup [Bo^{n'} B] ,B) $.
We prove this result  by induction on the number $m$ of unclaimed vertices of $P \cup  [Ao^n A] \cup [Bo^{n'} B]$. 
For initialization, if $m = 1$,  then $P$ has no unclaimed vertices, $n = 1$ and $n' = 0$. 
Thus,  Alice  dominates $P \cup  [A o^n A]  \cup [Bo^{n'}B]$ and the result is true, by definition.

Assume now that $m \geq 2$, 
If Alice  dominates $(P \cup  [A o^n A]  \cup [Bo^{n'}B])$ (which implies that  $n   \leq 2 $ and $n'  = 0$) then the result is true, by definition. 
Otherwise, we have to prove that, for each vertex $x$ claimed by Bob, there exists an answer $y$ for Alice such that 
$o(P \cup  [Ao^n A] \cup [Bo^{n'} B] ,B) = \oA$, the pointed position obtained after the two claim is \Awin.

%Bob does not dominate $P$ and $[A,A]_n$ (since $n\geq 1$), thus he cannot dominate at the next move. 
Denote by $(u_1, \dots, u_n)$ the unclaimed vertices of $[Ao^nA]$ and by $(v_1, \dots, v_{n'})$ the unclaimed vertices of $[Bo^{n'}B]$. 
We consider all the possible moves for Bob.

\begin{itemize}
    \item If Bob claims $u_1$ or $u_2$ in $[A o^n A]$ and $n\leq 2$, we consider two subcases. If $n'=0$, then Alice dominates $ [Ao^n A] \cup [Bo^{n'} B]$,  which obviously gives the result. 
   If $n' \geq 2$,  Alice claims $v_1$. By Observation~\ref{obs:removedomcomp}, the resulting game is equivalent to $P\cup [Ao^{n'-1}B]$, which is equivalent to $P$ by Lemma~\ref{lem:[A, B]}. 
    \item If Bob claims $u_1$ or $u_2$ in $[A o^n A]$ and $n\geq 4$, by Lemma~\ref{lem:sommetdomine}, we can suppose he claims $u_2$. then Alice claims $u_3$. Thus, the resulting position  is equivalent to $P\cup [Ao^{n-3}A]\cup [Bo^{n'}B]$, which gives the result using the  induction hypothesis. The case where Bob claims $u_n$ or $u_{n-1}$ is symmetric.
    \item If Bob claims a vertex $u_k$ in $[Ao^nA]$ with $k \notin \{1,2,n-1,n\}$, Alice answers by claiming a vertex $u_i$ with $i \in \{k-1,k+1\}$ such that the resulting component is $[Ao^jA]\cup[A o^{j'}B]$ with $j \not \equiv 0 [3]$. Note that, since $n \not \equiv 0 [3]$, this vertex always exists. By Lemmas~\ref{lem:split} and~\ref{lem:[A, B]}, the resulting position is equivalent to $P \cup [Ao^jA]\cup [Bo^{n'}B]$, which gives the result using the  induction hypothesis.
        \item If Bob claims some vertex $v_k$ in $[Bo^{n'}B]$, Alice answers by claiming a vertex $v_i$ with $i \in \{k-1,k+1\}$ such that the resulting component is $[Bo^jB]\cup[Ao^{j'}B]$ with $j \equiv 0 [2]$. Note that as $n' \equiv 0 [2]$, this vertex always exists. By Lemma~\ref{lem:[A, B]} and Lemma~\ref{lem:split}, the resulting position is equivalent to $P \cup [Ao^nA] \cup [Bo^jB]$, which gives the result using the  induction hypothesis.
        \item If Bob claims a vertex $x$ in $P$, we consider two subcases.  If Alice does not dominate $G$ yet, she claims the vertex she would have claimed as an answer to $x$ in her winning strategy in $P$. The resulting position has two vertices less,  which gives the result using the  induction hypothesis.

        If Alice already dominates $G$, and $n' \geq 2$, then Alice  claims $v_1$ in $[Bo^{n'}B]$, turning it into $[Ao^{n'-1}B]$. By Lemma~\ref{lem:[A, B]}, the game is now equivalent to $P \cup [Ao^nA]$, and Alice already dominates $G$, so by induction hypothesis, it is a winning position for her.
        For $n'  = 0$ and  $n \geq 4$ (cases where $n \leq 2$ are trivial), then Alice  claims $u_2$.  The resulting position  is better for Alice than the position  $P \cup [AoA] \cup   [Ao^{n-2}B]$, which, by Lemma~\ref{lem:[A, B]}, is equivalent to $P \cup [AoA]$  and therefore is  a winning position for Alice.  \qedhere
\end{itemize}
\end{proof}

We finish this subsection by proving that bounded paths $[Bo^nB]$, where $n$ is odd, and $[Ao^nA]$, when $n\equiv 0 \mod 3$, are not good for Alice when Bob starts. The first one is natural and give a condition for a \Draw whatever the rest of the position is. The second one is more surprising. Here Bob obtains a \Draw mostly by threatening Alice to dominate before her. Thus, one cannot add any position and maintain a \Draw.

\begin{lemma}  \label{lem:[BB]_odd}
For any position  $P$ and any odd integer $n$, $o(P\cup  [B o^n B], B) = \oD$.
\end{lemma}

\begin{proof}
We prove the result by induction on $n$.
If $n =1$, there is a $A$-trap: by Lemma~\ref{lemma:trap}, $o(P\cup  [BoB], B) = \oD$.
Now, let $n \geq 3$.  Bob claims $v_2$  which forces Alice to claim $v_1$. The position is then equivalent to $(P\cup  [Bo^{n-2}B], B)$ which is \Draw by induction.
\end{proof}

\begin{lemma}\label{lem:AA}
For any positive integer $k$, $o([Ao^{3k}A],B)=\D$.
\end{lemma}

\begin{proof}
We prove the result by induction on $k$.

If $k = 1$, Bob claims $v_2$ and directly dominates $[Ao^3A]$. Thus, $o([Ao^3A],B)=\D$.

If $k = 2$, Bob claims $v_2$. If Alice does not answer in $\{v_1, v_3,v_4\}$, Bob claims $v_3$ at his second turn and create two $A$-traps in $v_1$ and $v_4$ which ensures a \Draw. Thus, Alice should claim a vertex among $\{v_1, v_3,v_4\}$ and does not dominate the graph at her first claim. Then Bob can claim $v_5$ and dominates the graph. Hence $o([Ao^6A],B)=\D$.

Assume now that $k\geq 3$ and that the result is true for any positive $k'<k$. Consider the position $([Ao^{3k}A],B)$.
Bob claims $v_5$ at his first turn.
By Lemma~\ref{lem:P5} applied on $(v_3,v_4,v_5,v_6,v_7)$, Alice should answer on $v_4$ or $v_6$.
If she claims $v_4$, then by Lemma~\ref{lem:cut}, the position is equivalent to the position
$([Ao^3A]\cup [Bo^{3k-5}A], B)$, which is equivalent, by Lemma~\ref{lem:[A, B]} to $([Ao^3A], B)$ 
which is \Draw.

If she claims $v_6$, in the same way, the position is equivalent to $([Ao^4B]\cup [Ao^{3(k-2)}A], B)$, which is equivalent to $([Ao^{3(k-2)}A], B)$ which is \Draw by induction.
\end{proof}

\subsection{Paths}

Lemma~\ref{lem:[A, B]} can be used to prove that Alice always wins on paths.

\begin{theorem}
Let $P_n$ be the path of length $n$. Then $o(P_n)=\oA$.
\end{theorem}
\begin{proof}
Let $v_1,...,v_n$ be the vertices of the path. In $n\leq 3$, then Alice wins at her first turn, so we can assume that $n\geq 4$.
    Alice starts by claiming $v_2$. By Lemma~\ref{lem:tail}, we can assume that Bob claims $v_{n-1}$ and Alice answers by claiming $v_n$. 
    Let $(P, B)$ be the actual pointed position of the game. 
    Using Observation~\ref{obs:sameWS}, and since $v_1$ should be in any winning set of Bob, this position is equivalent to the position $([A o A]\cup[A o^{n-4} B],B)$.
     By Lemma~\ref{lem:[A, B]},  $o([A o A]\cup[A o^{n-4}  B],B) = o([Ao A],B) = \oA$,  which ensures that $o(G)=\oA$. 
   % Obviously  $o(G) \geq o(P, B) $.
    % One can extend $P$ by adding two vertices $v_{-1}$, claimed by Bob, and $v_0$, claimed by Alice. Let $P'$ be this new position.  $(P, B) $  is clearly equivalent to $(P', B) $ since the vertices to be dominated are the same for both players.
    % We have $o(P', B) \geq o([AoA]\cup[A o^{n-4} B], B)$. To see it,  first 
%     
     %Now add a vertex $u$  only linked to $v_2$ to $P'$ and assume that $u$ is claimed by Bob.  Obviously, the obtained position $P''$ is such that   $o(P', B) \geq o(P'', B)$. Using Lemma~\ref{lem:split}, with $X =\{u, v_2\}$, $o(P'', B) = o(([A o A]\cup[A o^{n-4}  B],B)$.  
 % One cannot directly apply Lemma~\ref{lem:split} to $P'$ to obtain the position $[AoA]\cup[A o^{n-4} B]$, since $v_2$ is not dominated by Bob. However, Bob will have to claim $v_1$ to dominate it, and this way will dominate $v_2$. Thus we can consider that Bob does not have to take into consideration the domination of $v_2$. Then the reasoning in the proof of Lemma~\ref{lem:split} can be repeated in this special case. We then obtain that the actual pointed position is equivalent to $([AA]_1\cup[AB]_{n-4},B)$. 
      %By Lemma~\ref{lem:[A, B]},  $o([A o A]\cup[A o^{n-4}  B],B) = ([Ao A],B) = \oA$,  which ensures that $o(G)=\oA$. 
\end{proof}

When playing in the Maker-Breaker convention, note that this result implies that Dominator also has a winning strategy on any path playing first. This result was already known since~\cite{makerbreaker}, but the strategy developed here is different from the other one.\\

\subsection{Cycles}

The case of cycles in the Maker-Maker convention is more subtle than for the Maker-Breaker convention, where Dominator always wins. More precisely, we will show that there are infinitely many \Awin and \Draw outcomes that depend on the size of the cycle modulo $3$. 

From now, we will denote by $C_n$ the cycle of order $n$. The vertices of $C_n$ will be denoted by $v_0$ to $v_{n-1}$.

\begin{theorem} \label{thm:cycle}
Let $n$ be an integer. We have $o(C_n) = \oD$ if and only if $n \ge 10$ and $n \equiv 1 \mod 3$. 
\end{theorem}

\begin{proof}
%\begin{lemma}\label{claim : Bob circle}
%Let $n$ be an integer. $o(C_n) = \oD$ if $n \ge 10$ and $n \equiv 1 \mod 3$.
%\end{lemma}
%\begin{proof}
We first prove the “if” part. Let $n = 3k+1$, with $k \geq 3$,  and let $C_n$ be a cycle of order $n$
By symmetry, we can assume that Alice first claims $v_0$ and thus $o(C_n) =  o(C_n, \{v_0\}, \emptyset, B)$. We give a strategy for Bob to obtain at least a draw. Bob first claims $v_5$. By Lemma~\ref{lem:P5} with $G'=(v_3,v_4,v_5,v_6,v_7)$, Alice has to answer $v_4$ or $v_6$.
\begin{itemize}
\item If Alice claims $v_6$. Then Bob can claim $v_3$ which forces Alice to claim the $A$-trap in $v_4$ and then Bob can claim $v_1$ which forces Alice to claim the $A$-trap in $v_2$. At this point, the position is equivalent to $([Ao^{3(k-2)}A],B)$ which is \Draw by Lemma~\ref{lem:AA} (remember that $k-2\geq 1$).

\item If Alice claims $v_4$. Bob  can claim successively all the $v_{2i+1}$ starting from $v_7$, creating a $A$-trap in $v_{2i}$ that Alice is forced to claim.
If $n$ is even, Bob follows this strategy until claiming $v_{n-1}$. Then Alice has to claim $v_{n-2}$ and does not dominate the cycle (she does not dominate $v_2$). Then Bob can dominate the cycle by claiming~$v_2$.

If $n$ is odd, Bob follows this strategy until claiming $v_{n-4}$. After Alice has claimed  $v_{n-5}$,  Bob claims $v_2$. 
Alice is forced to claim $v_{n-1}$, otherwise Bob wins by claiming it. Then Bob claims $v_{n-2}$, creating a trap in $v_{n-3}$, forcing Alice to claim it. At this point, Alice still does not dominate $v_2$. Then Bob can dominate the cycle by claiming $v_1$.
\end{itemize}
In all cases, the game either ends in a draw or Bob dominates the graph, thus $o(C_n) = \oD$.

\

We now prove the “only if” part. First consider that $n \not \equiv 1 \mod 3$.
%\begin{lemma}
%\label{claim:Alice cycle}
%    Let $n$ be an integer such that $n \not \equiv 1 \mod 3$. Then $o(C_n) = \oA$.
%\end{lemma}
%\begin{proof}
Without loss of generality, one can assume that Alice will claim $v_0$ at her first turn, and thus we consider the pointed position $(P,B)=((C_n,\{v_0\},\emptyset),B)$.
Consider now the position $P'$ obtained from $P$ by adding a vertex claimed by Bob adjacent only to $v_0$. This position is better for Bob since he dominates more vertices than in $P$. Thus, it is enough to prove that $(P',B)$ is \Awin to ensure that $(P,B)$ is \Awin. 
Using Observation~\ref{obs:sameWS}, $(P',B)$ is equivalent to $([Ao^{n-1}A],B)$. By Lemma~\ref{lem:AABB}, since $n-1\not \equiv 0 \mod 3$, $o([Ao^{n-1}A],B)=\oA$, and thus $o(C_n) = \oA$

It remains to prove that Alice wins on $C_4$ and $C_7$. On $C_4$, Alice wins by claiming any two vertices and as it is not possible to dominate in one move, she will dominate first. On $C_7$, Alice can claim $v_0$, and then use a pairing strategy with pairs $(v_2,v_3)$ and $(v_4,v_5)$. This way, she dominates the cycle in three moves, and Bob cannot dominate before since at least three vertices are required to dominate $C_7$. 
\end{proof}

\section{Forests}
%\section{Structure of the Algorithm for Acyclic Graphs}

In this section, we give the essential elements to solve the case of forests, given by Theorem~\ref{thm:trees} below. In particular, we will first reduce the problem to any {\em standard} forest, meaning that non-standard forests correspond to cases that can be solved easily. Then we will give necessary conditions about the first move of Alice, yielding to the introduction of the skeleton of a forest. From that definition, we will be able to present the general algorithm that computes the outcome of any forest, as depicted by the decision tree of Figure~\ref{fig:algo}. In the next section, we will give the full proof of the validity of the algorithm. 

\begin{theorem}\label{thm:trees}
Deciding the outcome of a forest can be done in linear time.
\end{theorem}

\subsection{Removing small components}

If there is an isolated vertex in the forest, this is like playing in the Maker-Breaker convention.

\begin{lemma}\label{lem:isolatedvertex}
    Let $F$ be a forest with an isolated vertex $v_0$. Then $o(F)=\oA$ if and only if $F\setminus \{v_0\}$ contains a perfect matching.
\end{lemma}

\begin{proof}
Alice has to claim first $v_0$. Then Bob is playing first in $F\setminus \{v_0\}$. Since he cannot dominate anymore, he has the same role as Staller in the Maker-Breaker Domination Game. In \cite{makerbreaker}, it is proved that Dominator playing second in a forest in the Maker-Breaker Domination Game wins if and only if there is a perfect matching, which gives the result.
\end{proof}

Isolated edges can be removed without changing the outcome.
\begin{lemma}\label{lem:isolatededge}
    Let $F$ be a forest with an isolated edge $e=uv$. Then $o(F)=o(F\setminus \{u,v\})$.
\end{lemma}

\begin{proof}
If Alice (respectively Bob) has a winning (resp. draw) strategy in $F\setminus \{u,v\}$, then she can apply her strategy in $F$ by pairing $u$ with $v$. The resulting strategy will still be winning (resp. leading to a draw) in $F$.
\end{proof}

By applying Lemma~\ref{lem:isolatedvertex} and Lemma~\ref{lem:isolatededge}, one can consider in what follows that all the connected components have at least three vertices.

\subsection{Bottom-to-top strategies for Bob}\label{sec:bottom-to-top}

In this subsection, we describe a strategy for Bob that will often be considered to obtain draws or to reduce trees. Let $T$ be a tree rooted on a vertex $r$, vertices of $T$ except $r$ are  labeled inductively with values $0$ and $1$, starting from the leaves as follows:
 \begin{itemize}
 \item  If all children of $v$ are labeled $1$, then $v$ is labeled $0$  (hence,  all the leaves are labeled $0$);
 \item   If at least a  child  of $v$ is  labeled by $0$, then $v$ is labeled by $1$. 
\end{itemize}

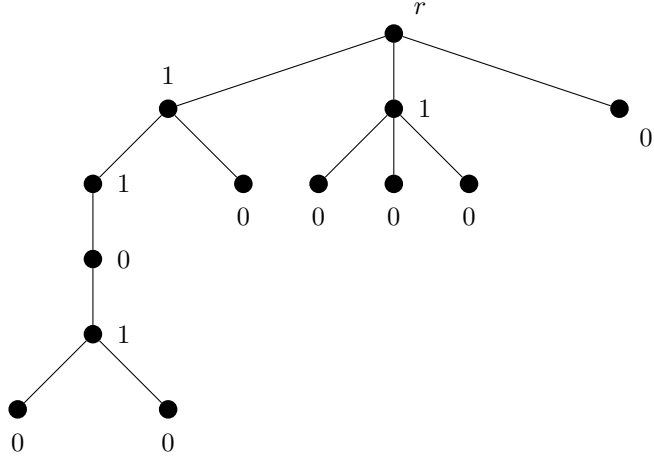
\begin{figure}
    \centering

\begin{tikzpicture}

    \draw (5,5) node[v](r) {} node[above right=2mm] {$r$};

    \draw (2,4) node[v](a1) {} node[above=2mm] {$1$};
    \draw (5,4) node[v](a2) {} node[right=2mm] {$1$};
    \draw (8,4) node[v](a3) {} node[below right=2mm] {$0$};

\draw[] (r) -- (a1) ;
\draw[] (r) -- (a2) ;
\draw[] (r) -- (a3) ;

    \draw (1,3) node[v](a11) {} node[right=2mm] {$1$};
    \draw (3,3) node[v](a12) {} node[below=2mm] {$0$};

\draw[] (a1) -- (a11) ;
\draw[] (a1) -- (a12) ;

    \draw (4,3) node[v](a21) {} node[below=2mm] {$0$};
    \draw (5,3) node[v](a22) {} node[below=2mm] {$0$};
    \draw (6,3) node[v](a23) {} node[below=2mm] {$0$};

\draw[] (a2) -- (a21) ;
\draw[] (a2) -- (a22) ;
\draw[] (a2) -- (a23) ;

    \draw (1,2) node[v](b1) {} node[right=2mm] {$0$};

\draw[] (a11) -- (b1) ;

    \draw (1,1) node[v](c1) {} node[right=2mm] {$1$};

\draw[] (b1) -- (c1) ;

    \draw (0,0) node[v](d1) {} node[below=2mm] {$0$};
    \draw (2,0) node[v](d2) {} node[below=2mm] {$0$};

\draw[] (c1) -- (d1) ;
\draw[] (c1) -- (d2) ;

\end{tikzpicture}

    \caption{Example of labelling}
    \label{fig:0-1labelling}
\end{figure}

Figure~\ref{fig:0-1labelling} gives an example of such a labeling.

Let $T$ be a tree rooted in $r$ and consider the pointed position 
$(P,B)$ with $P=(T,V_A,\emptyset)$ and $V_A\subseteq \{r\}$. 
A \emph{bottom-to-top strategy} for Bob on $(P,B)$ consists, at each step, in claiming a vertex $v$ labeled by $1$ such that all the successors of $v$ labeled by $1$ are already claimed by Bob. The following property is maintained during this process: Alice is forced to claim only vertices labeled by $0$, and any vertex claimed by Alice (except $r$) has his parent claimed by Bob.
Indeed, it is true before the first claim of Bob. Assume it is true before Bob's turn. Let $v$ be a vertex labeled by $1$ such that all the successors of $v$ labeled by $1$ are already claimed by Bob. By definition of the labeling and the assumption, $v$ has a child $u$ labeled by $0$ that is unclaimed. All the children of $u$ are by definition labeled by $1$ and thus already claimed by Bob. Thus, $u$ is a $A$-trap and Alice is forced to claim it, maintaining the property true. Such a strategy can also be applied when all the leaves of $T$ are adjacent to vertices already claimed by Bob.
In this strategy, Bob can claim all the vertices labeled by $1$. A particular interesting case for Bob is when there exists a vertex $v$ labeled by $1$ with two children labeled by $0$:

\begin{lemma}
 Let $T$ be a tree and consider the labeling of $T$ rooted in $v_0$. If there exists a vertex labeled by 1 with two children labeled by $0$, then the position  $(P,B)$ with $P=(T,\{v_0\},\emptyset)$ is \Draw.
\end{lemma}

\begin{proof}
    Bob uses a bottom-to-top strategy on $T$. When claiming $v$, he will create two $A$-traps on the two children of $v$ labeled by $0$, which concludes the proof by Corollary~\ref{cor:2traps}.
\end{proof}

Bob can also use bottom-to-top strategies to reduce the forest to a smaller one where he has a \Draw strategy.

\begin{lemma} \label{lem:cut}
Let $F$ be a forest and $v_0$ a vertex. Consider the labeling of $F$ rooted in $v_0$ (for the components not containing $v_0$, root on any vertex). Let $v \neq v_0$ be a vertex labeled by 1 and $S_v$ be the set of successors of $v$ in the rooted tree. Let $F_{\overline{v}}$  be the tree obtained by removing all vertices of  $S_v$ and adding a leaf $v'$ connected to $v$.  %(and edges whose at least an endpoint is in $S_v$). We have the implication: 
 Then, if $o(F_{\overline{v}},\{v_0,v'\},\{v\},B)=\oD$, then $o(F,\{v_0\},\emptyset,B)=\oD$. 
\end{lemma}

\begin{proof}
Bob follows by a bottom-to-top strategy on $S_v$ and can claim all the vertices of $S_v$ labeled by $1$ until claiming~$v$. Alice is always forced to claim a child of the claimed vertex.
Afterward, Bob plays his strategy to obtain a \Draw in
$(F_{\overline{v}},\{v_0,v'\},\{v\},B)$, leading to a \Draw  for  $(T,\{v_0\},\emptyset,B)$.
  \end{proof}

\subsection{Cherries}

A particular case where the bottom-to-top strategy will be useful is when there are some cherries in the forest. Recall that a cherry is a vertex $c$ connected to two leaves $\ell_1$ and $\ell_2$. It will be denoted by  the triple $C = (c, l_1, l_2)$.
If $F$ contains some cherries, the outcome of $F$ can be easily computed.

\begin{lemma}\label{lem:2cherries}
Let $F$ be a forest. If $F$ has two cherries  $C = (c, \ell_1, \ell_2)$ and $C' = (c', \ell_1', \ell_2')$, with $c \neq c'$, then $o(F) = \oD$. 
\end{lemma}

\begin{proof}
Let $F$ be such a forest. After Alice has claimed her first vertex, she cannot have claimed both  $c$ and  $c'$. Suppose without loss of generality that Alice has claimed $c'$. Then Bob claims $c$. 
The resulting pointed position contains two  $A$-traps and thus  is \Draw by Corollary~\ref{cor:2traps}.
\end{proof}

\begin{lemma}\label{lem:1cherry}
Let $F$ be a forest with  exactly one cherry $C = (c,\ell_1,\ell_2)$. Then  $o(F) = \oA$ if and only if there is a matching in $F\setminus\{c\}$ that covers $V(F)\setminus N[c]$.
\end{lemma}

\begin{proof}
%Let $T$ be a forest. %Suppose $T$ has exactly one cherry $(c,l,l')$
Suppose first that $F$ has a matching $M$ in $F\setminus\{c\}$ that covers $V(F)\setminus N[c]$. Then Alice claims $c$, which creates a double $B$-trap in $\ell_1$ and $\ell_2$. The matching $M$ is actually a $A$-pairing, and, we can suppose it contains neither $\ell_1$ nor $\ell_2$ since these two vertices have their neighborhood included in $N[c]$. Then by Lemma~\ref{lem:eyes_pairing}, Thus $o(F) = \oA$. \\

%Following this strategy, Alice ensures that she will play at least once in each pair of the matching and at least one of $l$ or $l'$ (if the game goes until she has to play one of them). Bob cannot dominate both $l$ and $l'$, indeed, as Alice has played $c$ the only one way for Bob to dominate them is to play them. But if he does, Alice plays the other one. So Bob cannot dominate. Moreover, as we supposed that $M$ covers $V(T)\setminus (\{c\} \cup N[c])$, after all the vertices are played, Alice dominates $V(T)\setminus (\{c\} \cup N[c])$ threw her moves on $M$ and she dominates $N[c]$ as she has played $c$. Thus, she has a winning strategy.\\

Suppose now that $F$ has no such matching $M$. %Consider the labeling defined in Subsection~\ref{labeling}.   
We define a strategy for Bob as follows. If Alice's first claim is not element of  the cherry, then Bob claims  $c$ and creates two $A$-traps, leading to a \Draw position.
Thus, we can assume that Alice's first claim $r$ is and element of $\{c, \ell_1, \ell_2\}$.
Let $T$ be the connected component of $F$ containing $r$. If there exists another connected component $T'$ of $F$ that has no perfect matching, then Bob can apply the strategy of Staller playing first in $T'$ that prevents Dominator to dominates $T'$ (see \cite{makerbreaker}). Thus, one can assume that there is a perfect matching in all the components of $F$ distinct from $T$. 
Now label the vertices of $T$ rooted in $r$ as defined in Section~\ref{sec:bottom-to-top}.
We want to prove that there exists a vertex labeled $1$ with two children labeled $0$, which will ensure a \Draw strategy for Bob by applying the bottom to top strategy.
If it is not the case, then consider the matching $M'$ where all the vertices labeled $1$ (except $c$) are paired with their unique child labeled $0$. We claim that $M'$ covers $V(T)\setminus N[c]$. Indeed, assume $x\in V(T)\setminus N[c]$ is not covered by $M$. Then $x$ must be labeled $0$ and its parent should be $r$ since it is the only vertex not labeled. But then $x\in N[c]$. Thus, there exists a matching in $F\setminus \{c\}$ that cover $V(F)\setminus N[c]$: take the union of $M'$ and the perfect matchings of all the other components.
%  \begin{itemize}
%     \item If such a component contains a vertex adjacent to a neighbor of $c$, then Bob plays a bottom to top  strategy in $T \setminus \{l, l'\}$ rooted in $c$. As we supposed that no matching of $T\setminus\{c\}$ covers $V(T)\setminus (\{c\} \cup N[c])$, there exists a vertex $v$ labeled by $1$  with (at least) two children, otherwise the pair of moves of Alice and Bob would create a matching that dominates the component.   %at a certain moment, Bob will have played all the vertices labeled one, but at least one vertex labeled $0$ will still be available, otherwise, all the pairs of moves from Bob and Alice would form this matching.
%     This ensures that the outcome is \Draw  from Remark~\ref{rem:bottomtotop}. 
%     %by construction of the labeling. 
%        % When Bob plays $v$, at the next move, he can play one of his  children labeled $0$ says $v'$, and hence, he prevents Alice to dominate $v'$, as all the neighbors of $v'$ are labeled $1$, and are already played by the bottom to top strategy. He ensures that the outcome is $\oD$ by construction of the labeling. 
% %
% %  
% 	\item otherwise, consider a connected component with no matching, and root it in an arbitrary vertex $v_0$. 
%         Bob plays the bottom to top strategy in  this rooted tree.  As in the previous case, there exists a vertex $v$ labeled by $1$  with (at least) two children,  which ensures that the outcome is \Draw  from Remark~\ref{rem:bottomtotop}.    
%       \end{itemize}
\end{proof}

\subsection{Definition of the skeleton and easy cases}

Considering Lemmas~\ref{lem:2cherries} and~\ref{lem:1cherry}, we will assume now that there is no cherry in $F$, no isolated vertex and no isolated $P_2$. 
From Lemma~\ref{lem:tail}, one can assume that after the first turn of Alice, Bob will claim all the unclaimed vertices of $F$ with a leaf as a neighbor. Alice is then forced to answer on each leaf.
This motivates us to define the {\em skeleton} of $F$, denoted by $\sk_F$, as the vertices that are not a leaf nor a parent of a leaf. 

\begin{figure}
    \centering
\begin{tikzpicture}[square/.style={regular polygon,regular polygon sides=4}, triangle/.style={regular polygon,regular polygon sides=3}]

    \draw (3,2) node[v](b2) {} ;
    \node at (3,2) [square,draw] (s) {};
    \draw (3,1) node[v](b1) {} ;
    \node at (3,1) [triangle,draw] (t) {};

    \draw (1,0) node[v](a2) {} ;
    \node at (1,0) [square,draw] (s) {};

    \draw (2,0) node[v](a1) {} ;
    \node at (2,0) [triangle,draw] (s) {};

    \draw (3,0) node[v](v2) {} ;
    \draw (4,0) node[v](v1) {} ;
    \draw (5,0) node[v](v0) {} ;
    \node at (5,0) [triangle,draw] (s) {};

    \draw (5,1) node[v](l1) {} ;
    \node at (5,1) [square,draw] (s) {};

\draw[] (v0) -- (v1) ;
\draw[] (v0) -- (l1) ;
\draw[] (v1) -- (v2) ;

\draw[] (a1) -- (a2) ;
\draw[] (a1) -- (v2) ;

\draw[] (b1) -- (b2) ;
\draw[] (b1) -- (v2) ;

\draw (5,-1) node[v](r) {};
\draw[] (v0) -- (r) ;

\draw (5,-2) node[v](f1) {};
\draw (5,-3) node[v](f2) {};

\draw[] (r) -- (f1) ;
\draw[] (f1) -- (f2) ;

\draw (4,-4) node[v](c1) {};
\draw (4,-5) node[v](c2) {};

\draw (4,-6) node[v](d11) {};
\draw (4,-7) node[v](d12) {};
\node at (4,-6) [triangle,draw] (s) {};
\node at (4,-7) [square,draw] (s) {};

\draw[] (f2) -- (c1) ;
\draw[] (c1) -- (c2) ;

\draw[] (c2) -- (d11) ;
\draw[] (d12) -- (d11) ;

\draw (5,-4) node[v](c1) {};
\draw (5,-5) node[v](c2) {};

\draw (5,-6) node[v](d11) {};
\draw (5,-7) node[v](d12) {};
\node at (5,-6) [triangle,draw] (s) {};
\node at (5,-7) [square,draw] (s) {};

\draw[] (f2) -- (c1) ;
\draw[] (c2) -- (c1) ;
\draw[] (c2) -- (d11) ;
\draw[] (d12) -- (d11) ;

\draw (6,-4) node[v](c1) {};
\draw (6,-5) node[v](c2) {};

\draw (6,-6) node[v](d11) {};
\draw (6,-7) node[v](d12) {};
\node at (6,-6) [triangle,draw] (s) {};
\node at (6,-7) [square,draw] (s) {};

\draw[] (f2) -- (c1) ;
\draw[] (c1) -- (c2) ;
\draw[] (c2) -- (d11) ;
\draw[] (d12) -- (d11) ;

 \filldraw[fill= white, fill opacity = .1, style = dashed] ($(v2)+(-.5,.5)$) 
        to[] ($(v1) + (.5,.5)$) 
        to[] ($(v1) + (.5,-.5)$)
        to[] ($(v0) + (1.5,-.5)$)
        to[] ($(c2) + (.5,-.5)$)
        to[] ($(v2) + (-.5,-5.5)$)
        to[] ($(v2) + (-.5,.5)$);

\node at (7,-1.5) {\Large $S_F$};

\end{tikzpicture}
    \caption{A forest $F$ and its skeleton $S_F$. Note that $F$ is connected but $S_F$ is not. $M_F$ is the set of vertices in triangles, and $L_F$ is the set of vertices in squares.}
    \label{fig:skeleton}
\end{figure}
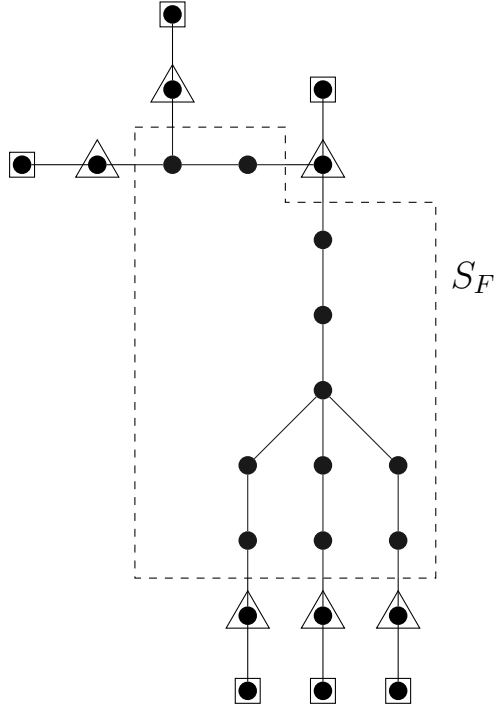

More formally, we denote by $L_F$ the leaves of $F$ and by $M_F$ their parents. Remark that, with the hypothesis that there is neither   cherry nor isolated vertex, we have $L_F \cap M_F = \emptyset$ and the mapping :   $L_F \to M_F$, which associates to each vertex $v$ of $L_F$ its parent,  is bijective.  Then, let $\sk_F$ be defined by   $\sk_F = V(F) \setminus (L_F \cup M_F)$. Figure~\ref{fig:skeleton} is an illustration of a forest with its skeleton. 

In some simple cases, we can directly give the outcome of $F$.

\begin{lemma}\label{lem:skvide}
If $\sk_F$ is empty, then $o(F)=\oA$.
\end{lemma}

\begin{proof}
Set $L_F = \{\ell_1, \dots, \ell_k\}$. By hypothesis, $F$ has no cherry, therefore, no vertex can be adjacent to two of these leaves. Thus, we can denote $M_F = \{m_1, \dots, m_k\}$ with $m_i$ adjacent to $\ell_i$ for $1 \le i \le k$.  
Note that $k$ vertices  (one on each  set  $\{\ell_k, m_k\}$) are necessary and sufficient to dominate $F$. Thus $\gamma(F)=k$ and there is a pairing dominating set of size $k$ (the set of pairs $\{(\ell_i,m_i),i\in \{1,...,k\}$). By Lemma~\ref{lemma pds=gamma}, Alice has a winning strategy in $F$.
%  The strategy for Alice is as follows. 
%  The first claim is $m_1$. 
%  Afterwards, if Bob claims a vertex of a pair $\{l_i, m_i\}$ and the other vertex of the pair is unclaimed, then   then Alice answers by claiming this vertex. 
% If the  other vertex of the pair has been  previously claimed, Alice claims $m_{i'}$, for a pair $\{l_{i'}, m_{iu'}\}$ whose  vertices are both unclaimed. 
% By this way, at each claim of Alice, the number of pairs dominated by Alice  is increasing.  Thus after $k$ claims, Alice dominates $T$, while Bob cannot dominate $T$ since he has claimed only $k-1$ vertices. Thus, we have $o(T) = \oA$. 
\end{proof}

\begin{lemma}\label{lem:rmvSvide}
 Let $T$ be a connected component of $F$ such that $S_T$ is empty. Then $o(F)=o(F\setminus T)$.   
\end{lemma}

\begin{proof}
Assume first that $o(F\setminus T)=\oA$ and let $x$ be the first claim of Alice in a winning strategy. Then Alice claims $x$ in $F$. By Lemma~\ref{lem:tail}, Bob will claim all the vertices of $M_T$ and Alice will answer all the vertices of $L_T$. After these moves, $T$ is completely claimed and dominated by both players. By Observation~\ref{obs:removedomcomp}, we can remove this component. The game is then equivalent to $(F\setminus T,\{x\},\emptyset)$ which is \Awin when Bob starts, leading to $o(F)=\oA$.

Assume now that $o(F\setminus T)=\oD$. Consider the game played in $F$ and a first claim $x$ of Alice. If $x\in V(F\setminus T)$, then as before, the position $(F,\{x\},\emptyset)$  is equivalent to the position $(F\setminus T,\{x\},\emptyset)$ which is \Draw when Bob starts.
If $x\in V(T)$, let $x'$ be the unique neighbor of $x$ if $x$ is a leaf or the leaf connected to $x$ if $x\in M_T$. Then we can assume by Observation~\ref{obs:removedomcomp} that Bob claims all the vertices in $M_T$ (except $x$ or $x'$)  and that Alice answers by claiming all the vertices in $L_T$ (except $x$ or $x'$). At this point, all the vertices of $T$ have been claimed except $x'$. Then Bob claims $x'$. Then both players dominate $T$ and the position is equivalent to $(F\setminus T,A)$ which is a \Draw position.
In conclusion, whatever Alice claims, Bob can ensure a draw in $F$. Thus $o(F)=\oD$
\end{proof}

%\begin{proof}
%We will prove that $T$ satisfies $\gamma(T) = pds(T)$, and Lemma~\ref{lemma pds=gamma} to obtain the result.

 %Moreover, as no vertex can dominate two leaves, we have $\gamma(T) \ge k$.

%Now, let prove that the $(l_i,m_i)$s for $1 \le i \le k$ are a pds of $T$. In fact, let $x$ be a vertex of $T$. As $\sk_F$ is empty, we have either $x \in L_F$ or $x \in M_F$. Therefore, there exists a $1 \le i \le k$ such that $x \in \{l_i,m_i\}$. Now, as $(l_i, m_i)$ are adjacent, $x$ is both adjacent to $l_i$ and $m_i$. Thus, the $(l_i,m_i)$s for $1 \le i \le k$ are a pds of $T$ of size $k$.

%Finally, we have $\gamma(T) = pds(T) = k$. So $o(T) = \oA$ by Lemma~\ref{lemma pds=gamma}.
%\end{proof}

 % {\color{red} Eric R: cette preuve me semble bien compliquée car elle renvoie à d'autres notions, alors que presenter une stratégie gagnante ne semble pas si difficile}
 
\begin{lemma}\label{lem:star}
If $\sk_F$ induces a star of center $c$ such that $c$ has no neighbor in $M_F$, then $o(F) = \oA$.
\end{lemma}

\begin{proof}
%If $\sk_F$ is a star of center $c$, and if $c$ has no neighbor in $M_F$, then Alice can start by playing $c$. 
Alice claims the center $c$ as a first move. Then as explained above, Bob claims all the vertices of $M_F$ and Alice answers all the leaves of $L_F$. After that, Alice dominates the whole graph and Bob does not, since he does not dominate $c$.
%all the vertices she does not dominate are in $M_F$ or $L_F$. Once again, by denoting by $l_i$ the vertices of $L_F$ and $m_i$ the vertices of $M_F$ (for $1 \le i \le k$). If Alice plays once in each pair $(l_i, m_i)$, she dominates the graph. Therefore, she has a strategy to dominate in $k+1$ moves. But, we have $\gamma(T) = k+1$ since  all vertices  $l_i$,  $1 \le i \le k$,  and $c$ has no common neighbor, thus they need to be dominated by different vertices. So Bob cannot prevent Alice to dominate in $k+1$ moves and cannot dominate in $k$ moves. Thus, we have $o(T) = \oA$. 
%Alice plays $c$ then she needs $|N(c)|$ moves to dominate the graph, Bob needs $|N(c)| + 1$.
\end{proof}

We say that that $F$ is {\em standard} if all its connected components have at least four vertices and a non-empty skeleton, if $F$ has no cherries, and if $S_F$ does not induce a star of center $c$ where $c$ has no neighbor in $M_F$. We now focus on the standard forests.

\subsection{First move of Alice}

If $F$ is standard, next lemmas say that it can be assumed that Alice must claim outside $\sk_F$ and must connect it. The main idea behind this result is that, if Alice claims in $\sk_F$, she claims too far from the leaves and Bob can win with a bottom-to-top strategy.

\begin{lemma}\label{lem:firstmovesk}
Let $F$ be a standard forest and $v_0 \in  \sk_F$. 
We have  $o((F, \{v_0\}, \emptyset), B) = \oD$.  
\end{lemma}

\begin{proof}
Assume first that the graph induced by  $\sk_F$ is  either a unique vertex, a unique edge,   or  a star (whose center has necessarily a neighbor in $M_F$).  Let $L_F = \{\ell_1, \ldots, \ell_k\}$ and $M_F = \{m_1, \ldots, m_k\}$
such that for $1\leq i \leq k$, $\ell_i$ and $m_i$ are neighbors. Starting from $(F, \{v_0\}, \emptyset)$, Bob successively claims $m_1, \ldots, m_k$,  which forces Alice to reply $\ell_1, \ldots, \ell_k$.  When Bob has just claimed $m_k$, he dominates the whole graph while Alice does not yet dominate $\ell_k$. This ensures that  $o((F, \{v_0\}, \emptyset), B) = \oD$.

It can now be assumed  that 
the graph induced by  $\sk_F$ is  neither a unique vertex, a unique edge,   nor  a star. 
 Thus, we can consider  the position $P'=(F, \{v_0\} \cup L_F, M_F )$, and from successive applications of   Lemma~\ref{lem:tail}, we just need to prove that $o(P',B)=\oD$. 

%By Lemma~\ref{lem:tail}, Bob claims one by one all the vertices of $M_F$ and Alice should answer by claiming the corresponding vertices in $L_F$, leading to the equivalent position $(P',B)$ with $P'=(F, \{v_0\} \cup L_F, M_F )$. We just need to prove that $o(P',B)=\oD$. 
%{\color{red} Note that since $F$ is standard, $\sk_F$ is not a star whose center has no neighbor in $M_F$. In the case where $\sk_F$ is still a star but has a neighbor in $M_F$,  in $P'$,  Bob will dominate the center of the star. All the leaves of the star will then be dominated by Bob Through the claims  in $M_F$, otherwise, they would be in $L_F$, and thus Bob will dominate first as he will dominate by playing the last vertex of $M_F$ while Alice won't dominate its private neighbor.  

%If $\sk_F$ is not a star. In particular, Alice does not dominate $F$ in $P'$ which means that at least one vertex of $\sk_F$ is not dominated.}

First focus on components of $\sk_F$ which do not contain $v_0$. 
Let $C$ such a component, if each vertex of $C$ is dominated by a vertex of  $M_F$, then Bob dominates C while Alice does not. 
Otherwise, let $v_1$ be a leaf of the subtree $T_C$ of $F$ induced by $C$. Bob plays the bottom-to-top strategy in $T_C$ rooted in $v_1$. It can be done since the leaves of $T_C$ are only connected to vertices of $M_F$ already claimed by Bob.  
If  two vertices labeled by 0 have the same parent, then the strategy creates a double trap, which ensures a \Draw. If $v_1$ has all its children labeled by 1, when Bob  claims the last vertex labeled by  1, a double trap is created, since $v_1$ will be a trap. This also ensures the \Draw.
The only non directly conclusive case  is when  $v_1$ has one child labeled by 0, the reached position after Bob has followed the bottom-to-top strategy is such that Bob dominates $C$ while Alice does not. Indeed, Bob dominates $v_1$ since it is a leaf of $C$ and thus should be connected to $M_F$ but Alice does not dominate $v_1$. 

Bob follows this strategy on each component of $\sk_F$ not containing $v_0$. All the answers from Alice are forced.  If a double trap appears, we are done. 
Otherwise, each component is dominated by Bob and not by Alice. Moreover, it is Bob's turn. In such a case , we now need to focus on the component $C_0$ of $\sk_F$ which contains $v_0$.  

\begin{itemize}
    \item  If the  diameter of $C_0$  at most 1,  (which implies that the graph induced by  $\sk_F$ is not connected), Bob already dominates $C_0$, and therefore $F$,  but Alice does not.
    \item If the diameter of $C_0$ is $2$,  the graph induced by $C_0$ is a star centered in a vertex $c$.  
    Since it is assumed that the graph induced by $\sk_F$ is not a star, $C_0$ is not the only component of the graph induced by $\sk_F$. 
    Thus, if Bob does not dominate $C_0$ yet, he claims $c$ or any of its neighbor and dominates the whole graph while  Alice does not.  
    \item If the diameter of $C_0$ is  at least $3$, there exists a vertex $v_2$ at distance exactly 2 from $v_0$ in the subgraph induced by $C_0$.  
    Let  $(v_0,v_1,v_2)$  be  the path from $v_0$ to $v_2$. Bob can then use a bottom-to-top strategy in the tree induced by $C_0$ rooted in $v_0$.
        \begin{itemize}
        \item  If $v_1$ is labeled by  $0$, then $v_2$ is labeled by $1$. Bob can play a bottom-to-top strategy claiming all the vertices labeled by 1 except $v_2$. Then he claims $v_1$. At this moment, Bob dominates  $F$  but Alice  does not dominate $v_2$.
        \item If $v_1$ is labeled by $1$, Bob plays a bottom-to-top strategy in all branches, finishing by claiming $v_1$. At this moment, he dominates $F$ but, by construction, Alice does not dominate $v_2$.
    \end{itemize}
\end{itemize}
Thus, in any case, Bob has a strategy which leads to a \Draw position, and thus $o((F, \{v_0\}, \emptyset), B) = \oD$.
\end{proof}

\begin{lemma}\label{lem:firstmove}
If $F$ is standard and $o(F)= \oA$, then there exists a vertex $v_0 \in M_F$ such that $o((F, \{v_0\}, \emptyset), B) =\oA$. Moreover,  $v_0$ satisfies:
\begin{enumerate}
\item  the subgraph of $F$ induced by $\sk_F\cup \{v_0\}$  is a tree ;
 \item  in the labeling of $F$ rooted in $v_0$, each vertex $v$ labeled by 1 has a unique child.
\end{enumerate}
\end{lemma}

\begin{proof}
    As it is   supposed  that $o(F,A) = \oA$, there exists 
    a vertex  $v_0 \in V(F)$ such that $o((F, \{v_0\}, \emptyset), B) =\oA$.
    From Lemma~\ref{lem:firstmovesk}, we have $v_0 \in M_F \cup L_F$. From Lemma~\ref{lem:sommetdomine}, it can be assumed that $v_0 \notin L_F$. Indeed, if $l \in L_F$ and $m \in M_F$ is its private neighbor, $N[l] \subset N[m]$.

%Assume that there exists such a $v_0 \in L_F$. Let $v'_0$ be the unique neighbor of $v_0$, and compare positions  $((F, \{v_0\}, \emptyset), B)$ and $((F, \{v'_0\}, \emptyset), B)$.  Since $o((F, \{v_0\}, \emptyset), B) = \oA$, there exists a winning strategy for Alice from this position. 
 %    Assume that the same strategy is used by Alice in $(F, \{v'_0\}, \emptyset)$, (with the light modification that if Bob claims 
   % $v_0$, them Alice answers as if Bob had claimed $v'_0$). At each step of the game, the set of vertices dominated by Alice starting from  $(F, \{v_0\}, \emptyset)$ is included in the set of vertices dominated by Alice starting from  $(F, \{v'_0\}, \emptyset)$  at the same step, and we have the reverse inclusion for the set of vertices dominated by Bob. 
   %It follows that $o((F, \{v_0\}, \emptyset), B)  \leq  o((F, \{v'_0\}, \emptyset), B)$, which implies that $o((F, \{v'_0\}, \emptyset), B)= \oA$.

  %  it can  be assumed that $v_0 \in M_F$ (since, when   $v_0 \in  L_F$, we have $o(T, \{v_0\}, \emptyset, B) =$ \Awin  $\implies o(T, \{v_1\}, \emptyset, B) =$ \Awin, where $v_1$ is the neighbor of $v_0$)
    %Alice has a winning move outside $\sk_F$. As any move in $M_F$ is better than a move in $L_F$, she has a winning move in $M_F$.

Thus, we can assume that $v_0\in M_F$. We now prove the two other properties.

\begin{enumerate}
    \item Since $F$ has no cherry, there exists a unique vertex $v_{-1} \in L_F$ which is a neighbor of $v_0$. Assume that the graph induced by $\sk_F \cup \{v_0\}$ is not connected. We prove that, under this hypothesis,   $o((F, \{v_0\}, \emptyset), B) =\oD$, which gives the result by contraposition.
    First, using Lemma~\ref{lem:tail}, we have 
    $o((F, \{v_0\}, \emptyset), B) = o((F, \{v_0\} \cup L_F \setminus \{v_{-1}\}, M_F \setminus \{v_{0}\}), B)$.

   With the same arguments as in the previous lemma, the implementation of a bottom-to-top strategy on each component of the graph induced by $\sk_F \cup \{v_0\}$ not containing $v_0$ leads to either a double trap (which gives the result), or a position where each of these components is dominated by Bob but not by Alice.
    Now, focus on the component $C$ containing $v_0$. Bob then plays a bottom-to-top strategy on the tree induced by $C$ rooted in $v_0$.
    Then Bob claims $v_{-1}$ after all the vertices labeled by 1 have been claimed. By this way, if no trap has appeared before, Bob dominates the whole forest $F$, while Alice does not totally dominate $F$.  This gives the result.     
   % Let $C_1, \dots, C_k$ be the connected components of $T(\sk_F \cup \{v_0\}$. Up top play bottom to top strategies in $C_j$ for $3 \le j \le k$, and so making these component dominated by Bob, we can suppose $k=2$. Suppose that we have $v_0 \in C_1$. First, 
    
    %Bob plays a bottom to top strategy in $C_1$ rooted in $v_0$. As $T$ has no cherry, he can do it by forcing all the moves of Alice until only $v_0$ and, if it exists, a leaf $v_{-1}$ connected to it are not dominated. Now, let $v_1$ be the vertex of $T$ that disconnected $C_1$ from $C_2$. By construction, $v_1 \in M_F$ as the leaves cannot disconnect a graph. let $x$ be the first neighbor of $v_1 \in C_2$. Bob plays in $C_2$ rooted in $x$ a bottom to top strategy stopping before playing $x$. At this moment, Bob has dominated $C_2$ by construction of the labeling, but Alice does not dominate $x$. He can now play $v_{-1}$ and he dominates $T$ before Alice. Note that if $v_{-1}$ does not exist, if he has not dominated yet, he can do it by playing any neighbor of $v_0$. Finally, if $o(T) = \oA$, then $v_0$ have to connect $\sk_F$.

    \item  For the second item, consider the labeling of $F$ (that is actually a tree) rooted in $v_0$. First note that all the leaves of $\sk_F$ are labeled by $0$. Thus, each vertex labeled by $1$ has at least one child. Assume that there exists a vertex $v$ labeled by 1 with has at least two children.  
    \begin{itemize}
        \item If two children $v'$ and $v"$ of $v$ are labeled by 0, Bob can then use a bottom-to-top strategy until $v$ is claimed. When he claims $v$, he creates a double $A$-trap in $v'$ and $v"$ and thus obtains a \Draw.
        \item Otherwise, there exists one child $v'$  of $v$ labeled by 0 and the other one, $v''$ labeled by 1. 
        Then Bob uses a bottom-to-top strategy but without claiming $v"$, until all the vertices labeled by 1 (except $v"$ are claimed). All replies of Alice remain forced.  After this is done,  Bob claims $v_{-1}$, the neighbor of $v_0$ which belongs to $L_F$, and then dominates the whole graph while Alice does not dominate $v"$, which ensures that $o((T, \{v_0\}, \emptyset), B) = \oD$. \qedhere
    \end{itemize}
    \end{enumerate}
\end{proof}

%Corollaire à garder ?
%\begin{corollary}
%Consider $F_A(T)$ $T$ where Alice has played her first move in $S_T$, then Bob has played all the other available moves in $S_T$ and Alice has played their neighbors.

%Alice wins in $T$ if and only if there exists $A \in S_T$ such that $o(F_A(T),\{A\}, S_T\backslash A,B) = 1$
%\end{corollary}

As a consequence of this result, if $F$ is standard and winning for Alice, then $F$ is necessarily a tree.

\subsection{Splitting the graph}

By Lemma~\ref{lem:firstmove}, it can be assumed that Alice first claims a vertex $v_0\in M_F$ that is connected to all the components of $\sk_F$  (note that if $\sk_F$ is not connected, there is at most one such vertex).
It can also be assumed that each vertex labeled by 1 in $F$ rooted in $v_0$ has degree exactly 2.
In all the remaining, $v_0$ will denote this first claim of Alice.
Let $v_{-1}$ be the leaf connected to $v_0$. After this first move, it can be assumed, using Lemma~\ref{lem:tail}, that  Bob will claim  all the other vertices of $M_F$ one by one. At each time, Alice must answer to the corresponding leaf in $L_F$. After this step, the free vertices are the vertices in $\sk_F$ with the vertex $v_{-1}$.  Formally, the obtained position is $ P = (F, L_F \setminus \{v_{-1}\} \cup  \{v_0\}, M_F \setminus \{v_0\})$ and we have 
$o(F)  = o(P, B)$. In that follows, we will split the graph into several components defined from the connected components of the skeleton.

\begin{definition}
For a connected component $C$ of $\sk_F$, let $T$ be the connected component of $F\setminus\{v_0\}$ that contains $C$. The position $P_C$ is defined as the position induced by $T\cup \{v_{0}, v'_{0}\} $ in the position $P = (F, L_F \setminus \{v_{-1}\} \cup  \{v_0\}, M_F\cup \{v'_0\} \setminus \{v_0\})$, where $v'_0$ is an additional leaf connected to $v_0$ and claimed by Bob.
\end{definition}

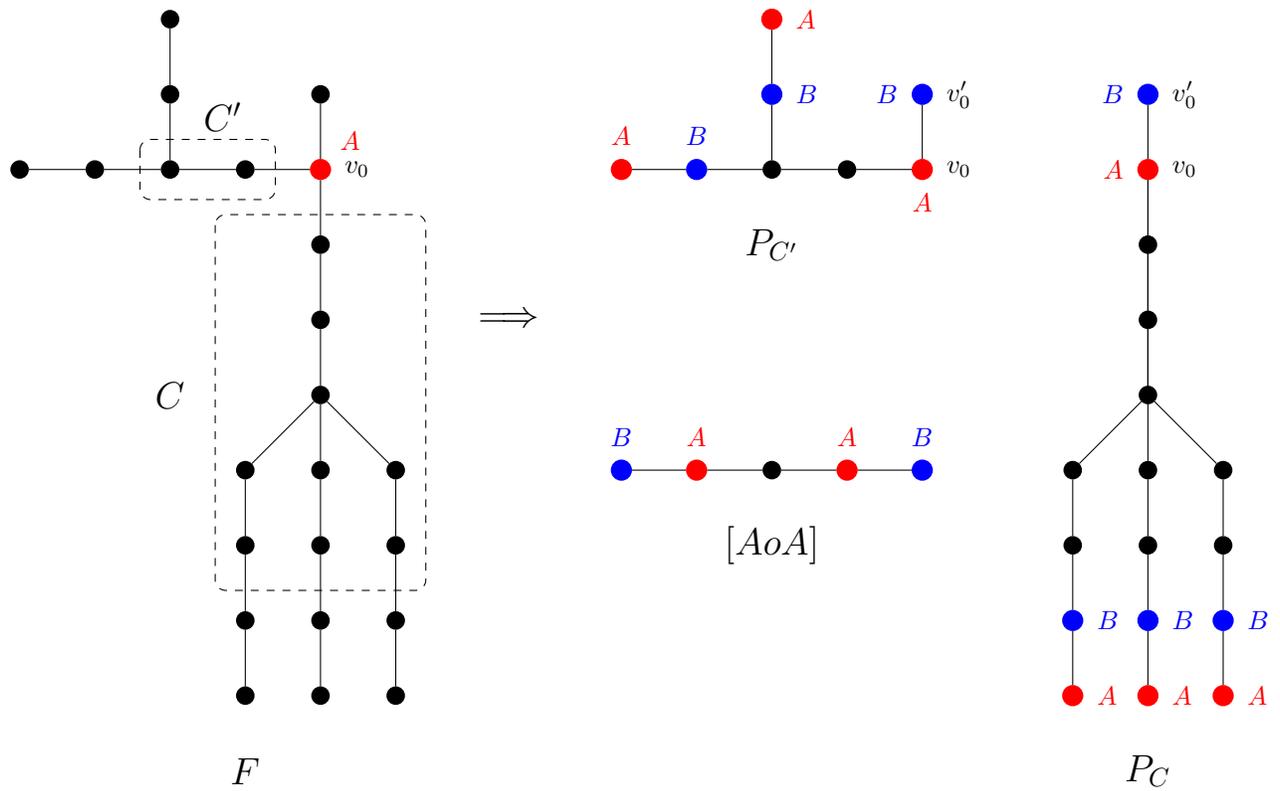
\begin{figure}
    \centering
\begin{tikzpicture}

%%Arbre général 

    \draw (3,2) node[v](b2) {} ;
    \draw (3,1) node[v](b1) {} ;

    \draw (1,0) node[v](a2) {} ;
    \draw (2,0) node[v](a1) {} ;

    \draw (3,0) node[v](v2) {} ;
    \draw (4,0) node[v](v1) {} ;
    \draw (5,0) node[R](v0) {} node[right=2mm] {$v_0$} node[above right=2mm] {\color{red}$A$};

    \draw (5,1) node[v](l1) {} ;

\draw[] (v0) -- (v1) ;
\draw[] (v0) -- (l1) ;
\draw[] (v1) -- (v2) ;

\draw[] (a1) -- (a2) ;
\draw[] (a1) -- (v2) ;

\draw[] (b1) -- (b2) ;
\draw[] (b1) -- (v2) ;

\draw (5,-1) node[v](r) {};
\draw[] (v0) -- (r) ;

\draw (5,-2) node[v](f1) {};
\draw (5,-3) node[v](f2) {};

\draw[] (r) -- (f1) ;
\draw[] (f1) -- (f2) ;

\draw (4,-4) node[v](c1) {};
\draw (4,-5) node[v](c2) {};

\draw (4,-6) node[v](d11) {};

\draw (4,-7) node[v](d12) {};

\draw[] (f2) -- (c1) ;
\draw[] (c1) -- (c2) ;

\draw[] (c2) -- (d11) ;
\draw[] (d12) -- (d11) ;

\draw (5,-4) node[v](c1) {};
\draw (5,-5) node[v](c2) {};

\draw (5,-6) node[v](d11) {};

\draw (5,-7) node[v](d12) {};

\draw[] (f2) -- (c1) ;
\draw[] (c2) -- (c1) ;
\draw[] (c2) -- (d11) ;
\draw[] (d12) -- (d11) ;

\draw (6,-4) node[v](c1) {};
\draw (6,-5) node[v](c2) {};

\draw (6,-6) node[v](d11) {};

\draw (6,-7) node[v](d12) {};

\draw[] (f2) -- (c1) ;
\draw[] (c1) -- (c2) ;
\draw[] (c2) -- (d11) ;
\draw[] (d12) -- (d11) ;

\draw[rounded corners, style = dashed] (3.6, -.6) rectangle (6.4, -5.6) {};

\draw[rounded corners, style = dashed] (2.6, .4) rectangle (4.4, -.4) {};

\node at (4,-8) {\Large $F$};

\node at (3,-3) {\Large $C$};

\node at (3.7,.7) {\Large $C'$};

%%%%% Arbre de C

\draw (16, 1) node[B](v0p) {} node[right=2mm] {$v_0'$} node[left=2mm] {\color{blue}$B$};
\draw (16, 0) node[R](v0) {} node[right=2mm] {$v_0$} node[left=2mm] {\color{red}$A$};
\draw (16,-1) node[v](r) {};
\draw[] (v0) -- (r) ;
\draw[] (v0p) -- (v0) ;

\draw (16,-2) node[v](f1) {};
\draw (16,-3) node[v](f2) {};

\draw (15,-4) node[v](c1) {};
\draw (15, -5) node[v](c2) {};

\draw (15,-6) node[B](d11) {} node[right=2mm] {\color{blue}$B$};

\draw (15,-7) node[R](d12) {} node[right=2mm] {\color{red}$A$};

\draw[] (v0) -- (f1) ;
\draw[] (f1) -- (f2) ;

\draw[] (f2) -- (c1) ;
\draw[] (c1) -- (c2) ;

\draw[] (c2) -- (d11) ;
\draw[] (d12) -- (d11) ;

\draw (16,-4) node[v](c1) {};
\draw (16, -5) node[v](c2) {};

\draw (16,-6) node[B](d11) {} node[right=2mm] {\color{blue}$B$};

\draw (16,-7) node[R](d12) {} node[right=2mm] {\color{red}$A$};

\draw[] (r) -- (c1) ;
\draw[] (c1) -- (c2) ;

\draw[] (c2) -- (d11) ;
\draw[] (d12) -- (d11) ;

\draw (17,-4) node[v](c1) {};
\draw (17, -5) node[v](c2) {};

\draw (17,-6) node[B](d11) {} node[right=2mm] {\color{blue}$B$};

\draw (17,-7) node[R](d12) {} node[right=2mm] {\color{red}$A$};

\draw[] (r) -- (f1) ;
\draw[] (f1) -- (f2) ;
\draw[] (f2) -- (c1) ;
\draw[] (c1) -- (c2) ;

\draw[] (c2) -- (d11) ;
\draw[] (d12) -- (d11) ;

\node at (16,-8) {\Large $P_{C}$};

%%% Arbre de C'

    \draw (11,2) node[R](b2) {} node[right =2mm] {\color{red}$A$};
    \draw (11,1) node[B](b1) {} node[right=2mm] {\color{blue}$B$};

    \draw (9,0) node[R](a2) {} node[above =2mm] {\color{red}$A$};
    \draw (10,0) node[B](a1) {} node[above =2mm] {\color{blue}$B$};

    \draw (11,0) node[v](v2) {} ;
    \draw (12,0) node[v](v1) {} ;
    \draw (13,0) node[R](v0) {} node[right=2mm] {$v_0$} node[below =2mm] {\color{red}$A$};
    \draw (13,1) node[B](l1) {} node[right=2mm] {$v'_0$} node[left=2mm] {\color{blue}$B$} ;

\draw[] (v0) -- (v1) ;
\draw[] (v0) -- (l1) ;
\draw[] (v1) -- (v2) ;

\draw[] (a1) -- (a2) ;
\draw[] (a1) -- (v2) ;

\draw[] (b1) -- (b2) ;
\draw[] (b1) -- (v2) ;

\node at (11, -1) {\Large $P_{C'}$};

\draw (9, -4) node[B](w0) {} node[above =2mm] {\color{blue}$B$};
\draw (10, -4) node[R](w1) {} node[above =2mm] {\color{red}$A$};
\draw (11, -4) node[v](w2) {} node[above =2mm] {};
\draw (12, -4) node[R](w3) {} node[above =2mm] {\color{red}$A$};
\draw (13, -4) node[B](w4) {} node[above =2mm] {\color{blue}$B$};

\draw[] (w0) -- (w1) ;
\draw[] (w1) -- (w2) ;
\draw[] (w2) -- (w3) ;
\draw[] (w3) -- (w4) ;

\node at (11, -5) {\Large $[AoA]$};

\node at (7.5, -2) {\Large $\Longrightarrow$};

\end{tikzpicture}
    \caption{On the left a forest $F$. On the right, the equivalent position obtained by splitting $F$ in Lemma~\ref{lem:bigsplit}. The component $C$ is a fork whereas $C'$ is a $P_2$.}
    \label{fig:position_PC}
\end{figure}

Figure~\ref{fig:position_PC} illustrates the two positions $P_C$ and $P_{C'}$ derived from the forest $F$ of Figure~\ref{fig:skeleton} when played on $v_0$. Lemma~\ref{lem:bigsplit} shows that this splitting (with an additional $[AoA]$ position) yields to an equivalent position.

\begin{lemma} \label{lem:bigsplit}
The position $ P = (F, L_F \setminus \{v_{-1}\} \cup  \{v_0\}, M_F \setminus \{v_0\})$ is equivalent to the position $$\left(\bigcup_{C \in CC(\sk_F)}P_C \right)\cup [AoA]$$
  where $CC(\sk_F)$ denotes the set of connected components of $\sk_F$.
\end{lemma}

\begin{proof}
% Let  $P'$ be the position $P$ with an additional leave $v'_0$ adjacent to $v_0$ and claimed by Bob.
% We have  $o(F)  = o(P, B) = o(P',B)$. The last equality is easy since, obviously,   $o(P', B) \leq o(P,B)$, and conversely, if Bob has a strategy on $P'$ leading to \Draw, then the same strategy also gives Draw in $P$. Indeed, if Bob 
% succeeds as a breaker in $P'$, he is also a breaker in $P$.  Otherwise, If Bob 
% dominates in $P'$ before Alice, then he has necessarily claimed $v_{-1}$. Thus 
% Bob dominates in $P'$ before Alice. 

% Then the result is   directly obtained by successive applications of  Lemma~\ref{lem:split} on $P'$,   
%  with $X = \{ v_{0}, v'_{0}\}$, the element  $[AoA]$ being equivalent to the component 
%    $(F(v_{-1},  v_{0}),  \{v_{0}\} , \emptyset)  $.  We obtain $o(P',B) = o( (\bigcup_{C \in CC(\sk_F) }P_C )  \cup [AoA], B) $, which gives the result. 

Let $P'$ be the position $(\bigcup_{C \in CC(\sk_F)}P_C) \cup [AoA]$.
Note that the unclaimed vertices are in a  one-to-one correspondence in the two positions ($v_{-1}$ is corresponding to the unclaimed vertex of $[AoA]$).
By Observation~\ref{obs:sameWS}, we just need to prove that $P$ and $P'$ have the same winning sets for both players.

A set $S$ of unclaimed vertices is winning for Alice in $P$ and in $P'$ if and only if it is dominating all the vertices of $\sk_F$ except the ones connected to (a copy of) $v_0$, which corresponds to the same condition in both positions.
A set $S$ of unclaimed vertices is winning for Bob in $P$ if and only if $v_{-1}\in S$ and $S$ is dominating all the vertices of $\sk_F$ that are not connected to a vertex of $M_F\setminus \{v_0\}$. In $P'$, $S$ is winning for Bob if it contains the unclaimed vertex of $[AoA]$ and if it is dominating all the unclaimed vertices not already dominated by Bob, that are exactly all the vertices of $\sk_F$ not connected to $M_F\setminus \{v_0\}$. Thus, the winning sets are in bijection and by Observation~\ref{obs:sameWS}, the positions are equivalent.
   \end{proof}

%, and two vertices   $v_{-2}, v_{-3}$, mutually connect, such that $v_{-1}$ is also connected to $v_0$. Moreover, in $P'$, $v'_0$ and  $v_{-2} $ are claimed by Bob, and $v_{-1}$  is 
%The next lemma uses the above result to give a sufficient condition for having a draw on $F$.
Using this decomposition, we now prove that Bob can just focus on a subset of components where he has a \Draw strategy.

\begin{lemma} \label{lem:component}
 Assume there exists a set $\mathcal S$ of connected components of $\sk_F$ such that $o((\bigcup_{C \in \mathcal S} P_C)\cup [A o A],B) =\oD$. Then
 $o(F)=  \oD$. 
\end{lemma}

\begin{proof} 
Using Lemma~\ref{lem:bigsplit}, 
it suffices to prove that $o((\bigcup_{C \in CC(\sk_F)}P_C) \cup [AoA], B)  =  \oD $.   
For each $C \notin \mathcal S $, Bob plays  a bottom-to-top  strategy on  $P_{C}$ rooted on $v_0$.   At each time Alice is forced to answer in the same component $C$,  and, at the end, Bob dominates the component $C$. 
This is successively done for all such components.
 
Afterward, the remaining position is equivalent to
$((\bigcup_{C \in \mathcal S} P_C)\cup [A o A], B)$  which is \Draw by hypothesis.  
 \end{proof}

%We first consider the case where $\sk_F$ induces a connected graph, as there is an affordable characterization for deciding whether a position is \Awin or not. For that purpose, we define a {\em fork} as follows:

\subsection{Favorable skeletons for Alice}
In this subsection, we give some necessary and sufficient conditions for a component $C$ to be winning for Alice.

\begin{definition}
A fork is a star with at least three branches where each is subdivided exactly once.
\end{definition}

On Figure~\ref{fig:position_PC}, the component $C$ is an example of a fork with four branches.

% \begin{figure}
%     \centering
% \begin{tikzpicture}

%     \draw (0,0) node[v](s) {} node[above=2mm] {};

%     \draw (1,1) node[v](a1) {} node[above=2mm] {};
%     \draw (2,1) node[v](a2) {} node[above=2mm] {};

%     \draw (1,0) node[v](b1) {} node[above=2mm] {};
%     \draw (2,0) node[v](b2) {} node[above=2mm] {};

%     \draw (1,-1) node[v](c1) {} node[above=2mm] {};
%     \draw (2,-1) node[v](c2) {} node[above=2mm] {};

% \draw[] (s) -- (a1) ;
% \draw[] (a1) -- (a2) ;

% \draw[] (s) -- (b1) ;
% \draw[] (b1) -- (b2) ;

% \draw[] (s) -- (c1) ;
% \draw[] (c1) -- (c2) ;
    
% \end{tikzpicture}
%     \caption{A fork}
%     \label{fig:fork}
% \end{figure}

\begin{lemma}  \label{lem:fork_or_path}
If $((F,\{v_0\},\emptyset),B)$ is \Awin, then all the connected components of $\sk_F$ must induce a path  or a fork that is connected to $v_0$ by a leaf.
\end{lemma}

\begin{proof}

 Proceeding by contraposition, assume that there exists a component  $C$ of $\sk_F$, which is neither a fork nor a path connected to $v_0$ by a leaf . We will prove  $o(P_C\cup [AoA], B) = \oD$, which gives the result using Lemma~\ref{lem:component}.
 
 Since $C$ is not a path connected by a leaf to $v_0$, there exists a vertex $c \in C$ of degree at least 3 in the tree induced by $C\cup \{v_0\}$. Let $P_a=(c, a_1, a_2, ...a_{p})$ be  the path linking $c$ to $v_0$ ($a_p$ is adjacent to $v_0$), $P_b=(c, b_1, b_2, ..., b_{q}) $ be a path of $C$ of maximal length such that $b_1  \neq  a_1$ and  $P_c=(c, c_1, c_2, ..., c_{r})$ be another maximal path of $\sk_F$ starting in $c$. Note that possibly $p=0$, but that $q\geq r\geq 1$.
 In the labeling of $F$ rooted in $v_0$, the vertex $c$ is necessarily labeled by 0 since it has degree 3. This implies that both $q$ and $r$ are even since $b_q$ and $c_r$, as leaves of $C$, are labeled by $0$.
 Moreover, all the vertices labeled by $1$ have degree 2. This implies that all the vertices of the three paths $P_a$, $P_b$ and $P_c$ that are connected to other vertices of $F$ must be labeled by $0$. In particular, using Lemma~\ref{lem:cut}, it is enough to prove that there is a \Draw strategy when $C$ is reduced to these three paths.
 
Assume first that $q \geq 4$. By Lemma~\ref{lem:cut},  it is enough to give a \Draw strategy for Bob for  $q = 4$ and $r = 2$ with $C$ reduced to the union of the three paths. 
The first claim of Bob is $b_2$. By Lemma~\ref{lem:P5}, Alice should claim either $b_1$ or $b_3$.

 \begin{itemize}
\item  If Alice replies by claiming $b_3$, then Bob claims $c$, which forces Alice to claim $b_1$. Then, he successively claims $a_{2}$, $a_{4}$, and so on until $a_{p-1}$ is dominated by $B$. Successive replies of Alice are forced: when Bob claims $a_{2i}$, Alice necessarily replies in $a_{2i-1}$. Finally, Bob claims in $[AoA]$ and gets a \Draw by dominating  before Alice (Alice does not dominate $c_1$ and $c_2$). 
 
\item If Alice replies by claiming $b_1$, then Bob claims $c_1$, which forces Alice to claim $c_2$. Then Bob claims the unclaimed vertex of $[AoA]$. Then the position is equivalent to the position $([Ao^{p-1}B]\cup[Bo^2B],A)$. By Lemma~\ref{lem:[A, B]}, this position is equivalent to $([Bo^2B],A)$ which is \Draw since Bob already dominates.
 \end{itemize}

Assume now that $r=q=2$. As before, it is enough to give a strategy for $C$ restricted to the three paths.
\begin{itemize}
 \item  If $p \geq  5$, then Bob claims $a_3$, which enforces Alice to reply either $a_2$ or $a_4$ by Lemma~\ref{lem:P5}.
	  \begin{itemize}
	 \item If Alice claims $a_4$, then Bob successively claims  $b_1$, $c_1$ (with forced Alice to claim $b_2$ and $c_2$), and then  $a_1$, which creates two $A$-traps in $a_2$ and $c$.
 
	 \item If Alice replies by $a_2$, then, first, Bob claims, $a_5, a_7, \dots$, and so on until $a_{p}$ is dominated by Bob, replies of Alice being forced on $a_4, a_6, \dots$.  Second, Bob claims $c$. If Alice does not answer in the set $\{ a_1, b_1, b_2, c_1, c_2\} $, Bob successively claim $b_1$ creating a $A$-trap in $b_2$, $c_1$, creating a $A$-trap in $c_2$ and $a_1$ isolating $c$. Thus, Alice must claim  a vertex in the set  $\{ a_1, b_1, b_2, c_1, c_2\} $. Then Bob claims the unclaimed vertex of $[AoA]$ and dominates the whole position before Alice.  
	\end{itemize}
	
\item  If $p = 4$,   the position can be treated as for $p = 5$ if Alice replies in $a_2$ or $a_4$. But she can also claim $a_1$. In this case, Bob claims the unclaimed vertex of $[AoA]$  with the threat to claim $c$ and dominate before Alice. Even if Alice replies in $c$,  Bob succeeds in  dominating before Alice by successively claiming $b_1$, $c_1$ and $a_2$ (Alice will not dominate $a_3$ during this time). 
\item  If $p = 1$ or $p=3$, then Bob  claims $c$. If Alice claims the unclaimed vertex of $[AoA]$, $a_2$ or $a_3$, then Bob can claim $b_1$ and $c_1$, forcing Alice to reply by claiming $b_2$ and $c_2$. Then Bob can isolate $c$ by claiming $a_1$. Thus, Alice should answer
by claiming a vertex in $\{a_1,b_1, b_2, c_1, c_2\}$. Then Bob can claim the unclaimed vertex of $[AoA]$. If $p=1$, he wins. If $p=3$, he can dominate in one move by claiming either $a_2$ or $a_3$ whereas Alice need at least two moves to dominate. 
\item If $p = 2$, then  we have a fork (since all the other branches of $C$ starting from $c$ must have length 2 and the vertices adjacent to $C$ are labeled by $1$, and thus of degree 2), which is not possible, by hypothesis.  

\item  If $p = 0$, then Bob claims the unclaimed vertex of $[AoA]$. Alice needs at least two moves to dominate. If Alice does not claim $c$, Bob wins at his second turn by claiming it. If Alice claims $c$ she still need two moves to dominate. Then Bob can dominate before by claiming $b_1$ and $c_1$. \qedhere
  \end{itemize}
\end{proof}

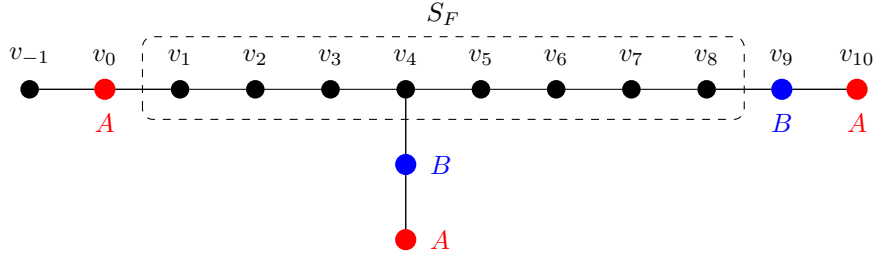
\begin{figure}
    \centering
\begin{tikzpicture}

\draw ( 0 , 0 ) node[v](s){} node[above=2mm] {$v_{-1}$};
\draw ( 1 , 0 ) node[R](s1){}  node[above=2mm] {$v_{0}$} node[below=2mm] {\color{red}$A$};
\draw ( 2 , 0 ) node[v](s2){}  node[above=2mm] {$v_{1}$};
\draw ( 3 , 0 ) node[v](s3){}  node[above=2mm] {$v_{2}$};
\draw ( 4 , 0 ) node[v](s4){}  node[above=2mm] {$v_{3}$};
\draw ( 5 , 0 ) node[v](s5){}  node[above=2mm] {$v_{4}$};
\draw ( 6 , 0 ) node[v](s6){}  node[above=2mm] {$v_{5}$};
\draw ( 7 , 0 ) node[v](s7){}  node[above=2mm] {$v_{6}$};
\draw ( 8 , 0 ) node[v](s8){}  node[above=2mm] {$v_{7}$};
\draw ( 9 , 0 ) node[v](s9){}  node[above=2mm] {$v_{8}$};
\draw ( 10 , 0 ) node[v](s10){}  node[above=2mm] {$v_{9}$};
\draw ( 11 , 0 ) node[v](s11){}  node[above=2mm] {$v_{10}$};
\draw ( 10 , 0 ) node[v, color = black!30](s10){};
\draw ( 11 , 0 ) node[v, color = black!30](s11){};

\draw ( 10 , 0 ) node[B](s10){} node[below=2mm] {\color{blue}$B$};
\draw ( 11 , 0 ) node[R](s11){} node[below=2mm] {\color{red}$A$};

\draw ( 5 , -1 ) node[v](j1){};
\draw ( 5 , -2 ) node[v](j2){};

\draw ( 5 , -1 ) node[v, color = black!30](j1){};
\draw ( 5 , -2 ) node[v, color= black!30](j2){};

\draw ( 5 , -1 ) node[B](j1){} node[right=2mm] {\color{blue}$B$};
\draw ( 5 , -2 ) node[R](j2){} node[right=2mm] {\color{red}$A$};

\draw[dashed,rounded corners] (1.5,-0.4) rectangle (9.5,0.7);
\node at (5.5,1) {$\sk_F$};

\path[draw] (s2) to (s9){};

\path[draw] (s) to (s1){};
\path[draw] (s1) to (s2){};
\path[draw] (s9) to (s10){};
\path[draw] (s10) to (s11){};
\path[draw] (s5) to (j1){};
\path[draw] (j2) to (j1){};

\path[draw] (s) to (s1){};
\path[draw] (s1) to (s2){};
\path[draw] (s9) to (s10){};
\path[draw] (s10) to (s11){};
\path[draw] (s5) to (j1){};
\path[draw] (j2) to (j1){};

%\draw ( 7 , 0 ) node[B](s6){};

%\draw ( 8 , 0 ) node[R](s8){};

%\draw ( 3 , 0 ) node[B](s3){};

%\draw ( 4 , 0 ) node[R](s4){};

%\draw ( 0 , 0 ) node[B](s){};

\end{tikzpicture}
    \caption{Example of a tree $F$ such that $o(F) = \oD$ but $\sk_F$ is a path. Alice should claim first $v_0$ or $v_{9}$. If Alice claims $v_0$, Bob claims the other vertices of $M_F$. Then Bob can win by claiming $v_2$ (Alice should answer $v_1$, $v_3$ or $v_4$), $v_6$(Alice should answer $v_5$,$v_7$ or $v_8$) and $v_{-1}$. The case where Alice claims first $v_9$ is similar.}
    \label{fig : draw path}
\end{figure}

Lemma~\ref{lem:fork_or_path} gives us the possible structures of the connected components of $\sk_F$  to have a position \Awin.  But actually, this condition is not sufficient: there are for example trees where $\sk_F$ is a path, but Bob can obtain a \Draw (see for example Figure~\ref{fig : draw path}). We need to consider  which vertices of $\sk_F$ are already dominated by Bob. 

\begin{definition}
Let $X, Y, \in \{A,B \}$, $n$ be a positive integer, and a subset $U \subseteq \{1, 2, \ldots , n\}$. We denote by $[X o^nY]^U$  the position obtained from the bounded path $[Xo^nY]$ where for each $i \in U$ a pendant edge $x_iy_i$ is added to the vertex $v_i$, with $x_i$ is linked to $v_i$ and claimed by Bob, and $y_i$ claimed by Alice.
\end{definition}

Informally, an  $[X o^nY]^U$  is a bounded path, where some vertices are already dominated by Bob.  
As an illustration, if we set that $v_{-1}$ is claimed by Bob on Figure~\ref{fig : draw path}, then we obtain the position $[Ao^8B]^{\{4\}}$.
When $U$ is empty,  $[Xo^nY]^U = [Xo^nY]$. %Note that adding $U$ can only by favorable to Bob.
If $Y=B$  (respectively $X=B$), one can assume that $n\notin U$ (resp. $1\notin U$). Indeed, $v_n$ (resp. $v_1$) is already dominated by Bob.

Finally, if $C$ is a connected component of $\sk_F$ that is a path connected to $v_0$ by a leaf, then $P_C$  is equivalent to a position $[Ao^nB]^U$, for a fixed $n$ and a fixed $U$. Indeed, set for $U$ all the integers $i$ such that $v_i$ is connected to a vertex of $M_F\setminus \{v_0\}$. 

The following three observations give natural properties about bounded paths. The first one is about the existence of a pairing.

\begin{observation} \label{obs:pairing}
A position $[X o^nY]^U$ contains a A-pairing, except if $X =Y= B$ and $n$ is odd. 
\end{observation}

Roughly speaking, the next two observations say that it is always better for Alice to play on a bounded path with the extremities claimed by her, and with fewer vertices dominated by Bob in $U$.

\begin{observation} \label{obs:paths}
For any integer $n$ and any set $U  \in \{1, ..., n\}$, and any $X, \in  \{A, B\}$, we have  
$$ o([A o^nA]^U,X) \geq  o([A o^nB]^U,X)  \geq  o([B o^nB]^U,X) . $$ 
\end{observation}

\begin{observation} \label{obs:include}
For any integer $n$ and any sets $U \subseteq U'  \subseteq \{1, ..., n\}$, and each $X, Y , Z \in  \{A, B\}$, we have  
$$ o([X o^nY]^U,Z) \geq  o([X o^nY]^{U'},Z) $$. 
\end{observation}

% In the following definition, we give all the winning components for Alice (they are proved to be winning positions in Theorem~\ref{maintree}.

% Some vertices in $\sk_F$ are already dominated by Bob. They correspond to the vertices that have a neighbor in $M_F\setminus \{v_0\}$. We denote the set of these vertices by $S_T$.\\

In the next definition, we define the components $C$ that are {\em favorable} for Alice, and among them, the ones that  are {\em strongly} and {\em weakly} favorable. Theorem~\ref{th:maintree} will justify this terminology.

\begin{definition}\label{def:fav}
%{\color{red}
%Let $T$ be a tree satisfying the conditions of Lemma~\ref{lem:firstmove} for a vertex $v_0\in M_F$. 

%Let $B_{T}$ be the set defined below. 
%$$B_{T}  =   \{v \in \sk_F, v \mbox{ is not a leaf of the tree } T [\sk_F \cup \{v_0\}], \exists v'  \in M_F \setminus \{v_0 \} \vert \{v, v' \} \in E \}. $$}
Let $v_0$ be the first move of Alice, satisfying Lemma~\ref{lem:firstmove}. Let $C$ be a connected component of $\sk_F$.
We say that $C$ is {\em favorable for Alice} if it satisfies one of the following case:
\begin{enumerate}
\item $C$ is a path connected to $v_0$ by a leaf, i.e, $P_C=  [A o^nB]^U$  with $ U \subseteq \{1, 2, \ldots, n-1\}$, and at least one of the following cases holds: 
\begin{enumerate}
\item  $n  \in \{1,2\}$ and $U = \emptyset$;
\item $n=3$ and $U\subseteq \{1\}$ or  $U\subseteq \{2\}$; 
\item $n\geq 4$ and $U \subseteq \{2,3,{n-2}\}$;
\item $n \geq 9 $, $n$  is odd and  $U \subseteq \{2,5, {n-2}\}$;
\item $n \in \{9,11\}$ and $\{3,5 \} \subseteq U \subseteq \{2,3,5,{n-2}\}$;

% \item  $n\leq 2$ and $U = \emptyset$;
% \item $n=3$ and $U\subseteq \{1\}$ or  $U\subseteq \{2\}$;
% \item $n \in \{4,5 \}$ and $U\subseteq \{2,3\}$;
% {\color{red} \item $n=6 $ and $U\subseteq \{2,3,4\}$;
%  \item $n =  7 $ and $U\subseteq \{2,3,5\}$;}
% \item $n\geq 8$ and $U \subseteq \{2,3,{n-2}\}$;
% \item $n \geq 9 $, $n$  is odd and $ \{5\}  \subseteq U \subseteq \{2,5, {n-2}\}$;
% \item $n \in \{9,11\}$ and $\{3,5 \} \subseteq U \subseteq \{2,3,5,{n-2}\}$.
\end{enumerate}
\item $C$ induces a fork and the only vertices of $C$ that can be connected to $M_F$ are the center $c$, the leaves except the one connected to $v_0$, and eventually the neighbor of $c$ between $c$ and $v_0$.
\end{enumerate} 

Moreover, if $C$ belongs to the cases (1.a), (1.b) or (1.c), we say that it is {\em strongly} favorable. On the opposite, if $C$ belongs to the cases (1.e) or (2), we say that it is {\em weakly} favorable.\footnote{If $C$ corresponds to the case (1.d), it is neither strongly nor weakly favorable.} 
\end{definition}

%{\color{blue}
%Ancienne classif
%\begin{itemize}
%\item (a)  $n\leq 2$ and $U = \emptyset$;
%\item (b) $n=3$ and $U\subseteq \{1\}$ or  $U\subseteq \{2\}$;
%\item (c) $n \in \{4,5 \}$ and $U\subseteq \{2,3\}$;
%{\color{red} \item (d) $n=6 $ and $U\subseteq \{2,3,4\}$;
 %\item (e) $n =  7 $ and $U\subseteq \{2,3,5\}$;}
%\item (f) $n\geq 8$ and $U \subseteq \{2,3,{n-2}\}$;
%\item (g) $n \geq 9 $, $n$  is odd and $ \{5\}  \subseteq U \subseteq \{2,5, {n-2}\}$;
%\item (h) $n \in \{9,11\}$ and $\{3,5 \} \subseteq U \subseteq \{2,3,5,{n-2}\}$.
%\end{itemize} 
%}

We can now state the final theorem that ends the characterization of trees. This theorem is actually technical to prove, and thus the proof is postponed to Section~\ref{sec:mainthm}. 

\begin{theorem}\label{th:maintree}
Let $F$ be a tree  and $v_0\in M_F$ be a first move of Alice that satisfies the condition of Lemma~\ref{lem:firstmove}.
The position $((F,\{v_0\},\emptyset),B)$  is \Awin if and only if all the components of $\sk_F$ are favorable to Alice and at most one of them is weakly favorable.\end{theorem}

\subsection{Proof of Theorem~\ref{thm:trees}}

\begin{figure}
\begin{center}
\scalebox{0.75}{\begin{tikzpicture}[node distance = 3cm, auto]
  \node (input) {Input $F$};
  \node[decision, below of=input, node distance=2cm, aspect=2](isolated) {Is there an isolated vertex $v_0$?};
   \node[sortieA, right of=isolated, node distance=14cm](Aiso) {$\oA$};
\node[sortieD, below of=Aiso, node distance=1.5cm](Diso) {$\oD$};
     \node[decision, left of=Aiso, node distance=7cm, aspect=2](PM) {Is there a perfect matching in $F\setminus \{v_0\}$};

\node[block, below of=isolated, node distance=2.5cm](rmvedge) {Remove isolated edges};

  \node[decision, below of=rmvedge, node distance=2.2cm, aspect=2](cherry) {How many cherries?};
  \node[sortieD, right of=cherry,node distance=14cm](D2ch) {$\oD$};
  \node[sortieA, below of=D2ch,node distance=2cm ](A1ch) {$\oA$};
  \node[decision, left of=A1ch,aspect=2, node distance=7cm](1cherry) {\small Is there a matching in $F\setminus \{c\}$ that dominates $V(F)\setminus N[c]$?};
  \node[sortieD, below of=A1ch,node distance=2cm](D1ch) {$\oD$};
  
  \node[block, below of=cherry, node distance =2.5cm](compute) {Computes $L_F$,$M_F$,$S_F$};

  \node[decision,below of=compute, node distance=2.8cm, aspect=2](empty) {Is $\sk_F$ empty?};
    \node[block, below of=empty, node distance= 2.4cm](rmv) {Remove components of $F$ with empty skeleton};
  
  \node[sortieA, right of=empty,node distance=14cm ](Aem) {$\oA$};
  
  \node[decision,below of=rmv, node distance=2.3cm,aspect =2](connected) {Is $\sk_F$ connected ?};
  
  \node[decision,below of=1cherry, node distance=8cm, aspect =2](star) {Is $\sk_F$ a star whose center is not linked to $M_F$?};
  \node[sortieA, right of=connected, node distance=14cm ](Astar) {$\oA$};
  \node[decision,below of=star, node distance=3cm,aspect=2](onecomp) {$\exists v_0\in M_F$ that makes $P_{\sk_F}$ favorable?};
  \node[sortieA, below of=Astar,node distance=3cm ](A1comp) {$\oA$};
  \node[sortieD, below of=A1comp,node distance=1.5cm ](D1comp) {$\oD$};

  \node[decision,below of=connected, node distance=6cm,aspect=2](v0nonconn) {$\exists v_0\in M_F$ that connects $\sk_F$?};
  \node[sortieD, below of=D1comp,node distance=1.5cm ](Ddisc) {$\oD$};

  %\node[block,below of=v0nonconn, node distance=3cm](v0) {$v_0$ be such vertex};
  \node[decision,below of=v0nonconn, node distance=3cm,aspect=2](allcc) {All cc of $\sk_F$ favorable to Alice ?};

  \node[decision,below of=onecomp, node distance=6cm,aspect=2](weak) {At most one cc weakly favorable ?};
  
  \node[sortieA, below of=Ddisc,node distance=3cm ](Aweak) {$\oA$};
  \node[sortieD, below of=Aweak,node distance=1.5cm ](Dweak) {$\oD$};
  \node[sortieD, below of=Dweak,node distance=1.5cm](Dnonfav) {$\oD$};

  \path[line] (star) -- node[pos=0.1] {yes} (Astar) ++ (-1.5,0.2) node {Lemma \ref{lem:star}};
  \path[line] (star) -- node {no} (onecomp);
  \path[line] (onecomp) -- node[pos=0.1] {yes} (A1comp) ++ (-1.5,0.2) node {Theorem\ref{th:maintree}};
  \path[line] (onecomp) |- node[pos=0.2] {no} (D1comp) ++ (-1.5,0.2) node {Theorem \ref{th:maintree}};
  \path[line] (v0nonconn) -- node[pos=0.1] {no} (Ddisc) ++ (-1.5,0.2) node {Lemma \ref{lem:firstmove}};
  \path[line] (v0nonconn) -- node {yes, fix $v_0$} (allcc);
  \path[line] (allcc) -- node {yes} (weak);
  \path[line] (weak) -- node[pos=0.1] {yes} (Aweak) ++ (-1.5,0.2) node {Theorem \ref{th:maintree}};
  \path[line] (weak) |- node[pos=0.3] {no} (Dweak) ++ (-1.5,0.2) node{Theorem \ref{th:maintree}};
  \path[line] (allcc) |- node[pos=0.1] {no} (Dnonfav) ++ (-1.5,0.2) node{Theorem \ref{th:maintree}};
  
  \path [line] (input) -- (isolated);
  \path [line] (isolated) -- node [near start] {no} (rmvedge);
 \path [line] (rmvedge) --  (cherry);
\path [line] (isolated) -- node[near start] {yes} (PM);
\path [line] (PM) -- node[near start] {yes} (Aiso) ++(-1.5,-0.8) node {Lemma \ref{lem:isolatedvertex}};
\path [line] (PM) |- node[near start] {no} (Diso);

\path (rmvedge) ++ (3,0) node {Lemma \ref{lem:isolatededge}};
\path (rmv) ++ (3,0) node {Lemma \ref{lem:rmvSvide}};
  \path [line] (cherry) -- node [pos=0.1] {$\geq 2$} (D2ch) ++ (-1.5,0.2) node {Lemma \ref{lem:2cherries}};
  \path [line] (cherry) -- (1cherry);
  \path (cherry) ++ (2,-1) node {1, of center $c$};
  \path [line] (1cherry) -- node[pos=0.1] {yes} (A1ch) ++ (-1.5,0.2) node {Lemma \ref{lem:1cherry}};
  \path [line] (1cherry) |- node [near start] {no} (D1ch) ++ (-1.5,0.2) node {Lemma \ref{lem:1cherry}};
  \path[line] (cherry) --node[near start] {0} (compute);
  \path[line] (compute) -- (empty);
 \path[line] (empty)--(rmv);
 \path[line] (empty) --node[pos=0.1] {yes} (Aem) ++ (-1.5,0.2) node {Lemma \ref{lem:skvide}};
 \path[line] (empty) -- node[near start] {no} (rmv);
 \path[line] (connected) -- node[near start] {no} (v0nonconn);
% \path[line] (v0nonconn) -- node[pos=0.1] {no} node[pos=0.8] {Lemma \ref{lem:conn}} (Ddisc);Input
 \path[line] (connected) -- node {yes} (star);
 \path[line] (rmv) -- (connected);
 
\end{tikzpicture}}
\end{center}
\caption{\label{fig:algo} The decision tree to compute the outcome of any forest.}
\end{figure}
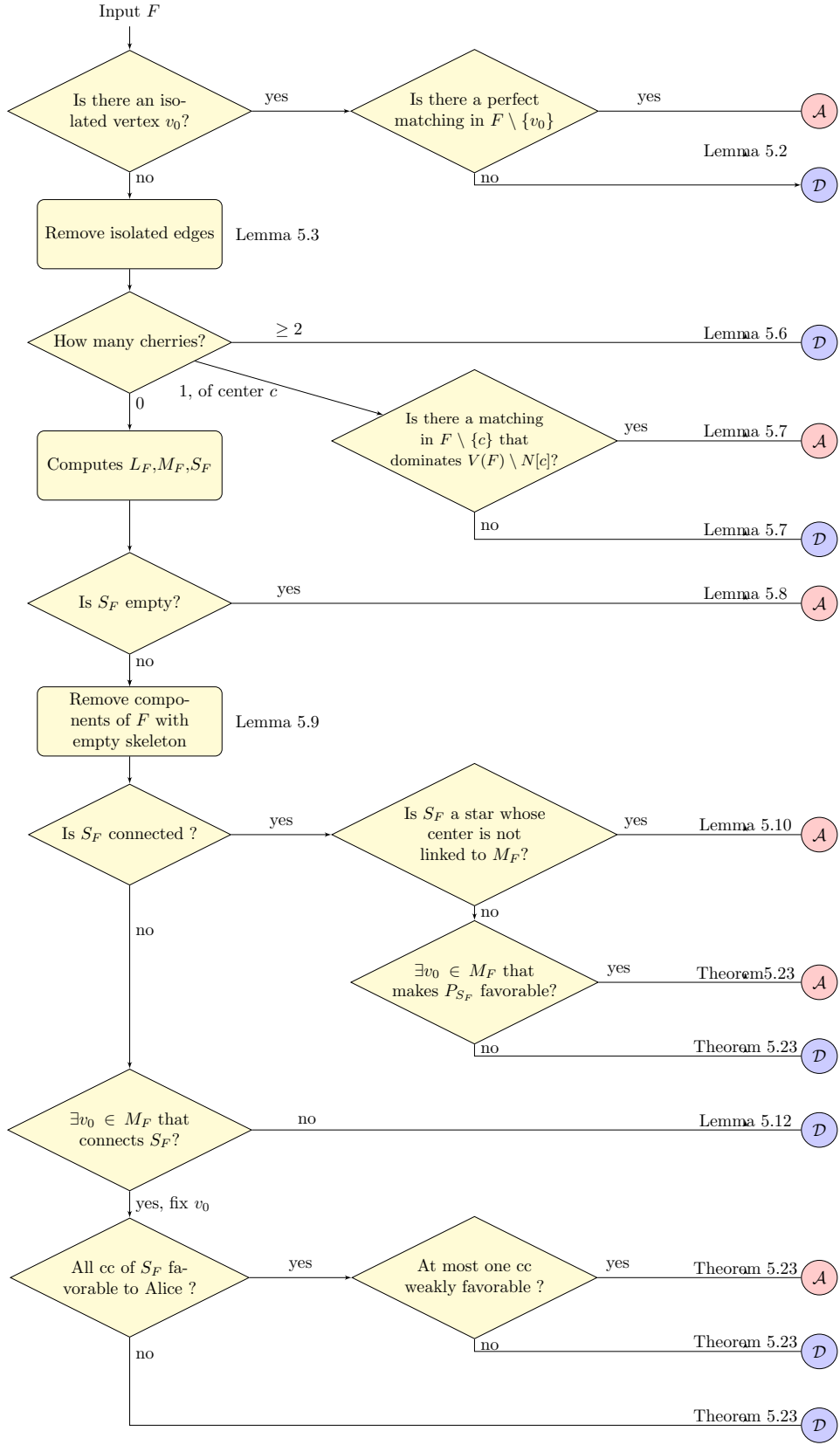

We now have all the ingredients to complete the study of forests and prove Theorem~\ref{thm:trees}. 
We summarize the algorithm into the diagram of Figure~\ref{fig:algo}.

\begin{proof}[Proof of Theorem~\ref{thm:trees}] Let $F$ be a forest. If $F$ has an isolated vertex $v_0$, then by Lemma~\ref{lem:isolatedvertex}, $o(F)=\oA$ if and only if $F\setminus \{v_0\}$ has a perfect matching.
If $F$ has an isolated edge $e=uv$, then by Lemma~\ref{lem:isolatededge}, $F$ has the same outcome as $F\setminus \{u,v\}$. Thus, we can assume that all the connected components of $F$ have at least three vertices. If $F$ has two cherries, then $o(F)=\oD$ by Lemma~\ref{lem:2cherries}. If it has one cherry $(c,l,l')$, then by Lemma~\ref{lem:1cherry}, $o(F)=\oA$ if and only if there is a matching in $F\setminus \{c\}$ that covers $V(F)\setminus N[c]$.
Thus, we can assume that $F$ has no cherry. If $\sk_F=\emptyset$, then $o(F)=\oA$ (Lemma~\ref{lem:skvide}). Otherwise, one can remove components of $F$ that have an empty skeleton (Lemma~\ref{lem:rmvSvide}).
If $\sk_F$ induces a star of center $c$ such that $c$ is not adjacent to $M_F$, then $o(F)=\oA$ (Lemma~\ref{lem:star}).
Otherwise, $F$ is standard. If $F$ is not connected, then $o(F)=\oD$.
If $\sk_F$ is not connected, Alice should claim the vertex $v_0$ of $M_F$ that connects $\sk_F$ (Lemma~\ref{lem:firstmove}) if it exists (otherwise $o(F)=\oD$). Then Alice wins if and only if all the components of $S_F$ are favorable to her and at most one is weakly favorable (Theorem~\ref{th:maintree}).
Assume now that $\sk_F$ is connected. Let $C$ be the tree induced by $\sk_F$. Alice should claim in $M_F$ and connected to a leaf of $C$ (Lemma~\ref{lem:firstmove}). If there exists $v_0\in M_F$ that is adjacent to a leaf of $C$ such that $P_C$ is favorable to Alice, then $o(F)=\oA$, otherwise $o(F)=\oD$ (Theorem~\ref{th:maintree}).

In terms of complexity, almost all the operations described in the algorithm are elementary (finding a matching in a forest, identifying the isolated vertices, edges and cherries, computing $L_F$, $M_F$ and $\sk_F$, deciding if $S_F$ is connected or a star) and can be done in linear time by examining the tree from the leaves. When $\sk_F$ is not connected, there exists at most one $v_0$ that can connect them and this is easy to check.
If $\sk_F$ is connected and thus reduce to a single component, one have to check all the possible $v_0$ that could make $P_{\sk_F}$ favorable. If $\sk_F$ induces a path, there are at most two possibilities for $v_0$ since it should be connected to an extremity of the path. If $\sk_F$ induces a fork (easy to check) of center $c$, to be favorable, there must be at most one neighbor of $c$ that is adjacent to $M_F$. If there is exactly one neighbor $a_1$ of $c$ adjacent to $M_F$, then the leaf of $\sk_F$ adjacent to $a_1$ must have only one neighbor in $M_F$ and this neighbor will be the vertex $v_0$. If all the neighbors at distance $1$ of $c$ in $\sk_F$ are not adjacent to vertices of $M_F$, then there must be at least one leaf of $\sk_F$ adjacent to exactly one vertex of $M_F$. This vertex of $M_F$ would be $v_0$. In all the other cases, $P_{\sk_F}$ cannot be favorable. Thus, finding $v_0$ can be done in linear time.
Finally, once $v_0$ is fixed, deciding whether a component is favorable or not can also be done in linear time (because it is a fork or a path with pending edges). 
% {\color{blue}
% \begin{itemize}
%     \item If $\sk_F$ is a fork, How many neighbors of the central node $c$  which are neighbors of with an element of $M_F$ ? 
%     \begin{itemize}
%         \item if there are more than 2, then $\oD$
%         \item  if there is a unique $v$, consider $w$ the other neighbor of $v$ in the fork. 
%          $d_G(w) = 2$ ?  if yes then $\oA$, if no then $\oD$
%         \item  if there is no such vertex,  is there a vertex at distance 2 from $c$ of degree 2 in $G$ ?  if yes then $\oA$, if no then $\oD$
%             \end{itemize}
%       \item If $\sk_F$ is a path, it suffices to  try a neighbor of each endpoint of the path as potential $v_0$.      
% \end{itemize}
% }
\end{proof}

\section{Proof of Theorem \ref{th:maintree}}\label{sec:mainthm}

We here prove Theorem~\ref{th:maintree} by considering successively both directions of the equivalence. The next subsection proves that the favorable components are necessary conditions. The sufficient condition is then proved in Subsection 6.2.

\subsection{The direct part}

In this subsection, we assume that $o((F,\{v_0\},\emptyset),B)=\oA$ and we prove that all the components are favorable and at most one is weakly favorable.
We already know by Lemma~\ref{lem:fork_or_path} that all the connected components of $\sk_F$ are paths or forks.
In Lemma~\ref{lem:fav_fork}, we prove that if a component is a fork, it should satisfy the conditions of Point 2 in Definition~\ref{def:fav}. In Lemma~\ref{lem:fav_path}, we prove that if a component is a path, it should satisfy the conditions of Point 1 in Definition~\ref{def:fav}. Finally, in Lemma~\ref{lem:twoweakly}, we prove that if two components are weakly favorable, then the outcome is $\oD$.

%%%%%%%%%%%%%%%%%%%%%%%%%%%%%%%%
%
%		LEMMES : WINNING   IMPLIES   FAVORABLE
%
%%%%%%%%%%%%%%%%%%%%%%%%%%%%%%%

\begin{lemma} \label{lem:fav_fork}
   Assume that  $o((F,\{v_0\},\emptyset),B) ) = \oA$. Let  $C$ be a connected component of $\sk_F$ which induces a fork. Then $C$ is favorable.
 \end{lemma}
 
 \begin{proof}
 Proceeding by contraposition, we 
assume that $C$ is not favorable.  We will prove that $o([A oA]  \cup P_C,B) = \oD$ which gives the result using Lemma~\ref{lem:component}.  

Let $c$  denotes the center of the fork. Let $a_1$ and $a_2$ be the two vertices on the path between $c$ and $v_0$ with $a_2$ adjacent to $v_0$, and  $v_{-1}$ be the unclaimed vertex of $[AoA]$. Since $C$ is not favorable, it means that either $a_2$ or a neighbor of $c$ but not $a_1$ is dominated by Bob in $P_C$. We consider the two cases.

\begin{itemize}
    \item
    If $a_2$ is already dominated by Bob, then Bob claims $c$. 
    \begin{itemize}
        \item If Alice replies in $v_{-1}$ then Bob can claim all the neighbors of $c$ distinct from $a_1$ and Alice is forced to answer the leaf of $C$ adjacent to  the vertex claimed by Bob.  Afterward,  Bob claims $a_1$, isolating the vertex $c$, which ensures that   Alice will not dominate the graph. 
        \item otherwise,  Bob claims $v_{-1}$ and, therefore,  dominates the graph before Alice. 
    \end{itemize}
    \item If Bob already dominates a neighbor $b_1$ of $c$ (with $b_1\neq a_1$), then Bob claims all the other neighbors of $c$ distinct from $a_1$ and $b_1$, forcing Alice to claim the leaves adjacent to them. Then Bob claims the unclaimed vertex of $[AoA]$. At this point, he needs only one move to dominate the graph that can be either $a_1$ or $a_2$, whereas Alice still need to dominate $a_1$ and $b_2$, which cannot be dominated in a single move. Thus, Bob will dominate before Alice. \qedhere
    \end{itemize}
%OLD PROOF
% If $a_2$ is adjacent to a vertex of $M_F$, 
%
% $a_1, b_1, c_1$ denote the neighbors of $c$, in such a way that such that $x$ is between $c$ and $v_0$, and $x_1, y_1 z_1$ respectively denote the neighbors of $x, y,  z$ different from $c$. 
% Some vertices in $\sk_F$ are already dominated by Bob. They correspond to the vertices that have a neighbor in $M_F\setminus \{v_0\}$. We denote the set of these vertices by $S_T$.
%
% Since  $C$ is not favorable,   $ S_T \cap  \{ x_1, y, z\} \neq \emptyset$. Up to symmetry and non triviality, we have to study the two following cases 
%
% \begin{itemize}
%
% \item If $x_1 \in S_T$, then Bob claims $c$ and, next,  achieves to dominate before Alice by claiming $v_{-1}$, 
%
% \item  If $y \in S_T$, then  Bob claims $z$, this forces Alice to play $z_1$. At this time, Bob claims $v_{-1}$ and needs one more move to dominate while  Alice needs two, as she cannot dominate both $x$ and $y_1$ in one move. 
% \end{itemize}
 \end{proof}

%%%%%%%%%%%%%%%%%%%%%%%%%%%%%
%
%   PATH ETUDE DRAW
%
%%%%%%%%%%%%%%%%%%%%%%%%%%%%%%%%%%%

\begin{lemma} \label{lem:fav_path}
Assume that  $o((F,\{v_0\},\emptyset,B) ) =   \oA$,  and let  $C$ be a connected component of $\sk_F$ which induces a path. Then $C$ is favorable.
 \end{lemma}
 
\begin{proof}
Let $P_C = [Ao^nB]^U$. We proceed by contraposition. Assume that $C$ is not favorable.  We will prove that $o([A oA]   \cup P_C,B) = \oD$, which gives the result by using Lemma~\ref{lem:component}. We split the proof  into a few claims. We denote by $v_{-1}$ the unclaimed vertex of $[AoA]$ and by $v_1,...,v_n$ the unclaimed vertices of $P_C$. We consider for the labeling that $P_C$ is rooted in $v_0$.

\begin{claim} \label{slem:draw_even}
If $n$ is even and $5 \in  U$, then $o([A oA]   \cup P_C,B) = \oD$.
 \end{claim}

\begin{proofclaim}   
  When $n= 6$,   Bob  claims $v_3$, which enforces  Alice  to   reply in  $v_4$  (if Alice replies in $v_2$,  then Bob claims $v_5$ and create a double trap). Then 
Bob  successively claims  $v_{1}$ and $v_{-1}$ and, therefore dominates before Alice. 

If $n \geq 8$, $v_7$ is labeled by $1$. By Lemma~\ref{lem:cut}, the position can be reduced to the case $n=6$, and thus is also \Draw.
\end{proofclaim}

\begin{claim} \label{slem:key_set}
If $n \geq 4$ and if there exists $i\in U$ such that  $i \notin  \{2, 3, 5,  n- 2 \}$, then $o([A oA]   \cup P_C,B) = \oD$.
\end{claim}

\begin{proofclaim} 
Assume first that $i=1$.

\begin{itemize}
\item If $n \in\{ 4, 5\}$, then Bob first claims $v_3$, (which forces Alice to answer $v_4$ or $v_2$). Then   Bob claims $v_{-1}$ and dominates before Alice (who does not dominate $v_2$ or $v_4$). 
%\item  If  $n = 5 $, then Bob successively claims $v_3$ and  $v_{-1}$ and dominates before Alice. 
\item  If $n \geq 6$,  either  $v_5$ or $v_6$ is labelled by 1. Then as before, we can apply Lemma~\ref{lem:cut} to reduce to the case $n\in \{4,5\}$.
\end{itemize} 

We can now assume that $i \geq 4$. There  
are  two cases, according to the parity of the value   $n-i$. In each case,  Lemma~\ref{lem:cut}  is used.  
\begin{itemize}
\item If $n-i$ is odd,

\begin{itemize}
\item If $n-i = 1$, that is $i=n-1$, then  Bob claims  $v_{n-3}$ (this is possible since  $4 \leq i < n$), 

\begin{itemize}
\item For $n \geq 6$, Alice  necessarily replies in  $v_{n-2 }$. Indeed, by Lemma~\ref{lem:P5}, Alice should answer in $v_{n-4}$ or $v_{n-2}$, but if she claims $v_{n-4}$ then a component $[Bo^{2k+1}B]$ will appear which is \Draw by Lemma~\ref{lem:[BB]_odd}.
After Alice has claimed $v_{n-2}$, the resulting position $P$ is then equivalent to a position $ [AoA] \cup [A o^{n-4}B] ^{U'} \cup [Ao^2B]^{\{1\}}$, for some set $U'$. We thus have

$$o(P, B) \leq o( [AoA] \cup ([A o^{n-4}B] \cup [Ao^2B]^{\{1\}}), B) =
o( [AoA]  \cup [Ao^2B]^{\{1\}}, B) = \oD. $$
The first equality comes from Lemma~\ref{lem:[A, B]}  and the second one is obvious since Bob can dominate in one move.  This gives the desired result in this case.

 \item For $n  =  5$, Alice  can also reply in  $v_{n-1}$. In such a situation, Bob  claims   $v_{-1}$ and directly dominates before  Alice. If Alice claims $v_{n-2}$ we are in the same situation as the previous case.
 
  \end{itemize} 
\item  If  $n-i \geq 3$,  then  $v_{i+2}$ is labeled by $1$ and Lemma~\ref{lem:cut}  applies to reduce the instance to the case where $i = n-1$. 
\end{itemize}

\item  Assume now that $n-i$ is even. Since $i\notin \{n-2,n\}$, we have $n-i\geq 4$. If $n-i  \geq 6$,  then $v_{i+5}$ is labeled by $1$,  and Lemma~\ref{lem:cut}  applies to reduce the instance to the case where $n-i = 4$. Thus, we can assume that $i=n-4$. Furthermore, since $i\notin \{1,2,3,5\}$, we can assume that $n\geq 8$ and $n\neq 9$.
\begin{itemize}
\item If $ n \geq 12$, Bob  claims  $v_{n-9}$, which, as before, enforces Alice  to reply in $v_{n-8}$ to avoid a component $[Bo^{2k+1}B]$.
Afterward, Bob successively claims   $v_{n-11}, v_{n-13}, \dots$  and so on until $v_1$ is dominated by Bob, the answers of Alice being forced. 

Bob   now  claims  $v_{n-6}$, which enforces Alice to reply $v_{n-7}$, $v_{n-5}$ or $v_{n-4}$. Then  Bob claims  $v_{n-2}$,  which enforces Alice  to reply  $v_{n-1}$, or $v_{n-3}$. 
In the resulting position,  Bob dominates each vertex $i$ , $1 \leq i \leq n$, while  Alice does not. Indeed, there exists a vertex $i \geq n -6$,  which is not dominated by Alice, as she has claimed only two moves that can dominate vertices from $v_{i}$ for $n-6 \le i \le n$, and each of them dominates at most three vertices. Finally, Bob  claims   $v_{-1}$ and  dominates before  Alice.

\item if $n = 11$,  Bob claims $v_9$. By Lemma~\ref{lem:P5}, Alice should answer in $\{v_8, v_{10}\}$.  If Alice claims $v_8$, Bob claims $v_5$. Alice has to claim either $v_4$ or $v_6$. Then, Bob claims $v_{-1}$. Alice now needs to claim at least two more vertices to dominate: one to dominate $\{v_4,v_6\}$ and one to dominate $v_{10}$, while Bob will dominate with his next claim in $\{v_1, v_2\}$.
    
    If Alice claims $v_{10}$, Bob also claims $v_5$. Alice has to claim $v_6$, otherwise $v_7$ creates a double trap. Now Bob claims $v_{-1}$ and will dominate with his next claim in $\{v_1, v_2\}$, while Alice needs two vertices to dominate.

\item if  $n = 10$ , then Bob  successively claims  $v_8$,  $v_4$. As, by Lemma~\ref{lem:P5}, Alice has to claim a vertex adjacent to the vertex that Bob has claimed, her moves are almost forced. If she has claimed $v_9$ and $v_3$, Bob claims $v_6$ and creates two traps in $v_5$ and $v_7$. Otherwise, he claims $v_{-1}$ and dominates before Alice by a final claim in $\{v_1, v_2\}$.

\item if  $n = 8$, then Bob  successively claims  $v_6$,  $v_2$,   $v_{-1}$ and dominates before Alice. Indeed, she has only claimed two moves to dominate vertices between $v_2$ and $v_7$, and the only dominating set of size two on this path contains $v_6$.
\end{itemize}
\end{itemize}
\end{proofclaim}

\begin{claim} \label{slem:short}
Assume that $n\geq 13$, $n$ is odd and $\{3,5\} \subseteq U$, then $o([A oA]  \cup P_C, B) = \oD$.
\end{claim}

 \begin{proofclaim}
It suffices to prove it for $n = 13$, since,  for $n > 13$,
 Lemma~\ref{lem:cut} applies. 
 Bob claims $v_{11}$, which enforces  Alice to reply in $v_{10}$  (if Alice   replies in $v_{12}$, then Claim~\ref{slem:draw_even}  applies). Then,  Bob claims $v_{8}$. 
 
 \begin{itemize}
 
\item If Alice  replies in $v_7$, then Bob successively claims $v_5$  and $v_{-1}$. Afterward, Bob
finally  succeeds in dominating before Alice by claiming either $v_1$ or $v_2$. 

\item If Alice replies in $v_6$, then Bob claims $v_2$, which enforces Alice to claim $v_{1}$,$v_{3}$, or $v_{4}$.
Then Bob claims  $v_{-1}$. At this step,   
Bob threatens to claims $v_5$  and dominate before Alice. Thus,  the reply of Alice is necessarily $v_5$. 
Then, Bob claims $v_7$, which enforces the reply $v_9$ from Alice. 
Finally,  Bob claims either $v_3$ or $v_4$ (one of these vertices is free)  and dominates before Alice. 
\end{itemize}
\end{proofclaim}

We can now finish the proof of the lemma. Since $C$ is not favorable, we can split the cases according to $n$ as follows:

 \begin{itemize}
     \item If $n = 1$, then the result is obvious, since  there is no unfavorable component. 
   \item If $n = 2$  and $1 \in U$, then Bob claims $v_{-1}$ and then  dominates the graph before Alice, which is a contradiction.
   \item If $n = 3$  and $1, 2 \in U$, then Bob claims $v_{-1}$ and then  dominates the graph before Alice,  which is a contradiction
   \item  If  $n  \in \{4, 5, 6, 7\} $  then Claim~\ref{slem:draw_even} and Claim~\ref{slem:key_set} give the result. 

   \item  If $n \geq 8$ and $n$ is even, the combination of Claims~\ref{slem:draw_even} and~\ref{slem:key_set} gives that $U \subseteq \{2, 3, n-2\}$.

   \item If $n$ is odd and $n \geq 13$, then   the combination of Claims~\ref{slem:key_set}     and~\ref{slem:short}  gives that  either 
$U \subseteq  \{2, 3, n-2\}$ or $U\subseteq\{2, 5, n-2\}$, which gives the result .

 \item If $n \in \{9, 11\}$, then Claim~\ref{slem:key_set}    allows to conclude. \qedhere
 \end{itemize}
\end{proof}

%%%%%%%%%%%%%%%%%%%%%%%%%%%%%%%%
%
%    LEMME  : TWO WEAKLY favourable
%
%%%%%%%%%%%%%%%%%%%%%%%%%%%%%%%%%%%

\begin{lemma}\label{lem:twoweakly}
Assume there are two connected components  $C, C'$ of $\sk_F$ that are weakly favorable. 

Then $o(F,\{v_0\},\emptyset,B) ) = \oD$. 
 \end{lemma}

 \begin{proof} 
Assume first that both $C$ and $C'$  both induce forks. We prove the result for forks with exactly three branches. Indeed, if Bob has a strategy in this case, he will have a strategy for forks with more branches using Lemma~\ref{lem:cut} since all the neighbors of $c$ are labeled by $1$. 
Let $c$  denote the center of the fork induced by $C$, $a_1, b_1,  c_1$ denote the neighbors of $c$, in such a way that such that $a_1$ is between $c$ and $v_0$, and $a_2, b_2, c_2$ respectively denote the neighbors of $a_1, b_1,  c_1$ different from $c$ . We define in the same way $c', a_1', b_1', c_1', a_2', b_2'$ and $c_2'$ for $C'$. Bob start by claiming  $v_{-1}$. By symmetry, we can suppose that Alice replies one vertex among $\{c, a_1, a_2, b_1, b_2\}$.

\begin{itemize}
\item  If  Alice replies in $c$ ,  then Bob successively claims $b_1$ and $c_1$  (the replies of Alice are forced). Now Bob claims $c'$, then Bob needs two more moves to dominate (one in $(a_1,a_2)$ and one in $(a_1',a_2')$) whereas Alice needs three moves.
\item   If   Alice replies in  $a_1$, then   Bob claims $a_2$. \begin{itemize}
 
\item  If Alice replies  in $a_1'$, then Bob claims $a'_2$. 
At this step, Alice needs at least four moves  to dominate  while  Bob can dominate in three moves in a lot of manners, with a center of a fork, and two  vertices of the other fork.
Alice cannot avoid  Bob  to dominate using three moves, and therefore, before Alice.

\item If Alice replies in $c$, then Bob  successively claims $b_1$ and $c_1$.  After the forced replies of Alice, Bob claims in  $c'$ and dominates before Alice with one last move in $(a_1', a_2')$ while Alice cannot dominate in one move.  The case where    Alice replies in $c'$  can be treated symmetrically.

\item If Alice replies in $b_1$, then Bob claims successively  $c_1$  and $c'$   and dominates before Alice by claiming one vertex in $(c,b_2)$.

\item If Alice replies in $b'_1$, then Bob claims $c'$. At this step, Alice needs at least four moves   to dominate  while  Bob can dominate in three moves by a pairing strategy with the pairs $(a_1', a_2'), (b_1, b_2)$ and $(c_1, c_2)$. %Alice cannot avoid  Bob  to dominates using three moves. 

\end{itemize}

\item   If   Alice replies in  $a_2$, then   Bob claims $a_1$, and afterward, all is similar to the previous case.

\item   If   Alice replies in  $b_1$ or $b_2$, then Bob claims $c'$. Bob can now dominate in three moves by pairing $(a_1',a_2'), (a_1, a_2)$ and $(c, c_1)$ while Alice needs at least four moves.

 \end{itemize}

Assume now that both $C$ and $C'$ induce paths. Since they satisfy the conditions $1.e$ of Definition~\ref{def:fav}, their length is 9 or 11. Bob first reduce the paths to length 9 if needed so that both paths have length~$9$. Let $(v_1, v_2, ...., v_9)$ (respectively $(v'_1, v'_2, ...., v'_9)$) the path induced by $C$ (resp. $C'$), with $v_1$ and $v'_1$ connected to $v_0$. Since the paths are weakly favorable, Bob already dominates $v_3$, $v_5$, $v'_3$ and $v'_5$.

First, Bob successively claims $v_7$ and $v'_7$. By Lemma~\ref{lem:P5}, Alice should answer first $v_6$ or $v_8$ and then $v'_6$ or $v'_8$. We have three cases according to the replies of Alice.
\begin{itemize}

\item If the replies are $v_8$ and $v'_8$, then Bob continues by claiming $v_3$ and $v'_3$ which enforces Alice to 
reply in $v_4$ and $v'_4$ (if, for instance, Alice does not reply $v_4$ , then Bob claims $v_5$ and forbids Alice to dominate at the following claim). Afterward,  Bob claims  $v_{-1}$, and achieves to dominate before Alice as he only needs one move in $(v_1, v_2)$ and one in $(v_1', v_2')$ while Alice needs to dominate $v_2, v_2', v_6$ and $v_6'$, each of them requiring a different move.

\item If the replies are $v_6$ and $v'_8$ (note that $v_8$ and $v'_6$ is symmetric), then Bob continues by claiming $v'_3$ which enforces Alice to 
reply in $v'_4$, Afterward,  Bob claims  $v_{-1}$. Now Bob has a paired dominating set of size 3 $\{(v_1,v_2), (v_1', v_2'), (v_3, v_4)\}$, but Alice needs to dominate $v_2, v_8, v_2'$ and $v_6'$, each of them requiring a different move. Thus, Bob can dominate first.

\item If the replies are $v_6$ and $v'_6$,  Bob claims  $v_{-1}$.  At this time Bob and Alice need four moves to dominate, but each set of four moves for Alice contains the pair $\{v_3, v'_3 \}$.  Thus,  after the reply of Alice, Bob can claim one element of  the pair $\{v_3, v'_3 \}$ and, by this way, dominate before Alice. 
\end{itemize}

Assume now that $C$ induces a path and $C'$ induces a fork. As before, $C'$ can be assumed to have exactly three branches and $C$ nine vertices. We  denote the vertices of the fork and the path as in the previous cases. 

First, Bob claim $v_7$. Once again, Alice has to claim either $v_6$ or $v_8$, otherwise Bob ensure a draw by Lemma~\ref{lem:P5} There are two cases according to the reply of Alice.

\begin{itemize}

\item If Alice replies $v_8$, then Bob successively claims $b_1, c_1$ and $a_1$ (replies of Alice being forced in $b_2, c_2$ and $c$). Afterward, Bob claims $v_{3}$. Alice has to claim either $v_2$ or $v_5$ by Lemma~\ref{lem:P5}. If she claims $v_2$, Bob claims $v_5$ and creates two traps. If she claims $v_4$, Bob claims $v_{-1}$. Now Bob only needs one move in either $v_1$ or $v_2$ to dominate while Alice needs at least two, as she does not dominate $v_2$ nor $v_6$. 

\item If Alice replies $v_6$,  then Bob claims $v_{-1}$

\begin{itemize}

\item If Alice replies $c$, then Bob  successively claims  $b_1$, $c_1$ (replies of Alice are again forced in $b_2$ and $c_2$) and then  $v_3$. At this time, Bob needs two moves to dominate (one in $v_1, v_2$ and one in $(a_1, a_2)$), while  Alice needs to dominate $v_2$, $v_4$ and $v_8$, and each of them requires a different move. Therefore, Bob dominates before Alice by this way.

\item If Alice replies $v_3$, then Bob claims $c$. Now $(a_1, a_2), (v_1, v_2)$ and $(v_4,v_5)$ is a paired dominating set for Bob of size three, while Alice needs to dominated $a_2, b_1, c_1$ and $v_8$, each of these vertices requiring a different move.

\item If Alice replies $v_1$ or $v_2$, by Lemma~\ref{lem:sommetdomine}, as $N[v_1]\setminus N[V_t] \subset N[v_2]\setminus N[V_t]$ for $t \in \{A,B\}$, we can suppose, she plays $v_2$. Then Bob  successively claims  $b_1$, $c_1$, $a_1$  (replies are forced in $b_2, c_2, c$ respectively) and then $v_1$. At this time Bob needs one  move to dominate in $(v_3, v_4)$, while Alice needs to dominate $v_4$ and $v_8$ which she cannot dominate in a single move.

\item If Alice replies $a_1$ or $a_2$, by Lemma~\ref{lem:sommetdomine}, as $N[a_2]\setminus N[V_t] \subset N[a_1]\setminus N[V_t]$ for $t \in \{A,B\}$, we can suppose she claims $a_1$. Bob claims $a_2$.  
\begin{itemize}

\item If Alice replies $c$,   then Bob  successively claims  $b_1$, $c_1$, (replies are forced) and  $v_3$. 
 At this time Bob can dominate in two moves, one in $(v_1, v_2)$ and one in $(v_4,v_5)$ while Alice cannot, as she still has to dominate $v_2, v_4$ and $v_8$.

\item If Alice replies $v_1$ or $v_2$, by Lemma~\ref{lem:sommetdomine}, as $N[v_1]\setminus N[V_t] \subset N[v_2]\setminus N[V_t]$ for $t \in \{A,B\}$, we can suppose she replies $v_2$. Then Bob  successively claims  $b_1$, $c_1$, (replies are forced for Alice) and  $v_1$. At this time Bob needs one move in $(v_3, v_4)$ to dominate, while Alice needs two moves to dominate $v_4$ and $v_8$.

\item If Alice replies elsewhere,   then Bob claims  $c$. At this time Bob needs two  moves to dominate,one in $(v_1, v_2)$ and one in $(v_3, v_4, v_5)$. By hypothesis, at least one will be available. Alice needs at least  three moves to dominate, as she does not dominate at least three of the four vertices $v_3, v_8, b_1$ and $c_1$, each of them requiring a different move.
\end{itemize} 
\item If Alice replies elsewhere, then Bob claims  $c$. At this time Bob needs three  moves to dominate,one in $(v_1, v_2)$, one in $(a_1, a_2)$ and one in $(v_3, v_4, v_5)$. By hypothesis, at least one will be available in each of these sets at any moment. Alice needs at least  four moves to dominate, as she does not dominate at least four of the five vertices $a_1, b_1, c_1,v_3$ and $ v_8$, each of them requiring a different move. \qedhere
\end{itemize} 
\end{itemize} 
\end{proof}

\subsection{The converse part}
  
  In this final subsection, we prove the reverse part of Theorem~\ref{th:maintree}: if all the connected components of $\sk_F$ are favorable and at most one is weakly favorable, then Alice has a winning strategy. This part is naturally harder than the other one, as Alice cannot force Bob to answer where she would like to. Hence, all the possible answers of Bob must be considered, which was not the case previously, as Alice was often forced to answer locally to a move of Bob. 

  We first prove in Subsection~\ref{sec:strongly}, that we can remove the strongly favorable components (i.e. corresponding to the cases 1.a to 1.c of Definition~\ref{def:fav}). In Subsection~\ref{sec:weakly}, we consider only one weakly favorable component and give a strategy for Alice in this case.
  Finally, in Subsection~\ref{sec:stable}, we define a class of positions $\mathcal C_1$ that contains the starting positions (without strongly favorable components) and for which Alice can ensure either to win or to stay in this class after a move of Bob. By induction, this will imply that Alice has a winning strategy in this class.
  Note that one can always consider $U$ to be maximal in Definition~\ref{def:fav}. Indeed, if Alice have a strategy for $U$ maximal, she will have a strategy with any subset of $U$.

\subsubsection{Removing strongly favorable components }\label{sec:strongly}

\begin{lemma}  \label{lem:[A, B]+}
Let $Q$ be a position where Bob is not dominating. Consider a strongly favorable component $[Ao^nB]^U$. 
We have $$ o(Q, B) \leq  o(Q \cup [A o^nB]^U, B) .$$ 
\end{lemma}

In other words, adding a strongly favorable component can only be favorable to Alice.

\begin{proof} First note that since $[Ao^nB]^U$ is strongly favorable, we have $U\subseteq \{2,3,n-2\}$.

We prove by induction of the number $p$ of unclaimed vertices of $Q \cup [A o^nB]^U$ that if  $o(Q, B) = \oA$, then $o(Q \cup [A o^nB]^U, B) = \oA$. Note that $p\geq n$
Let $v_1,v_2, ....v_n$ denotes the sequence  of unclaimed vertices of $[A o^nB]^U$. 
First we can assume that $n>1$. Indeed, if $n = 1$, then $U=\emptyset$ and by Observation~\ref{obs:removedomcomp}, $o(Q \cup [AoB]^U, B) =  o(Q, B)$.\\

Thus, we can assume that $p\geq 2$ and $n\geq 2$. If $p=n = 2$, $Q$ contains no unclaimed vertices. Since $o(Q, B) = \oA$, Alice dominates $Q$ while Bob does not.  Thus, in $Q \cup [A o^2B]^U$  Bob claims $v_1$ or $v_2$, and Alice answers by claiming the other one. By this way, Alice dominates $Q \cup [A o^nB]^U$ while Bob does not. Thus, the resulting position is \Awin. \\

Consider now that $p\geq 3$. Let $y$ be the vertex claimed by Bob.
Assume first that $y$ is an unclaimed vertex of $Q$. Note that Bob cannot dominate $Q$ in one move unless Alice already does, otherwise we would have $o(Q,B)=\oD$.
 	\begin{itemize}
	\item If  $Q$ is not  dominated by  Alice yet, then Alice claims $x$ according to a winning strategy in $Q$.   Thus, by definition,  $o(Q_{x, y}, B) = \oA$ and by  induction hypothesis $o(Q_{x, y} \cup [A o^nB]^U, B)= \oA$.

	\item  if  $G$ is already dominated by  Alice:
  		\begin{itemize}
	\item if $2 \leq n \leq 4$,  then Alice claims $v_3$ and dominates the whole graph before Bob;

        \item  if $n = 5$, the Alice claims $v_3$. We have: 
$$o((Q \cup [A o^nB]^U)_{v_3, y}, B) \geq  o([A o^nB]^U)_{v_3, y}, B)  \geq    o( [Ao^2A]^{\{{2}\}}  \cup [Ao^{2}B] ), B)= \oA, $$.  where the first inequality comes from the fact that $Q$ is dominated by Alice, the second inequality comes from Lemma~\ref{lem:split} and from the fact that $U\subseteq \{2,3\}$, which is strongly favorable according to Definition~\ref{def:fav}.1.a,  and the final equality is obtained by induction hypothesis, using $Q' = [Ao^2A]^{\{{2}\}}$.
        \item if $n=6$, then Alice  claims $v_3$.  We have: 
$$o((Q \cup [A o^nB]^U)_{v_3, y}, B) \geq  o([A o^nB]^U)_{v_3, y}, B)  \geq    o( [Ao^2A]^{\{{2}\}}  \cup [Ao^{3}B] )^{\{{1}\}}, B)= \oA, $$ 
where the first inequality comes from the fact that $Q$ is dominated by Alice, the second inequality comes from Lemma~\ref{lem:split} and from the fact that $U\subseteq \{2,3,n-2\}$, which is strongly favorable according to Definition~\ref{def:fav}.1.b,  and the final equality is obtained by induction hypothesis, using $Q' = [Ao^2A]^{\{{2}\}}$. 
        
		\item if  $n \geq  7$, then Alice  claims $v_3$.  We have: 
$$o((Q \cup [A o^nB]^U)_{v_3, y}, B) \geq  o([A o^nB]^U_{v_3, y}, B)  \geq    o( [Ao^2A]^{\{{2}\}}  \cup [Ao^{n-3}B]^{\{{n -5}\}}, B)= \oA, $$ 
where the first inequality comes from the fact that $Q$ is dominated by Alice, the second inequality comes from Lemma~\ref{lem:split} and from the fact that $U\subseteq \{2,3,n-2\}$, which is strongly favorable according to Definition~\ref{def:fav}.1.c,  and the final equality is obtained by induction hypothesis, using $Q' = [Ao^2A]^{\{{2}\}}$ . 
		\end{itemize}
	\end{itemize}

Assume now that Bob claimed a vertex $v_i$ of $[Ao^nB]^U$.
	\begin{itemize}
		\item If $2 \leq  i \leq n-1$,  then Alice  claims $v_{i+1}$. We  have 
		$$o((Q\cup  [A o^nB]^U)_{v_i, v_{i+1}}, B)  = o(Q \cup [Ao^{i-1}B]^{\{2, 3\}} \cup [Ao^{n-i-1}B]^{\{{n-i-3}\}}, B) = \oA,$$
		where the first equality comes from Lemma~\ref{lem:split} and the inequality comes from two applications of the induction hypothesis.    
	
		\item  If $i=1$,  and $n \geq 3$  then Alice  claims $v_3$  . We have 
		$$o((Q\cup [A o^nB]^U)_{v_3, v_{1}},B)  \geq o(Q \cup  [Ao^{n-3}B]^{\{{n-5}\}}, B) = \oA, $$  by induction hypothesis. 
  
  If $n= 2$ then Alice  claims $v_2$ . We have 
		$o((Q\cup [A o^2B]^U)_{v_2, v_{1}},B)  = o(Q , B)= \oA$, from Observation~\ref{obs:removedomcomp}
		\item  If $i = n$,  then Alice  claims $v_{n-1}$.
		We have 
		$$o((Q\cup [A o^nB]^U)_{v_{n-1}, v_{n}}, B)  \geq o(Q \cup  [Ao^{n-2}B]^{\{2, 3\}}, B) =\oA, $$  by induction hypothesis.  
	\end{itemize}

Thus, for each claim $y$ of Bob  on $Q \cup [Ao^{n}B]^U$, there exists an answer  $x$ of Alice such that   $o((Q\cup [A o^nB]^U)_{x, y},B)= \oA$.  This ensures that $o(Q \cup [A o^nB]^U,B) = \oA$. 
\end{proof}

%\begin{remark}  
%From Lemma~\ref{lem:[A, B]+}, 
%in order to prove Theorem~\ref{th:maintree}, it suffices to prove is for trees such that the set of components   $C$ is such  that  at  most one component is weakly   favorable and all the other ones are in case 1.d. 
%\end{remark}

 % From this step,  our goal is to prove that $o(T, B) = \oA$ for each position $T$ satisfying  hypotheses  described the previous remark. We need to study some larger class of positions, in order to take in account the some positions  which can possibly appear during the execution of the game.  But, to do it  it, we first need some  lemmas, below. 

%%%%%%%%%%%%%%%%%%%%%%%%%%%%%%%%%%
%
%  SINGLE COMPONENT
%
%%%%%%%%%%%%%%%%%%%%%%%%%%%%%%%%%%%%

\subsubsection{Dealing with the weakly favorable component}\label{sec:weakly}

\begin{lemma} \label{lem:weakly}
Let $C$ be a connected component of $\sk_F$ that is  weakly favorable, then   
$o(P_C \cup [AoA], B)  = \oA$. Moreover, Alice can play to ensure that after each of her move, there is a $A$-pairing.
\end{lemma}

\begin{proof} 
Recall that by Observation~\ref{obs:pairing}, there is $A$-pairing in any bounded path $[Xo^nY]$ except if $X=Y=B$ and $n$ is odd.

Assume first that $C$ is a fork. Let $c$ be the center of $C$, $a_1$ the neighbor of $c$ on the path to $v_0$ and $a_2$ the vertex between $a_1$ and $v_0$. Let $b_1$ be another neighbor of $c$ and $b_2$ the other neighbor of $y$. We prove the result by induction on the number of branches in the fork. Note that, by Lemma~\ref{lem:sommetdomine}, as $N[a_2]\setminus N[V_t] \subset N[a_1]\setminus N[V_t]$ for $t \in \{A,B\}$, $a_1$ is always a better move than $a_2$, thus we can suppose that it will not be played as next move.

\begin{itemize}
\item if Bob claims  $c$, then Alice replies in $a_1$. Since there is a double $B$-trap and a $A$-pairing, disjoint from the traps, Alice wins by Lemma~\ref{lem:eyes_pairing}. Moreover, by following the pairing strategy, she will keep the fact there is always a $A$-pairing.
 \item If Bob claims $b_1$, then Alice replies in $b_2$. 
 If there are strictly more than three branches, the position is still a fork where the center is dominated by Bob (which is still weakly favorable). The result is true by induction. If there were exactly three branches, then 
 the actual position is $((P_C \cup   [AoA])_{ b_2, b_1},B) = ([Ao^5B]^{\{2, 3\}} \cup   [AoA],B)$. By Lemma~\ref{lem:[A, B]+}, 
 $o((P_C \cup   [AoA])_{ b_2, b_1},B)\geq o( [AoA],B) = \oA$. Furthermore, there is $A$-pairing by Observation~\ref{obs:pairing}.
  \item If Bob claims $b_2$, then Alice replies in $b_1$. If there are strictly more than three branches, then Alice follows the strategy with one less branch. Since she dominates one more vertex, it can only be better for her. If there are exactly three branches, we have as before:
 $$o((P_C \cup   [AoA])_{b_1 b_2},B)  \geq  o( [Ao^5B]^{\{2,3\}} \cup   [AoA],B)  \geq o( [AoA],B) = \oA.$$
 \item If Bob claims $a_1$ , then Alice replies by claiming $c$. We have 
$$o((P_C \cup   [AoA])_{ c, a_1},B) =  o( [Ao^2B] \cup ... \cup [Ao^2B] \cup  [AoA],B)  \geq o( [AoA],B) = \oA.$$ There is a $A$-pairing by pairing together the two vertices belonging to the same paths.
\item If Bob claims  $v_{-1}$, then Alice replies by claiming $c$.
 We have  $$o((P_C \cup   [AoA])_{c, v_{-1}},B) \geq   o( [Ao^2A]^{\{2\}} \cup   [Ao^2B]  \cup  ...\cup [Ao^2B],B)  \geq o( [Ao^2A]^{\{2\}},B) = \oA.$$ Again, there is a $A$-pairing by pairing together the two vertices belonging to the same paths.
\end{itemize}

Assume now that $P_C$ is a bounded path. Alice will never let components of the form $[Bo^{2\ell}B]$ to Bob and thus there will always be a $A$-pairing.
Assume first that $P_C$ has length $9$. We just need to prove that Alice has a strategy in the worst case, that is for $U=\{2,3,5,7\}$. Thus, let $P_C = [Ao^9B]^{\{2,3, 5, 7\}}$. Note that, by Lemma~\ref{lem:sommetdomine}, as $N[v_1]\setminus N[V_t] \subset N[v_2]\setminus N[V_t]$ for $t \in \{A,B\}$, $v_2$ is always a better move than $v_1$, thus we can suppose that it will not be played as next move.

\begin{itemize}
\item If Bob claims $v_1$, then Alice replies in $v_3$. We have 
$$o(([Ao^9B]^{\{2,3, 5, 7\}})_{ v_3, v_1}\cup   [AoA]), B)  \geq  o( [Ao^{6}B]^{\{2,4\}} \cup   [AoA], B)  \geq o([AoA], B) = \oA,  $$
where the last inequality comes from Lemma~\ref{lem:[A, B]+}. 
  
\item If Bob claims $v_3$, then Alice replies in $v_2$. Then $o((P_C\cup [AoA])_{w_2, w_3 }, B)=  \oA$ since there are two $B$-traps and a $A$-pairing.
\item If Bob claims $v_i$,  $4 \leq i \leq 8$,  then Alice replies in $v_{i+1}$.  
We have 
$$o(([Ao^9B]^{\{2,3, 5, 7\}})_{ v_{i+1}, v_i} \cup   [AoA], B)  =  o( [Ao^{i-1} B]^{\{ 2, 3,5,7 \}\cap \{1, ...i-2\}} \cup    [Ao^{9-i-1} B]^{7 -i-1}  \cup [AoA],B)  \geq o( [AoA], B) = \oA,  $$ 
where the last inequality comes from two applications of Lemma~\ref{lem:[A, B]+} (since both  created paths satisfy hypotheses of Lemma~\ref{lem:[A, B]+} ).

%\item if Bob claims in  $v_7$, then Alice claims $v_6$.  At this step,  

%\begin{itemize}
%\item  if now, Bob claims in $v_3$, then Alice replies in $v_2$ And $o(([Ao^9B]^{\{2,3, 5, 7\}})_{ v_{2}, v_3} \cup [AoA], B) = \oA$,  by Lemmas~\ref{lem:eyes_pairing},   since $([Ao^9B]^{\{2,3, 5, 7\}})_{ v_{2}, v_3} \cup [AoA]$ has two traps 
%\item if now, Bob claims in $v_i$   $ i \in \{-1, 1,2, 4, 5\}$, then Bob claims in $v_3$, and then $v_8$ or $v_9$, and, by this way, dominates before Bob. 
%\item if now, Bob claims in $v_i$,    $ i \in \{8, 9\}$, then Bob claims in $v_j$, $ j \in \{8, 9\}$. The resulting position is equivalent to  $[Ao^{5} A]^{\{ 2, 3, \}}   \cup  [AoA]$. and we have 
% $$o ( [Ao^{5} A]^{\{ 2, 3, \}}   \cup  [AoA]) \geq  o( [Ao^{5} B]^{\{ 2, 3, \}}   \cup  [AoA])\geq  o( [AoA]) =\oA$$.
%\end{itemize}

\item If Bob claims $v_9$ , then  Alice replies in $v_8$. We have
$$o(([Ao^9B]^{\{2,3, 5, 7\}})_{ v_{8}, v_9}\cup   [AoA], B) \geq o ( [Ao^{7} A]^{\{ 2, 3, 5, 7 \}}   \cup  [AoA], B)   \geq o ( [Ao^{7}B]^{\{ 2, 3, 5 \}}   \cup  [AoA], B)   \geq o( [AoA], B) = \oA, $$ 
where the last inequality comes from two applications of Lemma~\ref{lem:[A, B]+}.
\item  If Bob claims $v_{-1}$, then Alice replies in $v_3$. We have 
$$o(([Ao^9B]^{\{2,3, 5, 7\}}\cup   [AoA])_{v_3,v_{-1}}, B)  \geq  o( [Ao^{6}B]^{\{2,4\}} \cup   [Ao^2A], B)  \geq o([Ao^2A], B) = \oA,  $$
where the last inequality comes from Lemma~\ref{lem:[A, B]+}.  
\end{itemize}

Assume now that $P_C$ has length 11. As before, we can assume that $P_C = [Ao^{11}B]^{\{2,3, 5, 9\}}$.
\begin{itemize}
\item If Bob claims $v_i$, $1 \leq i  \leq 8$, then it can be done like for a path of length 9.
%\item if Bob claims in  $v_7$, then Alice claims $v_6$.  at this step,  

%\begin{itemize}
%\item  if now, Bob claims in $v_3$, then Alice replies in $v_2$ . And $o(Q_{w_2, w_3 })$ = \Awin,  by Lemmas~\ref{lem:eyes_pairing} and~\ref{lem:C0_pairing} and Remark~\ref{rem:path_pairing}.
%\item if now, Bob claims in $v_i$   $ i \in \{-1, 1,2, 4, 5\}$, then Bob claims in $v_3$, and then $v_8$ or $v_9$, and, by this way, dominates before Bob. 
%\item if now, Bob claims in $v_i$,    $ i \in \{8, 9\}$, then Bob claims in $v_j$, $ j \in \{8, 9\}$. 
    
%We have, $o((((P_C \cup   [AoA]))_{ v_{6}, v_7})_{ v_{8}, v_9} \cup  [AoA]) \geq o ( [Ao^{5} A]^{\{ 2, 3, \}}  \cup [Bo^2B]^{\{2\}}) \geq o(([Ao^9B])_{v_6, v_7}) \cup o( [AoA], B) = \oA$, where the last equality comes from the previous case.  
%\end{itemize}

\item If Bob claims $v_9$, then  Alice replies in $v_8$,  
    \begin{itemize}
    \item If now Bob claims  $v_3$, then Alice replies in $v_2$, creating two $B$-traps. Thus  $o(Q_{v_2, v_3 }, B) = \oA$,  by Lemma~\ref{lem:eyes_pairing}, 
    \item If now Bob claims  $v_5$, then Alice replies in $v_6$, creating two $B$-traps. Thus  $o(Q_{v_6, v_5 }, B) = \oA$,  by Lemma~\ref{lem:eyes_pairing}
    \item If now Bob claims  $v_{10}$  (or $v_{11}$) , then Alice replies in $v_{11}$ (or $v_{11}$).  The resulting position is equivalent to  $[Ao^{7} A]^{\{ 2, 3,5 \}} \cup [AoA]$, and, from Lemma~\ref{lem:[A, B]+},  $$o([Ao^{7} A]^{\{ 2, 3,5 \}} \cup [AoA], B)   \geq o( [AoA], B) = \oA.  $$ . 
    \item If now Bob claims  $v_{i}$ , with $i \in \{-1, 1, 2, 4\} $,  then Alice replies in $v_{3}$. Afterward  Alice can dominate with two more claims,  one in $\{ v_5, v_6\}$, one in $\{ v_{10}, v_{11}\}$, whatever the strategy of Bob;   therefore Alice dominates before Bob
    \item If now Bob claims  $v_{i}$ , with $i \in \{6, 7\} $,  then Alice replies in $v_{5}$. Afterward  Alice can dominate with two more claims,  one in $\{ v_2, v_3 \}$, one in $\{ v_{10}, v_{11}\}$, whatever the strategy of Bob;   therefore Alice dominates before Bob. 
    \end{itemize}
\item if Bob claims in  $v_{10}$, then  Alice replies in $v_{11}$, which reduces the problem to the previous case with 9 vertices.  

\item If Bob claims $v_{11}$, then  Alice replies in $v_{10}$.  The resulting position is equivalent to  to  $[Ao^{9} A]^{\{ 2, 3,5, 9 \}} \cup [AoA]$,  and  by Observation~\ref{obs:paths}, we have 
$$o([Ao^{9} A]^{\{ 2, 3,5, 9 \}} \cup [AoA], B) \geq o([Ao^{9} B]^{\{ 2, 3,5 \}} \cup [AoA], B) = \oA, $$
as seen before. 

\item If Bob claims $v_{-1}$, then Alice replies in $v_3$. We have 
$$o(([Ao^{11}B]^{\{2,3, 5, 9\}}\cup   [AoA])_{v_3,v_{-1}}, B)  \geq  o( [Ao^{8}B]^{\{2,6\}} \cup   [Ao^2A], B)  \geq o([Ao^2A], B) = \oA,  $$
where the last inequality comes from Lemma~\ref{lem:[A, B]+}.  \qedhere
\end{itemize}
\end{proof}

%%%%%%%%%%%%%%%%%%
% 		PARTIE :      C 0
%%%%%%%%%%%%%%%%%%%
\subsubsection{A stable class of positions} \label{sec:stable}
%{\color{red} We will now introduce three classes $\mathcal{C}_0$, $\mathcal{C}_1$, and $\mathcal{C}_1$ of positions which helps us to prove the fact that $o(Q_C) \oA$,  which is the desired result.  

%\begin{itemize}
    %\item We first,  by construction,  that for each $Q \in \mathcal{C}_0$, we have $o([AoA] \cup Q, B) = \oA$.
    %\item Afterwards,  we prove, in Lemma~\ref{lem:C_1} that, for any position  $Q \in \mathcal{C}_1$,  with Bob's turn, Alice has a strategy such that,  after some claims, a  position $Q'$ is reached such that either $o(Q', B) \oA$,  or  $Q' = [AoA] \cup Q''$, with  $Q'' \in \mathcal{C}_0$ (which also implies that $o(Q', B) \oA$).  This proves that $o(Q, B) \oA$. 
    %\item Finally, we prove, in Lemmas~\ref{lem:C_2_standard} and~\ref{lem:C_2_trap} that, for any instance $Q$ of  $\mathcal{C}_1$ with Bob's turn, Alice has a strategy such that,  after some claims, a  position $Q'$ is reached such that either $o(Q', B) \oA$,  or  $Q'   ±in \mathcal{C}_1$, (which also implies that $o(Q', B) \oA$). This proves that $o(Q, B) \oA$.
%\end{itemize}
%Since $Q_C \in \mathcal{C}_1$, we get the result. 
%}

In order to define our stable class $\mathcal C_1$, we first define the class $\mathcal C_0$, that informally corresponds to the positions derived from a weakly favorable position when Alice follows a winning strategy.

\begin{definition}[$\mathcal{C}_0$]
The class $\mathcal{C}_0$ is defined recursively as follows. A position  $P$ is an element of $\mathcal{C}_0$ if:
\begin{itemize} 
\item $P = P_C$ where $C$ is a weakly favorable component;
\item there exists  $P' \in \mathcal C_0$, with two unclaimed vertices $x$ and $y$, such  that $P = P'_{x, y}$  and $o(P \cup [A oA], B) = \oA$.  
%an  instance of  $\mathcal{C}_0$ is an instance   $Q= (G_C, V_a \cup \{x_1, x_2..., x_p\},  V_b \cup \{y_1, y_2..., y_p\}, B) $ such that for each $i$, $0 \leq i \leq p$, 
%$o((G_C , V_a \cup \{x_1, x_2..., x_{i}\},  V_b \cup \{y_1, y_2..., y_i\}, B) \cup  [AoA] )$ =  \Awin. 
\end{itemize}
\end{definition}

Next corollary is a direct application of Lemma~\ref{lem:weakly}.
\begin{corollary}
    For each $P \in  \mathcal{C}_0$,  we have  $o(P \cup [A oA], B) = \oA$.
\end{corollary}

%\begin{proof}
%First assume that Assume first that  the set of free vertices of $P_C$  form a path. 
%Assume first that  the set of unclaimed vertices of $P_C$  form a path  $(v_1, ...,v_n)$.  since $o([AoA] \cup  Q) = \oA$, thus  there is no subsequence  $(v_i, v_{i+1}, ..., v_j)$ , with $j-i$ even  such that both $v_{i-1}$ and $v_{j+1}$ are in the set of elements previously claimed by Bob,  by Lemma~\ref{lem:[BB]_odd}.  Thus,  there exists an A-pairing  of $Q$ by Observation~\ref{obs:pairing}. 

%Otherwise,    the set of unclaimed vertices of $P_C$  form a fork, we also obtain that there exists an  A-pairing  of $Q$, by an easy case by case analysis. 
%\todo[inline]{Eric R:  A préciser?? A mettre dans le lemme 6.11 ?}
%\end{proof}

%%%%%%%%%%%%%%%%%%%%%%%%%%%%%%%%%%%%%%

%        PARTIE :      C 2

%%%%%%%%%%%%%%%%%%%%%%%%%%%%%%%%%%%%%%

We can now define the stable class of positions $\mathcal C_1$
 
\begin{definition}[$\mathcal C_1$]\label{def:c1}
 A position  $P$ belongs to the class  $\mathcal{C}_1$ if $P$  is the union of the following positions:
 \begin{enumerate}
 %\item one position $[A o A]$, whose  unique unclaimed vertex is denoted by  $v_{-1}$,  
\item  at most one position in  $\mathcal{C}_0$. 
\item  a union of positions of the type   $[Ao^nB]^{\{2,5, n-2\}}$  $n  \geq 9$, $n$ odd. 
\item  a union of $k$  positions of the type   $ [B o^{2t} B]^{\{{2t -2}\} }$, $t \geq 1$,
\item   a union of $k'$  positions of the type   $ [A o^{n} A]^{\{{2, 5}\} }$, with $k' \geq k$ and  $n\geq 1$.
\end{enumerate}
 \end{definition}

\begin{remark}\label{rem:C1}
By definition, each position formed by at most one weakly favorable position and some  positions of the alternative 1.d of Definition~\ref{def:fav} is an element of $\mathcal{C}_1$.
\end{remark}

%\begin{corollary}
%If $Q \in \mathcal{C}_0$ and $Q'$ which admits a A-pairing. 
%Then Alice has a strategy to win the Maker-Breaker Dominating Game on $Q \cup Q'$, even when  Bob starts.
%\end{corollary}

Next corollary is a direct application of Lemma~\ref{lem:weakly} and Observation~\ref{obs:pairing}.
\begin{corollary}\label{cor:pairing}
    For each $P \in  \mathcal{C}_1$, $P$ admits an $A$-pairing.
\end{corollary}

%The idea for the class $\mathcal{C}_1$ is to encode some positions which can happen starting from  a position $P \cup [AoA], B)$ when $P$ is  formed by at most a weakly favorable position and some  positions of the alternative 1.d of   Theorem~\ref{th:maintree},   Bob does not claim $v_{-1}$ (the unclaimed vertex of $[A o A]$) and  Alice plays a correct strategy. 

We will now prove that the class $\mathcal{C}_1$ (with a trap $[AoA]$ adjoined) is either stable after Alice's answer to Bob's claim, or directly winning for Alice. We will consider the two cases according to whether Bob claims the unclaimed vertex $v_{-1}$ of $[AoA]$ or not. The first case where Bob does not claim $v_{-1}$ is proved by Lemma~\ref{lem:C_2_standard} that requires Claim~\ref{lem:[B, B]_pair}. The second case is when Bob claims $v_{-1}$ and is solved by Lemma~\ref{lem:C_2_trap} that requires Claim~\ref{lem:[AooAooA]bis} as a particular case. 

\begin{claim} \label{lem:[B, B]_pair}
Let  $Q $ be any position, and $t$ be  a positive  integer. Let $(w_1, w_2, ..., w_{2t})$ be the unclaimed  vertices of  $[B o^{2t} B]^{\{{2t -2}\}} $. 
For each $w_i$, there exists $w_j\neq w_i$  and $0 \leq t' < t$ such that 
$$ o(Q \cup [B o^{2t'} B]^{\{{2t' -2}\} }, B)  \leq  o(Q    \cup [B o^{2t} B]^{\{{2t -2}\} }_{w_j, w_i}, B )$$ 
(with the convention that  $[Bo^0B]^{\{-2\}}$ is empty). 
\end{claim}

The idea behind this claim is that if a position contains a bounded path $[B o^{2t} B]^{\{{2t -2}\}} $ and Bob claims a vertex in it, then Alice can reply in the same bounded path, preserving the global structure of the  position.
%{\color{blue}
\begin{proof} We have three cases. 
\begin{itemize}
\item  If $i = 2t$ then Alice claims $w_{2t-1}$, using Observations~\ref{obs:paths} and~\ref{obs:include}, we have 
 $$o(Q \cup ([B o^{2t} B]^{\{{2t -2}\} })_{w_{2t-1}, w_{2t}},B) = o(Q \cup [B o^{2t-2} A]^{\{2t-2\} }, B) \geq  o(Q \cup [B o^{2t-2} B], B) \geq  o(Q \cup [B o^{2t-2} B]^{\{2t-4\} }, B)$$
 %\todo[inline]{ER ? Finir la preuve et la mettre avec le cas pair
 %ER : J'ai fait trois cas plutôt q'un cas et deux sous-cas. Pour moi la preuve est complète.  }
\item If $i = 2k$, $1 \leq k < t$,  
then Alice claims $w_{2k-1}$. We have 
$$o(Q \cup [B o^{2t} B]^{\{{2t -2}\} }_{w_{2k-1}, w_{2k}}, B) = o(Q \cup [Bo^{2k-2}A] \cup    [Bo^{2(t-k )}B]  ^{\{{2(t-k) -2}\}}, B) .$$ 

Lemma~\ref{lem:[A, B]+} applies,  and we get 
$$o(Q \cup [B o^{2t} B]^{\{{2t -2}\} }_{w_j, w_i}, B) \geq o(Q \cup    [Bo^{2(t-k )}B]  ^{\{{2(t-k) -2}\}}, B). $$
%\todo[inline]{Eric D: le cas $k=t$ ne me parait pas géré. ça donne une union de Q avec un $[BO^{2t-2}A]$ perturbé en $(2t-2)$  

%Eric R, j'ai proposé une alternative sans Lemma 2. 3, qui me semble bien complexe. }
%Note that we do not consider the case $k=t$ since from Lemma~\ref{lem:sommetdomine}, the vertex $w_{2t}$ is dominated by $w_{2t-1}$. 

\item  If $i = {2k-1}$, $1 \leq k \leq t$, then Alice claims $w_{2k}$. We have 
$$o(Q \cup [B o^{2t} B]^{\{{2t -2}\} }_{w_j, w_i}, B) = o((Q  \cup [Bo^{2k-2}B] \cup    [Ao^{2(t-k) }B]  ^{\{{2(t-k) -2}\}}, B), $$ 

Lemma~\ref{lem:[A, B]+} applies,  and we get 
$$o(Q \cup (B o^{2t} B]^{\{{2t -2}\} }_{w_j, w_i}, B) \geq o(Q \cup    [Bo^{2k-2}B] , B)\geq  o(Q \cup    [Bo^{2k-2}B]^{\{2k-4\}} , B). $$ \qedhere
\end{itemize}
\end{proof}

\begin{lemma}[Bob does not claim $v_{-1}$] \label{lem:C_2_standard}
Let $P  \in \mathcal{C}_1$,  such that  $P$ is not dominated by Alice. Let $Q = P \cup [AoA]$ and $y$ be an  unclaimed vertex of $P$. 
There exists an  unclaimed  vertex  $x \neq y$ such that at least one of the following alternatives holds: 
\begin{itemize}
\item $o(Q_{x, y}, B ) = \oA$;
\item  or there exists  $P' \in \mathcal{C}_1$  with at least two unclaimed vertices less than $P$, such that $o(P'\cup[AoA], B) \leq  o(Q_{x, y}, B)$.
%\item  there exists  $Q' \in \mathcal{C}_1$ such that $o(Q') \leq  o(Q_{x, y})$
%\item  there exists  $P' = (G, C_a, C_b, B) \in C_{next}$  with at least two  free vertices than $P$, such that $o(P') \leq  o(G, C_a \cup \{y \}, C_b \cup \{x \}, B)$.
\end{itemize}
 \end{lemma}
 
\begin{proof} 
We denote by $v_{-1}$ the unclaimed element of $[AoA]$ in $Q$. We assume that Bob claims $y$.
We have several cases, according to the component that contains $y$.
\begin{itemize}
    \item Assume first that $y$ is an element of the component $P'$ of $P$ that belongs to $\mathcal{C}_0$.
If there exists  a free vertex $x$ of  $P'$ such that $P'_{x, y} \in \mathcal{C}_0$, then we are done, since  $P_{x, y} \in \mathcal{C}_1$. 
	
Otherwise, by definition of $\mathcal C_0$, $(P'\cup [AoA], B)$ is \Awin. Since there is no winning answer in $P'$ for Alice to the claim $y$ of Bob, we necessarily have $o((P'\cup [AoA])_{v_{–1}, y}, B)  = \oA$. From the position $(P'\cup [AoA])_{v_{–1}, y}$, Bob cannot dominate the whole graph. Alice can follow her strategy on $(P'\cup [AoA])_{v_{–1}, y}$ and a pairing strategy on the rest of the graph (that exists by Corollary~\ref{cor:pairing}).
    This ensures that $o(Q_{v_{–1}, y}, B) = \oA$

\item Assume now that $y$ is an unclaimed vertex of $[Ao^nB]^{\{2,5, n-2\}}$, with $n$ odd, $n\geq 9$. Let $(w_1, w_2, ...,  w_n)$ be the sequence of free vertices of   $[Ao^nB]^{\{2,5, n-2\}}$. By Lemma~\ref{lem:sommetdomine}, as $N[w_1]\setminus N[V_t] \subset N[w_2]\setminus N[V_t]$ for $t \in \{A,B\}$, one can assume that $y\neq w_1$.
Let $P'$ be the position such that $P=P'\cup [Ao^nB]^{\{2,5, n-2\}}$ and let $Q'=P'\cup [AoA]$. Note that $P'\in \mathcal{C}_1$.
	\begin{itemize} 
	\item If $y=w_2 $ then take $x = w_3$. By Lemma~\ref{lem:split}, Observation~\ref{obs:removedomcomp} and Lemma~\ref{lem:[A, B]+},  $$o(Q_{x, y}, B) \geq   
		o(Q'\cup [Ao^{n-3}B]^{\{2, n-5\}}, B) \geq   o(Q', B).$$
	\item If $y = w_3 $,  then take $x= w_2$. There is double $B$-trap in $Q_{w_2, w_3 }$. Thus, by Lemmas~\ref{lem:eyes_pairing} and~\ref{cor:pairing}, we have  $o(Q_{w_2, w_3 },B)= \oA$.
	\item If $y   = w_i$  $i \geq 4$ with $i$ even,  then take $x = w_{i+1}$	
	 We have  $$o(Q_{x, y}, B) =   
		o(Q'\cup [Ao^{i-1}B]^{\{2, 5\}} \cup [Ao^{n-i}B]^{\{n-i-2\}}, B) \geq   o(Q'\cup [Ao^{i-1}B]^{\{2, 5\}}, B) $$ by  application of Lemma~\ref{lem:[A, B]+}. Note that $P'\cup[Ao^{i-1}B]^{\{2, 5\}}$  is an element of $\mathcal{C}_1$. 
		
		%For $i \geq 9$, $(Q \setminus [Ao^nB]^{\{2,5, n-2\}}) \cup [Ao^{i-1}B]^{\{2, 5\}} \in \mathcal{C}_1$, and for $i \in \{ 7, 8 \}$, Lemma~\ref{lem:[A, B]+} again applies to get that $o(Q_{x, y}) \geq  o(Q \setminus [Ao^nB]^{\{2,5, n-2\}})$. 
		
		\item if $y   = w_i$ with $i$ odd  and $i \geq 5$   then take $x = w_{i-1}$.
		We have    $$o(Q_{x, y}, B) =  
		o(Q'\cup [Ao^{i-2}A]^{\{2, 5\}}\cup [Bo^{n-i}B]^{\{n-i-2\}}, B). $$	Note that $P'\cup [Ao^{i-2}A]^{\{2, 5\}}\cup [Bo^{n-i}B]^{\{n-i-2\}}$ is an element of $\mathcal{C}_1$, since $n-i$ is even and since we add a position $[Ao^{i-2}A]^{\{2, 5\}}$.
	\end{itemize}

\item If $y$ is an unclaimed  vertex of a component of the form $[B o^{2t} B]^{\{{2t -2}\} }$, then the position can be reduced by Claim~\ref{lem:[B, B]_pair}. 

\item  Assume now that $y$ is an unclaimed vertex of  $[Ao^nA]^{\{2,5\}}$. As before, we denote by $(w_1, w_2, ...,  w_n)$ the vertices of $[Ao^nA]^{\{2,5\}}$, by $P'$ the position $P$ without the component $[Ao^nA]^{\{2,5\}}$ and by $Q'$ the position $P'\cup [AoA]$.
If $n \leq 3$, then  take $x = v_{-1}$. 
 We have $o(Q_{ v_{-1}, y},B) = \oA$ since
 the position   $Q_{ v_{-1}, y}$ admits a  $A$-pairing  by Corollary~\ref{cor:pairing} and Bob.  
 Thus, it can be assumed that $n\geq 4$.
		\begin{itemize}
		  \item If $y \in \{w_1, w_2\} $, then take $x = w_{3}$. 
            We have $$o(Q_{w_{3}, y},B) \geq   
		  o(Q' \cup [Ao^{n-3}A]^{\{2\}},B).  $$ Note that $P'\cup [Ao^{n-3}A]^{\{2\}}$ is in $\mathcal C_1$.
	
		  \item If $y = w_3 $,  then take $x= w_2$.  The position $Q_{w_2, w_3 }$  admits an $A$-pairing, by Corollary~\ref{cor:pairing}. Thus, by lemma~\ref{lem:eyes_pairing}, we have   $o(Q_{w_2, w_3 }) = \oA$. 
		  %\item If $y =  w_4 $,
    %we distinguish the case depending on the value of $n$.
	%\begin{itemize}
  % \item If $n  =  4$, take  $x= w_ {3}$. We have $o(Q_{x, y}, B) \geq   o(Q'\cup [Ao^{2}A]^{\{2\}},B)$. The position $P'\cup [Ao^{2}A]^{\{2\}}$ is still in $\mathcal C_1$.
   %\todo[inline]{A:J'ai l'impression qu'il faut laisse assez de place dans le $[Ao^{n-5}A]$, du coup faudrait aussi traiter $n=5$ ? peut être fait avec le cas d'avant ?}
			%\item if $n \geq 5$, then take 
			%$x= w_ {5}$. We have 
			%$$o(Q_{x, y}, B) =   o(Q'\cup [Ao^{3}B]^{\{2\}} \cup [Ao^{n-5}A],B)  \geq  o(Q'\cup [Ao^{n-5}A],B)  \geq  o(Q'\cup [Ao^{n-5}A]^{\{2,5\}},B)$$, according to Lemma~\ref{lem:[A, B]+}.   The position $P'\cup Ao^{n-5}A]$ is still in $\mathcal C_1$. \todo{A:On a besoin de $n\geq 6$ la non ??}
			%\end{itemize}
	 	\item if $y  =  w_i$, with  $i  \geq 4$,   
         then take $x= w_ {i-1}$. We have
		$$o(Q_{x, y}, B) =   o(Q'\cup [Ao^{i-2}A]^{\{2, 5\}}  \cup [Ao^{n-i}B],B)   \geq   o(Q'\cup [Ao^{i-2}A]^{\{2, 5\}}, B) $$   according to Lemma~\ref{lem:[A, B]+}. We are done since $P'\cup [Ao^{i-2}A]^{\{2, 5\}}$ is in $\mathcal C_1$. \qedhere
	 	\end{itemize}
\end{itemize}
\end{proof}

\begin{claim} \label{lem:[AooAooA]bis}
Let  $P$ be any position, and $t$ be  a nonnegative   integer. 
 We have 
$$   o(P \cup [AoA], B)  \leq  o(P \cup   [AooAooA]^{\{2, 5\}}  \cup [B o^{2t} B]^{\{{2t -2}\} }, B )  $$ 
where  $[AooAooA]^{\{2, 5\}}$ is the position obtained from $[Ao^5A]^{\{2, 5\}}$ by adding the central vertex $v_3$ in the set $V_A$. 
\end{claim}

\begin{proof}
Let $(w_1, w_2, ..., w_{2t})$ be the sequence of unclaimed  vertices of  $[B o^{2t} B]^{\{2t -2\}}$, $(v_1, v_2, v_4, v_5)$  the sequence 
of unclaimed  vertices of $[AooAooA]^{\{2, 5\}}$, from left to right, and $v_{-1}$ denote 
the unclaimed vertex of $[AoA]$. To lighten notation,  we also state  
$ R = P \cup   [AooAooA]^{\{2, 5\}}  \cup [B o^{2t} B]^{\{{2t -2}\} }$

Assume that  $o(P  \cup [AoA], B) = \oA$. We prove that   $o(R, B ) = \oA$  by induction on the number $p$ of unclaimed vertices of $R$. For initialization, if $p  \in   \{4, 5\}$, then $t = 0$ and all the vertices of $P$ are claimed except possibly one. Thus, since  $o(P \cup [AoA], B) = \oA$, Alice dominates $P$, and therefore Alice dominates  $R=P \cup [AooAooA]^{\{2, 5\}}$. 

 Now assume that $p \geq 6$. We have several alternatives according to the claim $y$ of Bob in $R$.
 
\begin{itemize}
    \item  If $y$ is an unclaimed vertex of $[B o^{2t} B]^{\{{2t -2}\}}$, then Claim~\ref{lem:[B, B]_pair} applies, and the induction gives the result. 
    \item Assume now that $y$ is an unclaimed vertex of $P$. If there exists an unclaimed vertex $x$ of $P$ such that $o(P_{x,y} \cup [AoA], B) = \oA$.  We conclude by the induction hypothesis.  If it is not possible, it means that once Bob has claimed $y$ in the position $P \cup [AoA]$, Alice is forced to claim $v_{-1}$.
    Then, in $R$, Alice claims $v_2$, creating a double $B$-trap in $v_1$ and $v_4$. 
    
    In the position $((P \cup [AoA])_{v_{-1},y}, B)$, Alice has a strategy
        which allows her to dominate $P$.  Thus, Alice can  win in  $R_{v_2,y}$ playing component by component,  until,  the game ends  as follows:  
            \begin{itemize}
            \item If Bob claims in $P$, Alice also claims in $P$ according to the strategy in  $((P\cup [AoA])_{v_{-1},y}, B)$.
            \item If Bob claims  in    $[AooAooA]^{\{2, 5\}}$,  then Alice isolates one vertex using one of the $B$-trap. 
            \item  If Bob claims in $[B o^{2t} B]^{\{{2t -2}\} }$, Alice follows a pairing strategy with the $A$-pairing of this component.
            \end{itemize}
       This way, Alice will prevent Bob to dominate since Bob has to claim at some point in $[AooAooA]$. The two other items ensure that Alice while dominate the whole graph.
    \item Assume finally that $y$ is  an unclaimed vertex of $[AooAooA]$. The position $(P, A )$ is \Awin since it corresponds to the position $(P\cup [AoA],B)$ where Bob has claimed $v_{-1}$. Let $x$ be a unclaimed vertex of $P$ such that $o(P_x, B) = \oA$, where $P_x$ is the position obtained from $P$ when Alice has claimed $x$. Alice answers $x$ in the game $R$. We prove that $o(R_{x,y},B)=\oA$. Alice plays as follows:
    \begin{itemize}
    \item If Bob claims a vertex in $P_x$, Alice answers in $P_x$ according to her strategy in $(P_x, B)$. 
    \item If Bob  claims a vertex in $[B o^{2t} B]^{\{2t -2\}}$, then Claim~\ref{lem:[B, B]_pair}]  can be used.
    \item If Bob claims  again a vertex in $[AooAooA]$, Alice claims $w_1$.
    \end{itemize}
    The first item ensures that Alice dominates $P$  before Bob.  
    The second item ensures that the fact that Bob claims  in the component $[B o^{2t} B]^{\{2t -2\}}$ is irrelevant. 
    In the case of third item,  one can consider using Lemma~\ref{lem:[A, B]+} that the component $[B o^{2t} B]^{\{2t -2\}}$ disappears.   The component $[AooAooA]$ becomes also irrelevant  since   each player dominates all vertices of this component.  Thus, it remains only the component derived from $P_x$ where Alice dominates before Bob. \qedhere
\end{itemize}
      \end{proof}

\begin{lemma} [Bob claims $v_{-1}$]  \label{lem:C_2_trap}

Let $P \in \mathcal{C}_1 $  such that  $P$ is not dominated by Alice and $Q = P \cup [AoA] $.

There exists an unclaimed vertex $x$ of $P$ such that  at least one of the following alternatives holds: 

\begin{itemize}
\item $o(Q_{x, v_{-1}}, B) = \oA$, 

\item  there exists  $P'  \in \mathcal{C}_1$  with at least one unclaimed  vertex less than $P$, 
such that $o(P'\cup[AoA], B) \leq  o(Q_{x, v_{-1}}, B)$. 
\end{itemize}

 \end{lemma}

\begin{proof} We have different cases according to the structure of $P$. 

\begin{itemize}
\item If $P$ contains a component $[Ao^nB]^{\{2,5, n-2\}}$, $n  \geq 9$, $n$ odd, let $v_1, v_2, ..., v_n$ denote the sequence of free vertices of this component. 
Take $x = v_2$.  
Let $P'$ be the position such that $P=P'\cup [Ao^nB]^{\{2,5, n-2\}}$ and let $Q'=P'\cup  [AoA]$. Note that $P'\in\mathcal{C}_1$.

We now have
$$o(Q_{v_2, v_{-1}}, B) \geq o(Q' \cup  [Ao^{n - 2}B]^{\{3, n-4\}},B)\geq o(Q', B )  $$ from Lemma~\ref{lem:[A, B]+}.

 \item If $P$ contains a component $ [A o^{n} A]^{\{{2, 5}\} }$, let $w_1, w_2, ..., w_n$ denote the sequence of  vertices of $ [A o^{n} A]^{\{{2, 5}\} }$. 
 Let $P'$ be the position such that $P=P'\cup [A o^{n} A]^{\{{2, 5}\} }$ and let $Q'=P'\cup  [AoA]$. 
 
 \begin{itemize}
\item  If $n \geq 6$, then take  $x = w_{n-1}$. 
 We have

$$o(Q_{w_{n-1}, v_{-1}}, B) \geq o(Q'  \cup   [A o^{n-2} A]^{\{{2, 5}\} }, B )$$ and note that $P'\cup [A o^{n-2} A]^{\{{2, 5}\} }$ is an element of $\mathcal{C}_1$.

%\geq o(Q \setminus   [A o^{p} A]^{\{{2, 5}\} }  )$, from Lemma~\ref{lem:[A, B]+}. 

 \item If 
 $3 \leq   n \leq 4$, then take  $x = w_{2}$. 
 We have

$$o(Q_{w_{2}, v_{-1}}, B) \geq o(Q' \cup   [A o^{n-2} A] , B  )$$ and we conclude as previously.

\item If $n= 1$, then take $x = w_1$.  The position $Q_{w_1, v_{-1}}$ admits an $A$-pairing from Corollary~\ref{cor:pairing}. Thus, from~\ref{lem:trap_pairing}, we have $o(Q_{w_1, v_{-1}}, B)= \oA$.

\item If $p= 5$ and  $k \geq 1$, take $x = w_3$, which  creates an $[AooAooA]$. We define $P''$ such that  
$$P'=P''\cup   ([B o^{2t} B]^{\{{2t -2}\} }).$$ 
Note that $P'' \in \mathcal{C}_1$ as one component of each type 3 and 4 (in Definition~\ref{def:c1})) has been removed. Lemma~\ref{lem:[AooAooA]bis} applies, thus  
$$o(P'' \cup [AoA], B)\leq o(P'' \cup [AooAooA]^{\{{2, 5}\} } \cup [B o^{2t} B]^{\{{2t -2}\} }, B)= o(Q_{w_3,v_{—1}}, B)$$
    \end{itemize}

\item If $P$  cannot be treated in one of the previous cases, then $P$ only contains some subpositions of the type    $[A o^{5} A]^{\{{2, 5}\} } $ and at most one 
position $P_0 \in  \mathcal{C}_0$.  In this case,  we will prove that $o(Q, B) = \oA$, which implies the result. 

First we  have $o([A o^{5} A]^{\{2, 5\}}, B)  =  \oA$.
Indeed, Bob is forced to claim $w_3$, since otherwise Alice claims $w_3$ and wins. After this claim, there exists an $A$-pairing of size 2, while Bob needs at least two more claims to dominate the graph. 
Thus $Q = P \cup [AoA]$ is formed by a union  of winning positions (when Bob starts):  one is $P_0  \cup [AoA] $ (or simply  $[AoA]$ if $P_0$ does not exist) and the other ones are positions
$[A o^{5} A]^{\{{2, 5}\}}$.  Thus,  by Observation~\ref{obs:union}, we get $o(Q, B)= \oA$.  \qedhere
%\item if no item above is saisfied, then there is  no $[Ao^nB]^{\{2,5, n-2\}}$, and all   $ [A o^{p} A]^{\{{2, 5}\} }$ are such that $p = 5$. 
%\begin{itemize}
%\item if there is such a $ [A o^{5} A]^{\{{2, 5}\} }$ as component of $Q$, then take $x = w_3$. By this way,  $Q_{x, v_{-1}} \in \mathcal{C}_1$. 
%\item otherwise,  i.e if no case before is satisfied,   $Q \setminus [A o A]  \in \mathcal{C}_0 $, thus,  by definition,  there exists $x$ such that $o(Q_{x, v_{-1}}, B)= \oA$. 
%\end{itemize}
\end{itemize}
\end{proof} 
 
\begin{corollary}\label{cor:final}
For each position $P \in \mathcal{C}_1$, we have $o(P \cup [AoA], B)= \oA$.
\end{corollary}

\begin{proof}We prove the result by induction on the number of unclaimed vertices in $Q=P\cup[AoA]$.
If Alice dominates $P$, then we directly have $o(P\cup[AoA],B)=\oA$.
Otherwise, consider a move $x$ of Bob.
If $x\neq v_{-1}$ (respectively $x=v_{-1}$), then by Lemma~\ref{lem:C_2_standard} (resp. Lemma~\ref{lem:C_2_trap}), there exists an unclaimed vertex $y$ of $P$ such that either Alice wins or there exists a position $P'$ of $\mathcal C_1$ with less unclaimed vertices such that $o(P'\cup[AoA],B)\leq o(Q_{x,y},B)$. By induction, $o(P'\cup[AoA],B)=\oA$ and thus $o(Q_{x,y},B)=\oA$. 

Thus, for any claim $x$ of Bob, Alice can claim a vertex $y$ such that $o(Q_{x,y},B)=\oA$. Therefore, $o(Q,B)=\oA$.
\end{proof}

\subsubsection{Conclusion}

Putting together all the previous results, we can prove the reverse part of Theorem~\ref{th:maintree}:

\begin{corollary}
If all the components of $\sk_F$ are favorable to Alice and at most one of them is weakly favorable, then $((F,\{v_0\},\emptyset),B)$ is \Awin. 
\end{corollary}

\begin{proof}
    Let $Q=(F,\{v_0\},\emptyset)$ be a position such that all the components of $\sk_F$ are favorable to Alice and at most one of them is weakly favorable. By Lemma~\ref{lem:[A, B]+}, one can assume that there is no strongly favorable component in $\sk_F$.  Let $P$ such that $Q=P \cup [AoA]$. By Remark~\ref{rem:C1}, $P$ is an element of $\mathcal{C}_1$. Corollary~\ref{cor:final} leads to the desired result.
\end{proof}

\bibliographystyle{plain}
\bibliography{Journal}

\end{document}